\documentclass[12pt]{article}

\usepackage{amsmath,amssymb, amscd}

\def\CC{{\mathbb C}}

\def\RR{{\mathbb R}}
\def\QQ{{\mathbb Q}}

\def\NN{{\mathbb N}}
\def\PP{{\mathbb P}}

\def\bd{{\partial}}
\def\del{{\partial}}

\newcommand {\eps} {\epsilon}
\newcommand {\ga} {\gamma}
\newcommand {\te} {\tilde}
\newcommand {\yd} {{Y^\partial}}
\newcommand {\xd} {{X^\partial}}
\newcommand {\fd} {{f^\partial}}
\newcommand {\gd} {{g^\partial}}
\newcommand {\hd} {{h^\partial}}
\newcommand {\ks} {{\kappa_\sigma}}
\newcommand {\BOX} {\rule{2mm}{2mm}}

\newcommand {\m} {{\setminus}}
\newcommand {\dv} {{\rm Div}}

\newcommand {\be} {\begin{equation}}
\newcommand {\eeqn} {\end{equation}}
\newcommand {\bea} { \begin{eqnarray}}
\newcommand {\eea} {\end{eqnarray}}
\newcommand {\beas} { \begin{eqnarray*}}
\newcommand {\eeas} {\end{eqnarray*}}
\newcommand {\ra} {\rightarrow}
\newcommand {\dra} {\dashrightarrow}
\newcommand {\hra} {\hookrightarrow}
\newcommand {\lra} {\longrightarrow}

\newtheorem {lemma} {LEMMA} [section]
\newtheorem {theorem}[lemma]{THEOREM}
\newtheorem {prop}[lemma]{PROPOSITION}
\newtheorem {lem}[lemma]{LEMMA}
\newtheorem {cor}[lemma]{COROLLARY}
\newtheorem {defn}[lemma]{DEFINITION}
\newtheorem {conj}[lemma]{CONJECTURE}

\newtheorem {rem}[lemma]{Remark}
\newtheorem {ex}[lemma]{Example}
\newtheorem {definition}[lemma]{DEFINITION}

\numberwithin{equation}{section}

\begin{document}


\title{A refined Kodaira dimension and its canonical fibration\\
\large{\it Dedicated to Professor Hirzebruch on his 75th Birthday}}
\author{Steven S. Y. Lu
\thanks{Partially supported by an NSERC grant and the Max Planck
Institute of Mathematics in Bonn.}}
\date{}
\maketitle


\section{Introduction}

{} From many earlier works on the classification theory of 
compact complex manifolds and their intrinsic geometric
structures, it has been
clear that ``positivity'' or ``nonpositivity''
properties of subsheaves of exterior
powers  of the tangent or
cotangent bundle provide some of the most important global
information concerning the manifold in question. The
programs of Iitaka and Mori on a general
classification scheme for algebraic varieties testify
to this with a special focus on the top exterior
power of the cotangent bundle, the canonical bundle.  
This paper is motivated by
some  of our earlier studies to bring
the other subsheaves of the cotangent bundle, such as
those defined by foliations and fibrations
(studied by F. Bogomolov and Y. Miyaoka for example), 
into focus as important objects of study for 
a general classification theory in birational geometry.
Especially relevant here are the line subsheaves of 
exterior powers of the cotangent bundle that correspond
to the canonical ``bundles'' of the orbifold bases
of fibrations.
\\

Our main motivation comes from the papers \cite{BL1,BL2}
and \cite{Lu01}. In these papers, it was shown that the
properties of having identically vanishing Kobayashi 
pseudometric (or holomorphic connectedness), 
having a Zariski dense entire holomorphic curve
and being dominated  by a holomorphic map from
$\CC^n$ (in the sense of having maximal
rank somewhere) are the same for most algebraic varieties in low
dimensions as well as for rationally connected fibrations and
abelian fibrations over curves and that this property
is characterized by birational invariants. In \cite{BL1},
we studied the compact complex surface case, among others,
and found that the Kodaira dimension and the fundamental
group are sufficient to characterize these properties.
\cite{BL2} showed that this characterization fails
for quasiprojective surfaces but a more intrinsic
birational invariant related to Campana's $\kappa^+$ in
\cite{Ca} was introduced in \cite{Lu01} that
completely characterized the above properties, even for
the quasiprojective case, for the above mentioned examples
and more (namely, those for which any of the above properties
can be verified via the techniques we have introduced in
these papers).\\

Our invariant is defined simply as follows. 
Given a compact complex manifold $X$ with a normal crossing 
{\em boundary} divisor $D$,
we denote the log manifold $(X,D)$ by $\xd$.
We define a birational invariant $\kappa'_+(\xd)$ 
of $X\smallsetminus D$ as the maximum integer $p$ 
for which there is a dominant meromorphic map $f: X\dra Y$ 
to a $p$-dimensional manifold $Y$ 
such that $\kappa (L)=p$ where $L=L_\fd$ is the saturation
of $f^*K_Y$ in $\Omega^p_X(\log D)$. Such line sheaves for
projective $X$ were first studied
by F. Bogomolov in \cite{Bo}. Saturated subsheaves 
corresponding to foliations were also studied by Y. Miyaoka 
in connection with classification theory, see \cite{MP}.\\

We call the orbifold $\xd$ special
if $\kappa'_+(\xd)=0$ and call a fibration from $\xd$ special
if its orbifold general fibers are.  We note that all
examples of compact complex manifolds known to us
that do not have vanishing Kobayashi pseudodistance
are non-special.
If $\dim X\leq 3$ and $X$ is not
birationally equivalent to a Calabi-Yau threefold, we have
verified in \cite{Lu01} that a projective $X$ is special if
and only if it is holomorphically chain connected in the sense 
that  there is a holomorphic family of entire 
holomorphic curves in $X$ such that any general pair of points
in $X$ are connected by a chain of holomorphic curves in this
family. Also any holomorphically connected projective manifold $X$
is special and hence so are the rationally connected ones. 
An important fact about special varieties, 
first obtained by Campana for a slightly different 
definition from above but more recently for the
definition above (\cite{Ca02}, see also \cite{Lu02}), is that
a projective manifold of Kodaira dimension zero is special.
As rationally connected manifolds are also
special, the two basic types of fibrations in the 
classification program of algebraic varieties, namely 
rationally connected fibrations
and Iitaka fibrations
are examples of special fibrations. 
Thus, special objects are basic building blocks in birational
geometry.
We mention in this connection
that Campana, Koll\'ar, Miyaoka and Mori
(see \cite{KoMM, Ca, Kol1, Db}) have defined (and
developed extensively the theory of) the maximal rationally 
connected (meromorphic) fibration (or MRC-fibration 
\cite[IV.5]{Kol1}) of a smooth projective variety.
The structure of this MRC-fibration, which is a proper fibration 
when restricted to some open subset (i.e., almost holomorphic),
has just been resolved completely in
\cite{GHS} using the
theory developed by the above mentioned people
and the theory of compactification of 
the appropriate space of stable curves 
\`a la Kontsevich:

\begin{theorem}[Graber-Harris-Starr]\label{ghs} 
If a smooth projective variety has a rationally 
connected (meromorphic) fibration over a  
rationally connected 
base, then it is rationally connected. Hence,
the base of the MRC-fibration 
of a projective manifold is not uniruled, or  
equivalently, the base admits no 
rationally connected (meromorphic) fibrations with positive
dimensional fibers. Equivalently, every rationally connected
fibration over a curve admits a section with prescribed 
behavior (in the sense given by \cite{BL2} for example)
over any finite number of points of the curve.
\end{theorem}

Campana has
also constructed in \cite{Ca02} (from a viewpoint
that is independent from ours) a canonical 
special fibration for any projective manifold
without boundary divisor 
and made a number of conjectures about its 
geometric structure and relevance to 
the classification program.  (His original construction 
\cite{Ca01} was slightly different, compare also \cite{Lu02}.)
Here, we give an alternative
more direct construction of this special fibration, 
generalizing the theorem above and
resolve, among others, most of his conjectures on 
it relating to their structure and classification
outside of the main conjectures in Mori's Program.
We obtain an analog of the above theorem for
our canonical fibration in the sense that our 
fibration has an orbifold base 
(see the next sections for the definitions) of 
general type. It follows in particular 
(see Corollary~\ref{tech}) that this orbifold base
has no orbifold special fibrations with positive dimensional
fibers. Another analog we obtain (of no lesser importance)
is that if a manifold has a special 
fibration over a special orbifold base, then it is
itself special, the content of Proposition~\ref{spgt}
which is a consequence of the easy addition law of Kodaira
dimension. \\

We now state our main theorem 
as follows (see Theorem~\ref{Main} for
the full statement).

\begin{theorem} Let $\xd=(X,D)$ be a compact complex log-manifold.
Then a saturated subsheaf $L$ of $\Omega^p_X(\log D)$
with $\kappa(L)=\kappa'_+(X)=p$ as defined above is unique and 
gives rise, via its Iitaka fibration, to a special almost
holomorphic fibration ${b\,\!c}_{\!\xd}^{}:\xd\dashrightarrow Y$ 
with $Y$ $p$-dimensional.
Moreover, the line sheaf 
$\ {b \,\! c}_{\!\xd}^*K_U\hra \Omega^p_X(\log D)$ factors through 
(i.e., is contained in) $L$ and 
${b \,\! c}_{\!\xd}^{}$ factors through any special fibration
from $\xd$. In particular, if $\xd$ is not of (log) general type
(i.e., $p\neq \dim Y$), then $\dim Y=p< dim X$ and so
${b\,\!c}_{\!\xd}^{}$ has positive dimensional fibers.
\end{theorem}

As for terminology, an almost holomorphic fibration
$f:X\dra Y$ is a meromorphic map that restricts to a proper 
holomorphic fibration on a Zariski open subset of $X$. 
Hence the general fibers of $f$ make sense and
are log manifolds when endowed with the log structure from 
the $\xd$ given, which we call the general fibers of the 
corresponding orbifold fibration  $\xd\dra Y$. If these orbifold
general fibers are special, then we call the correponding 
orbifold fibration to be special.
Also, we call line sheaves $L$ with the above 
properties Bogomolov sheaves as in \cite{Ca02}
since they have been well studied by Bogomolov 
in \cite{Bo}. The notation ${b\,\!c}_{\!X}^{}$ is an 
adaptation of Campana's $c_X: X\dra C(X)$, which is called the 
core of $X$ in \cite{Ca02}. Although our construction
of ${b\,\!c}_{\!X}^{}$ is different from that of $c_X$ in
\cite{Ca02}, it follows directly from our main theorem that
they are in fact the same fibration. The terminology ``core''
is quite appropriate except for the ambiguity
in the definition of the orbifold base. The
ambiguity disappears once it is defined by Bogomolov sheaves
as we amply demonstrate in the next section. Thus our notation
${b\,\!c}_{\!X}^{}$, which we heartily thank Bogomolov for
his agreement for adopting, is meant to indicate that the 
fibration is defined by Bogomolov sheaves and can be called the 
Bogomolov-Campana fibration or the Bogomolov-Campana core of $X$
(or more simply the basic fibration (or reduction) of $X$). 
We then call $\kappa'_+(X)$ the Bogomolov-Campana dimension 
or simply the basic (or base core) dimension of $X$.\\

We actually resolve the full orbifold version 
of the above theorem, where 
the coefficients of the component decomposition of $D$ 
are allowed to be rational numbers not greater than one.
In this generality, we show for example
that the algebraic reduction is a special fibration
(irrespective of the orbifold structure imposed)
which in turn follows from an existence result for
a relative fibration that we state as follows
(see Proposition~\ref{GTR}  for the full statement):
\begin{theorem}
Let $\xd$ be a log manifold as above.
Let $h:X\ra Z$ be a fibration such that either $X$ is K\"ahler
or $Z$ is Moishezon. Then there exist a fibration 
$Y\ra Z$ and a special almost holomorphic fibration $f: X\dra Y$
such that $f_z={b\,\!c}_{\!X_z^\bd}: X_z\dra Y_z$, 
where $f_z$ is the restriction of $f$ to the 
general fiber $X_z$ of $h$. Consequently, the
algebraic reduction of $X$ is a special fibration
irrespective of the log-structure imposed on $X$. 
\end{theorem}
We note that in the compact complex category, meaningful relative 
fibrations were not known (indeed, for the relative Albanese 
and the relative algebraic reduction, known not) to exist except
for the relative Iitaka fibrations, even for a projectively
parametrized family. Also, relatively few (if any) facts 
were known for the structure of the general fibers of the
algebraic reduction besides having non-positive Kodaira
dimension. Our result indicates that the
general fibers of algebraic reductions have the
strictest non-hyperbolicity properties. 
Using this proposition, we
show via Proposition~\ref{spgt} that the fibration given
above in the general compact complex category
is nothing but the algebraic reduction composed 
with iterations of the orbifold
Iitaka and Mori's fibrations assuming that the latter
fibrations exist in the negative Kodaira dimension range,
see our discussion right after Conjecture~\ref{Mo}. 
This last assumption on the existence of a 
particular kind of special fibrations
in this situation is in fact one of the central 
problems in Mori's classification program and
we give a solution to this problem here in
the sense that any variety with negative Kodaira
dimension has a special fibration with positive
dimensional fibers.\\

A major aim of this paper 
is to introduce a framework for a theory of the
birational geometry of orbifolds in an elementary and
self-contained way as an accessible extension of
an important part of Mori's Program. We do this in the
most general setting of orbifolds as compact complex
manifolds endowed with boundary $\QQ$-divisors 
(not necessarily with
normal crossing support). Such an extension
does not seem to exist in any coherent or consistent
fashion in the literature especially with respect to
the birational geometry. This lack is the topic
of Section~3.30 of the book  
\cite{KM} by Koll\'ar and Mori.
We tackle this 
in the general compact complex setting
for the purpose of further studies
on them in complex geometry.  This is in fact a 
tall task for general orbifolds for which this paper
only hopes to have done some initial justice (our
original motivation to do this was to resolve some confusion
in the literature about the birational geometry of general 
orbifolds) even though the notions are fairly straightforward 
geometrically for smooth orbifolds (i.e., orbifolds whose
boudary $\QQ$-divisor has normal crossing support),
see Remark~\ref{orbrem} and Appendix A.
A possibly more natural but less elementary category 
(or notion) for this is that of ``stacks'' 
which has recently become popular and which the paper of S. Mori
and S. Keel (\cite{KeM}) and its followups seem to suggest. A
good theory for the birational geometry of stacks should be
of much interest for which this paper has suggested certain
formal approaches along the lines of Andr\'e Weil and will be
the topic of a forthcoming paper.
Still another possible direction would be A. Connes' theory of
noncommutative geometry since
our approach is essentially that of foliations
(compare the suggestive terminology used for 
the title of \cite{McQ}). \\

Our notion of orbifold is based on subsheaves
of the cotangent sheaf corresponding to fibrations
and includes log pairs of Mori's program except that
we restrict to nonsingular varieties. 
By allowing all rational multiplicities not greater
than one, one can ``approximate'' an orbifold 
by neighboring orbifold structures which are better
in some way. A particular instance of this is the following
theorem, which solves a conjecture of Campana in \cite{Ca02}
(and therefore also of Demailly-Peternell-Schneider in
\cite{DPS, DPS0} concerning the Albanese map) that 
projective manifolds with nef anticanonical bundle are 
special (see Theorem~\ref{-K}), 
of importance to their fundamental group and the
structure of their universal cover, see also \cite{DPS1, DPS2}.
We note that similar theorems of this type were always obtained via
a theorem of Miyaoka that employs positive characteristic 
techiques but our method does not involve positive 
characteristic and our result is more refined.

\begin{theorem}
Let $\xd=(X,D)$ be a projective log manifold whose log
anticanonical bundle $K_X(D)^\vee$ is nef. Suppose 
$f:X\smallsetminus D\dra Y\smallsetminus B$ is a strictly
rational dominant map (see \cite{Ii}) where $(Y, B)$ is a log
manifold. Then 
$$
\bar\kappa(Y\smallsetminus B):=\dim H^0(\Omega_Y(\log B))
\le \kappa(L_\fd) \le 0.
$$
Hence $\xd$ is special. Also,
the quasi-Albanese map of $X\smallsetminus D$
is a dominant morphism with connected general fibers which
are special and without multiple fibers outside 
a codimension two subset of $Alb(X\smallsetminus D)$.
\end{theorem}

The last part of the theorem is a simple fact derived 
via the structure theorem of quasi-abelian varieties
of \cite{Ka, Ue} (see also \cite{KV}).
A corollary of this is a complementary path
for a result of Mihai Paun which 
says that such manifolds without orbifold structures
have almost abelian fundamental group. 
This result has been obtained by Mihai Paun
(\cite{P}) using, among other strong results,
a result of Qi Zhang (\cite{ZQ}) concerning the
Albanese map that entails characteristic $p>0$ techniques
from \cite{MM1}. Our result refines that of Qi Zhang
but without involving characteristic $p>0$ 
thereby showing that this result of Mihai Paun can be derived 
more directly in line with complex geometry with 
improved understanding on the structure of such varieties.
We generalize this result of Paun to the log case as follows,
see Theorem~\ref{pi} for the full orbifold statement.

\begin{theorem}
Let $(X,A)$ be a projective 
log manifold. Suppose that there is
an effective divisor $A'$ with $A'_{red}=A$ such that
$K_X(A-\eps A')^\vee$ is nef for all sufficiently small
$\eps>0$, which holds in particular when $K_X(A)^\vee$
and $A'$ are nef.   Then $\pi_1(X\smallsetminus A)$ 
is almost abelian, or more specifically, a finite extension
of a finitely generated free abelian group.
\end{theorem}

This type of theorem for orbifold surfaces has been
quite intensively studied recently by different
techniques from that of ours and has been quite useful for
classification even in the two dimensional case of open surfaces
(by D. Q. Zhang for example). A similar consideration to the
derivation of this theorem also provides an alternative proof
of the simply connectedness of log-Fano varieties found in \cite{Ta}.\\

The key ingredients that go into the proofs of our main
results, apart from a well-suited generalized theory of 
orbifolds and their ``birational geometry'', are 
generalizations of the theorems of Kawamata, Viehweg
and Koll\'ar on the subadditivity of Kodaira
dimension for a fiber space (general type
case) to the orbifold category. 
We expect these generalizations to shed some light on the full 
subadditivity conjecture of Iitaka and hence to
the classification program, see Remarks~\ref{rm10} and \ref{rm11}
where the problem is reduced to special fibrations from special
varieties. 
\\

Many of the facts derived here about orbifolds are natural
generalizations of facts and resolutions of conjectures 
stated in \cite{Ca01} (also generalized to the full 
orbifold situation here). The perspective in that paper
is broader than ours and was motivated completely 
independently from ours. We urge the reader to consult
this very motivating paper for a different perspective 
and many more facts and questions
concerning the geometry of orbifolds. The 
preliminary sections are elementary in nature but we 
decided to include them as we cannot
find any satisfactory reference in the literature.\\

The author would like to thank J. P. Demailly, A. Fujiki,
B. Shiffman, E. Viehweg  and M. Zaidenberg for their 
keen interest and guidance in this problem, 
especially to B. Shiffman
to whom this paper owes its initiation and 
existence. A very special thanks is owed to F. Campana for 
communicating his inspiring preprint in November of 2001. 
It has served as an important force behind 
many of our results here.
He also expresses his thanks to Professors Y. I. Manin, 
F. Hirzebruch and Dr. Y. Holla for some useful discussions 
at the Max Planck Institute and finally to Y. Kawamata and 
F. Bogomolov for their suggestions and encouragements.
The first version of this paper was written in Osaka
University while this revised version is written at
the Max Planck Institute for Mathematics in Bonn,
the author would like to thank these respective
places for their very generous hospitality and support.
Finally the author would like to thank 
a most wonderful couple Barbara and Spyros
for their tolerance and kindness in making
the place like a home to him.

\section{The birational geometry of fibrations and orbifolds}

All objects in this paper
are compact complex manifold unless clear from the
context or otherwise specified. In this category, we adopt the
terminology of identifying holomorphic maps with morphisms, 
meromorphic maps with rational maps, and bimeromorphic maps 
with birational morphisms. We also identify vector bundles
with locally free sheaves. We assume that the reader is 
familiar with the rudiments of the theory of Kodaira dimension 
as found in Ueno \cite{Ue}, Mori
\cite{Mori} or Iitaka \cite{Ii} for example. 
The analytic version of the resolution of singularity 
theorem of Hironaka (\cite{Hi0}) proved by more 
elementary methods recently in \cite{BM}  
is used freely throughout this paper.\\

Let $f:X\ra Y$ be a fibration, i.e., $f$ is a surjective 
morphism with connected fibers. Then $f$ imposes an orbifold 
structure on $Y$ a part of which 
can be visualized by putting in the discriminant (branched)
$\QQ$-divisor $D(f)$ on $Y$ inherited from $f$ so that $Y$ 
becomes a branched orbifold as follows.

Let $\dv'(X)$ be the set of codimension one subvarieties of $X$
(as Cartier divisors), i.e.,
$$\dv'(X)=\{\ D\in \dv (X)\ |\ 
D\ \mbox{a reduced irreducible divisor on}\ X,\ D\neq 0\ \},$$
also known as the set of {\bf prime divisors}.
Given $D\in \dv'(Y)$, we may write $f^*D=\sum_i m_iD_i$ 
for $m_i\in \NN$ and 
$D_i\in \dv'(X)$. Then we define 
the (minimum) multiplicity of $f$ over $D$ by
$$m(D,f)=\min \{\ m_i\ |\ f(D_i)=D\ \}.$$
We note that the classical multiplicity $m_f(D)$ is defined
by replacing $\min$ above by $\gcd$ and is the more natural
multiplicity to consider for questions concerning the 
fundamental group, see \cite{Ca01}.
We define the $\QQ$-divisor
$$D(f)=\sum \Big\{\ \Big(1-\frac{1}{m(D,f)}\Big)D\ \ \Big|\ 
D\in\dv'(Y)\ \Big\}.$$
which is a finite sum 
supported on the discriminant locus $\Delta(f)$ of $f$.
We call the pair $(Y;D(f))$ an orbifold which, in order 
to avoid ambiguity of notation, we will now always
write as $Y\m D(f)$ and, in case no ambiguity arises,
also write as $Y^\bd$ for simplicity. We define 
the canonical $\QQ$-bundle and the Kodaira dimension of 
$Y^\bd=Y\m D(f)$ by
$$K_{Y^\bd} := K_Y( D(f))\ \ \mbox{and}\ \ 
\kappa(Y^\bd)=\kappa(Y\m D(f)):=\kappa(K_Y(D(f))),$$   
respectively. We  
also write $K_{Y^\bd} = K_Y + D(f)$ as a $\QQ$-divisor
by abuse of notation.\\

Note that all of these definitions above are clearly  
unchanged if we change the birational model of $X$
and therefore these definitions  make sense 
for meromorphic fibrations, i.e., surjective meromorphic
maps whose desingularizations have connected fibers.\\

Now, $Y^\bd=Y\m D(f)$ captures only part of the orbifold 
structure from $f$ and ignores, for example, all 
codimension two structures imposed on $Y$, the most 
serious being that the divisors on $X$ that project to
codimension two subsets of $Y$ are being ignored. To remedy this,
we consider a related object on $X$ which is one of
the most basic geometric invariants of $f$.

\begin{defn} Let $f:X\ra Y$ be a fibration. We define
$L_{(Y,f)}$ (or $L_f$ for short) to be the saturation 
of $f^*K_Y$ in $\Omega_X^r$ where $K_Y$ is the canonical
(or dualizing) sheaf of $Y$ and $r=\dim Y$. We set
$\kappa(Y,f)= \kappa(L_f)$ and call it the Kodaira
dimension of the orbifold $(Y,f)$ defined by $f$.
\end{defn}

We recall that the saturation of a rank $r$ subsheaf $L$ in a 
coherent sheaf $\Omega$ over $X$ means taking the largest 
rank $r$ subsheaf of $\Omega$ containing $L$. 
If $\Omega$ is locally free and $L$ is of rank one, then
the saturation of $L$ in $\Omega$ is locally free (but not
necessarily a vector subbundle except outside a codimension
two subset) and is uniquely determined
by any nontrivial meromorphic section of $\Omega$, see
\cite{MP} and the references therein. Furthermore, 
the saturation remains saturated on any open subset.\\

Since nontrivial global meromorphic differential forms are
preserved by birational maps, a saturated line subsheaf
$L$ of $\Omega_X^i$ defines a unique such subsheaf 
$L'\hra \Omega_{X'}^i$ on any birational model $X'$ of 
$X$ and we have $\kappa (L) = \kappa (L')$. So we may define
$L_{f}$ for any meromorphic fibration 
$f: X\dashrightarrow Y$. In fact, we may even
allow $Y$ here to be singular and then this notion of orbifolds, 
as far as we know, includes all previously known notions of 
orbifolds in complex geometry, see however \cite{KeM}. 
Note that $\kappa(Y,f):=\kappa(L_{f})$
depends only on the birational class of $f$ in this sense.
Moreover, $L_f\hra \Omega_X^{\bullet}$ continues to make sense for any rational
map $f:X\dashrightarrow Y$ as it is determined uniquely by its
associated foliation on $X$ and one can still make
sense of $(Y,f)$ as an ``orbifold'' or at least of 
$\kappa(Y,f):=\kappa(L_f)$ as its Kodaira dimension (in 
line with Grothendieck's philosophy that every map
is a ``fibration''). 

\begin{lem}[Birational Invariance]\label{A} Consider a 
commutative diagram of meromorphic maps (i.e., where 
compositions are defined as meromorpic maps and the
diagram commutes as meromorphic maps):
\begin{equation}\label{dm1}
\begin{CD}
X' @>v>> X\\
@Vf'VV @VVfV\\
Y' @>>u> Y
\end{CD}
\end{equation}
where  $u$ is generically
finite and surjective. Then there is a natural
inclusion
$$v^\dagger: H^0(X, L_f^m)\hra H^0(X', L_{f'}^m)$$
for all $m\in \NN$ and thus
$$\kappa(Y',f')\ge \kappa(Y,f).$$

If moreover $v$ is birational, then
$v^\dagger$ is an isomorphism and so
$$\kappa(Y',f') = \kappa(Y,f).$$

\end{lem}

\noindent
{\bf Proof:} By taking the Stein factorization of $f$
and $f'$, we may assume that they are meromorphic fibrations.
Note that $f\circ v=u\circ f'$ is a well defined
surjective rational map by commutativity of the 
diagram. Hence there is a Zariski open set $U$ in $X'$
where $v$ is defined and where $L_{f'}=v^*L_f$. Let
$V$ be the open set where $v$ is defined. Then,
for all $m\in \NN$, we have
\begin{equation}\label{in1}
L_{f'}^m \supset v^*L_f^m
\end{equation}
as sheaves on $V$ by the commutativity of diagram~\ref{dm1}. 
Let $r:\tilde X\ra X$ be a birational modification 
so that $\tilde f=f\circ r$ is a morphism. Then sections
of powers of $L_f$ correspond to those of $L_{\tilde f}$.
Since $L_{\tilde f}$ is trivial on the general fibers of 
$\tilde f$ and $f\circ v$ is surjective, we have an injection 
$H^0(X, L_f^m)\hra H^0(V, v^*L_f^m)$. 
As $X'\smallsetminus V$ is codimension
two or more, the Hartog extension theorem gives the first
part of the lemma. The second part follows similarly.
$\BOX$\\

If the condition of the last part of the lemma is satisfied,
we say that $f$ and $f'$ are 
{\bf birationally equivalent fibrations}.\\

Now given an $f'$ as in Lemma~\ref{A}, it is clear from 
the Hartog extension theorem applied to $Z'$ that 
$\kappa(Y,f)=\kappa(Y',f')\le \kappa(Y'\m D(f')).$
It is in fact always possible to achieve equality here
by choosing an appropriate $f'$ as follows:

\begin{defn} A fibration $f':X'\ra Y'$ is called {\bf admissible}
if there is a birational morphism $v: X'\ra X$ with $X$ 
smooth such that every reduced divisor $E$ in $X'$ with
$f'(E)$ of codimension two or higher is exceptional with respect to $v$.
The union of such divisors will be denoted by $E(f')$.
\end{defn}

\begin{lemma}\label{B}
Let $f: X\ra Y$ be a fibration. Then there is a natural
inclusion $$ H^0(X, L_f^m)\hra H^0( Y, K_\yd^m)$$
for $\yd=Y\m D(f)$ and all $m>0$ divisible by
the multiplicities in $\yd$. In particular,
$$\kappa(Y, f) \leq \kappa(Y\m D(f)).$$
Also, one can find a commutative
diagram as in diagram~\ref{dm1} with $u,v$ birational and
onto so that the birationally equivalent fibration $f'$
is admissible.  If $f':X'\ra Y'$ is an admissible fibration
and $m$ is divisible by the multiplicities of $D(f')$, then
$$H^0(X', L_{f'}^m)= H^0( Y, K_{Y'\m D(f')}^m)\ \ \text{and }\ \ 
\kappa(Y',f')= \kappa(Y'\m D(f')).$$
\end{lemma}

\noindent
{\bf Proof:} The construction of a birationally 
equivalent admissible fibration is a basic starting point in
the usual classification theory via Kodaira dimension and the
Iitaka fibration. It is
achieved by resolving the singularities of the flattening of
$f$ (by grace of \cite{Hi} and, in the projective category,
also by \cite{GRay}). For the last statement, 
we have by definition (with $r=\dim Y'$) the inclusion
\begin{equation}\label{sat1}
{f'}^*K_{Y'}(D(f'))^m\hra L_{f'}^m\ 
\Big(\hra (\Omega^r_{X'})^{\otimes m}\Big)
\end{equation}
outside $O\cup E(f')$ where $O$ is a subset of
$X'$ of codimension two or higher
contained above $\Delta(f')$ and
$m$ is a positive multiple of all relevant multiplicities.
Moreover, this inclusion is an isomorphism on an open subset
of $X'$ that surjects to the complement of a 
subset of $Y'$ of codimension two or higher. Hence 
$$H^0(Y', K_{Y'}(D(f'))^m) \hra H^0(X, L_f^m)
=H^0(X', L_{f'}^m)$$
by the Hartog extension theorem applied to $X$ and the 
reverse inclusion by the Hartog extension theorem applied
to $Y$. The first two statements now follow using the 
inclusion (\ref{sat1}) for the fibration $f$ (rather than $f'$)
and by applying the Hartog extension theorem to $Y$ since
the said inclusion is an isomorphism on an open subset
of $X$ that surjects to the complement of a comdimension
two subset of $Y$.
$\BOX$\\

The inclusion (\ref{sat1}) comes from the fact, true 
outside the inverse image of a codimension two subset 
$S$ of $Y$, that $L_{f'}={f'}^*K_Y((df'))$
where $(df')$ is a divisor on $X$ given by (cf. \cite{Re})
$$\sum \{\ ({f'}^*D)-({f'}^*D)_{red}\ |\ 
D\in\dv'(Y)\ \}$$ and the easily verified 
fact that $(df')\ge {f'}^*D(f')$ with equality
achieved on an open set that surjects 
to $Y\smallsetminus S$, cf. Lemma~\ref{MR}.

\begin{rem}
$\bullet$ An $f'$ satisfying the conclusion of Lemma~\ref{B} is 
called admissible by Campana in \cite{Ca01}. (His definition of
$\kappa(Y,f)$ is different, but is equivalent to ours by
our result above.) In that paper, he posed the 
problem of finding an explicit condition satisfying his 
notion of admissibility. Our lemmas above give a 
solution to this problem. \\
$\bullet$ Readers not interested in the
full orbifold version of our results can safely skip to
Section~5 by ignoring the $\bd$'s on the total space
(of fibrations) and some related terminologies.
\end{rem}

\section{Log pairs as orbifolds and their fibrations}

We now generalize our notion  of orbifolds to the full log
version as in Mori's Program by allowing the multiplicities 
to be rational numbers as well as $\infty$. 
That is, $Y\m B$ is a log orbifold (also known as a log pair
or log variety, see \cite{KM, KMM})
or simply an orbifold if $B$ is a $\QQ$-divisor
of the form $\sum_i (1-1/m_i) D_i$ where 
$0< m_i\in \QQ\cup\{\infty\}$ and $D_i\in \dv'(Y)$ for all $i$.
It is called a standard pair or a {\bf standard orbifold }
if $m_i\in \NN\cup\{\infty\}$ for all $i$. We call $Y\m B$ 
an {\bf effective orbifold} if $B$ is an effective (possibly
the trivial) divisor. Note that standard orbifolds are 
effective. Although many (if not most) of the results in
this paper are worked out for arbitrary orbifolds, the
important results in this paper essentially only deal 
with effective orbifolds.

It is also
possible to consider real multiplicities $m_i>0$ as in
Mori's Program at the cost of introducing more technical
machinery for real divisors in addition to notations 
such as the round-up. We will restrict to $\QQ$-divisors
in this paper for simplicity but the
$\RR$-divisor version is obtained similarly.
 
Given an orbifold structure $B$ on $Y$, we call $B$
the boundary divisor of $Y^\bd=Y\m B$ and write it also
as 
$$B=\bd (\yd) =\sum \Big\{\ \Big(1-\frac{1}{m(D\cap B)}\Big)D\ \ \Big|\ 
D\in\dv'(Y)\ \Big\}.$$ 
We also write $m(D)=m(D,Y^\bd)=m(D\cap B)$.
Note that for $D\in \dv'(Y)$ we have
$$0< m(D\cap B)\leq \infty.$$
Also $m(D\cap B)=1$ if 
$D\not\subset B_{red}$, and $m(D\cap B)=\infty$ if $D\leq B$
as $\QQ$-divisors. We define 
$$\mu(D\cap B)=1-\frac{1}{m(D\cap B)}\leq 1$$ 
In general, if $B'$ is another $\QQ$-divisor on $Y$, we define
$$ B\cap B' = \sum \Big\{\ (\min\{\mu(D\cap B),\mu(D\cap B')\})D
\ \ |\ D\in\dv'(Y)\ \Big\}.$$ 
We also set, as before,
$$K_{Y^\bd} := K_Y(B)=K_Y(\bd (\yd))\ \ \mbox{and}\ \ 
\kappa(Y^\bd)=\kappa(Y\m B):=\kappa(K_Y(B)),$$ 
respectively.
If $m_i=\infty$ for all $i$ and $B=\sum_i D_i$
is a (simple) normal crossing divisor, then $Y^\bd$ is a
{\bf log-manifold} in the usual sense and one can define the 
log-tangent sheaf $T(Y,-\log B)$ (i.e., the orbifold tangent sheaf) 
of  $Y^\bd$ as the sheaf 
of vector fields on $Y$ leaving $B$ invariant. It is
locally free with dual
the log-cotangent sheaf $\Omega_{Y^\bd}=\Omega(Y,\log B)$.
We note $K_{Y^\bd} = \det\Omega_{Y^\bd}$.
Hence, this notion of orbifolds includes open subsets
in the sense of log-manifold.

To be consistent with the theory of log-manifolds, we say that  
$Y\m D$ is a {\bf smooth orbifold} or simply smooth
if $D_{red}$ is simple normal 
crossing. This allows, among others, a good notion of orbifold
tangent and cotangent sheaves. A smooth orbifold $X^\bd=X\m A$ 
gives rise to a log-manifold in a natural way, which we denote
by $$ X^{|\bd|}:= X\m A_{red}.$$

\begin{rem}\label{orbrem}
The above definition of orbifolds via boundary
with arbitrary singularities is only provisional 
since only smooth orbifolds make good sense for
the purpose of geometry. This is because the tangent and
cotangent sheaves for them exist (on some space) and 
are locally free in this case (see the next section
for a view to this).
Our purpose in allowing arbitrary singularities
in the above definition is to clarify
certain issues in the literature as well as to
explore some limits of such a definition.
However, it should be mentioned that these sheaves,
for example in the case the orbifold boundary is reduced, 
can be locally free without being normal crossing in the
case of ``free boundaries'' introduced and studied 
by K. Saito (\cite{Saito}). The free boundary orbifolds,
which allow for example ordinary triple points in the
boundary divisor of a complex surface, 
do play a role in the search for a good model. 
But {\bf for all practical purposes, 
it is sufficient to restrict attention to smooth
orbifolds} whenever it
comes to orbifold maps or morphism. 
This is what we will do implicitly 
from now on as it is certainly no loss of generality
from the viewpoint of birational geometry.
\end{rem}

Smooth orbifolds are important in 
considering  the birational geometry of orbifolds,
for which we first give a preliminary discussion.

\begin{defn}\label{dk} Let $u: Y_1\ra Y$ be a 
generically finite surjective morphism where
$Y_1^\bd$ and $Y^\bd$ are orbifolds. We then write 
$u^\bd:  Y_1^\bd\dashrightarrow \yd$ as an orbifold 
rational map.
We say that $u$ gives rise to a {\em  K-map}
$$u^\bd: Y_1^\bd\dra \yd$$ if 
$K_{Y_1^\bd}-(u^* K_{Y^\bd})$ as a $\QQ$-divisor
defined by $(du)$ is effective.
If further $u$ is birational, we say that $u^\bd$ is 
K-birational if 
$$K_{Y_1^\bd}=u^* K_{Y^\bd}+ \sum_i a_iE_i$$
where each $E_i$ is an exceptional divisor of $u$ and $a_i\ge 0$
is a rational number. If $a_i=0$ for all $i$, then we say
that $u^\bd$ is strictly {\em K-birational}.
\end{defn}

It should be understood above that the trivial divisor is
effective.\\

K-birational maps are used in the characterization
of log-terminal and of log-canonical 
singularities of log-pairs in Mori's Program
(our notion of orbifold can be restricted 
to that of log-canonical pairs and most of our results
for smooth orbifolds hold in this more general context),
see \cite{KM} for example. It is a good definition of orbifold 
birational morphisms outside the singularities
of the boundary divisors of the target orbifold, 
see Lemma~\ref{K-map}. If one is only interested
in preserving the Kodaira dimension for birational
maps, it is quite sufficient.
However, further conditions must be 
imposed in general to ensure the morphism condition.
The problem is already present in the 
absolute case since a K-map  $\fd: X^\bd\dra \yd$ between
log-manifolds need not be a log-morphism, even if it is 
K-birational. A simple example is the following:

\begin{ex}
Let $f:X\ra Y$ be the blowup of $Y=\CC^2$
at $0\in \CC^2$, $B$ is the $y$-axis and $A$ is the 
strict transform of $B$. Then $\fd: X\m A\dra Y\m B$
is a K-map but not a log morphism.\\
\end{ex}

We now consider orbifold fibrations
(as a prelude to orbifold 
morphisms).
 
\begin{defn}
Let $f:X\ra Y$ be a fibration
with $X^\bd=X\m A$ an orbifold. We denote this by
$f^\bd : X^\bd\ra Y$ and call it an orbifold fibration, or
fibration for short. We define
$L_{(Y,\fd)}$ (or $L_\fd$ for short) to be $L_f(A\cap f)$ where
$$
{A\cap f}=\sum \Big\{\ \Big(1-\frac{1}{m(D\cap A)}\Big)D\ \ \Big|\ 
f(D)\neq Y,\ D\in \dv'(X)\Big\}.
$$ 
We set
$$\kappa(Y,\fd)= \kappa(L_\fd),$$
which we call the 
{\bf Kodaira dimension of the orbifold} $(Y,\fd)$.
\end{defn}

The $f^\bd$ as given in the definition imposes  
a branched orbifold structure on $Y$ as follows.
Given $D\in \dv'(Y)$, we may write $f^*D=\sum_i m_iD_i$ 
for $m_i\in \NN$ and 
$D_i\in \dv'(X)$. Then we define 
the (minimum) multiplicity of $f^\bd$ over $D$ by
$$m(D,f^\bd)=\min \{\ m_i\, m(D_i\cap A)\ |\ f(D_i)=D\ \}.$$
Then the $\QQ$-divisor 
$$D(f^\bd)=D(f,\bd (X^\bd))=D(f,A):=\sum \Big\{\ 
\Big(1-\frac{1}{m(D,f^\bd)}\Big)D\ \ \Big|\ 
D\in\dv'(Y)\ \Big\}$$
defines the {\bf orbifold base} $\yd=Y\m D(f^\bd)$ of $\fd$. 
If $B\leq D(f^\bd)$ as $\QQ$-divisors, then we say that
the orbifold fibration $\fd$ (and $Y\m D(f^\bd)$)
dominates $Y\m B$.\\

We remark that we could have arrived at the full log version 
of orbifolds by considering the branched orbifold base of
a fibration from a log-manifold this way 
(see Definition~\ref{orf}). The following are
immediate consequences of the definition.

\begin{lem} \label{df}
Let $\fd:\xd\ra Y$  be a fibration with $\xd=X\m A$. 
Let $B=D(\fd)$ and $\yd=Y\m B$. Then
above the complement $U$ of a codimension two subset of $Y$,
we have
$$(df) + A\ge f^*B\ \ \ 
\mbox{where}\ \ \ (df)=f^*B-(f^*B)_{red}$$
as $\QQ$-divisors on $X$ with equality achieved 
everywhere on $U$ (compare Lemma~\ref{K-map},
the paragraph after Lemma~\ref{B} and the slightly
more general equations~\ref{a} and \ref{b}). \BOX
\end{lem}

\begin{lem} \label{comp} Let
$B=D(\fd)=D(f,A)$ for a fibration $\fd:X\m A\ra Y$, 
$\gd: Y\m B\ra Z$ a fibration and $h=g\circ f$.
Then $D(\gd)\geq D(h^\bd)$, i.e., $D(g, D(f,A))\geq D(g\circ f, A)$.
If the exceptional part of $A$ with respect to $f$ is reduced, then
equality holds. \BOX
\end{lem}

Now, all the above definitions for an orbifold 
fibration actually make sense for an arbitrary (rational)
map $f$ from an orbifold, such as $A\cap f$ 
and $L_\fd$ for example,
because $f$ is a morphism outside a codimension two subset.
We will adopt this more general usage of notation and, 
in the next definition, even allow $\fd$ to have noncompact
open sets as domains, where our definitions above generalize
straightforwardly. \\

Observe that if $\fd:\xd\ra Y$, $\gd:\yd\ra Z$ 
are surjective and 
$\hd=g\circ \fd$, then $f^*L_g\subset L_h$ is a 
nontrivial subsheaf so that $L_h-f^*L_g$ is an 
effective Cartier divisor on $X$, playing in some sense
the role of ``$(df)\cap h$.'' Hence $L_\hd-f^*L_\gd$ is
uniquely given by a $\QQ$-divisor. A $\QQ$-divisor $D$ is
called effective, written $D\geq 0$, if it is a 
nonnegative linear combination of reduced divisors.
More generally, we will drop the surjectivity condition
on $f$, the compactness assumption on $X$ and $Y$
and even allow $g$ to be meromorphic and $Z$ to be
singular. We only 
require that $\hd=g\circ \fd$ is defined on an open
set and extends to a holomorphic map onto 
a non-empty open subset of $Z$ and that the graph of 
$g$ restricted to the image of $f$ is irreducible. 
We then say that the map $g$ {\bf composes holomorphically}
with $f$.

\begin{defn} \label{defmo}
Let $u: Y_1\ra Y$ be a morphism where
$Y_1^\bd$ and $Y^\bd$ are orbifolds. We then write 
$u^\bd:  Y_1^\bd\dashrightarrow \yd$ as an orbifold 
rational map. We say that $u$ gives rise to an
{\bf orbifold morphism} 
$$u^\bd: Y_1^\bd\ra \yd$$ 
(or simply that $u^\bd$ is a morphism) if 
for every open subset $U$ of $Y$ and 
a meromorphic map $g: U\dra Z$
that composes holomorphically with $u_U:=u|_{u^{-1}(U)}$, 
we have $$L_{g_1^\bd}-(u^* L_{g^\bd})\geq 0.$$ 
If equality is attained everywhere, then we say that
$Y_1^\bd$ has the {\bf induced orbifold structure} from 
$\yd$ via $u$, denoted by $u^\circ \bd \yd$, and 
if further $u$ is an inclusion, then we say that
$u^\bd$ is an {\bf orbifold inclusion}. 
If we have equality outside a codimension 
two subset of $Y$ independent of $g$ above and $u$ is 
birational, then we say that $u^\bd$ is a 
{\bf birational morphism} or simply that $u^\bd$ is birational.
An orbifold morphism $u^\bd: Y_1^\bd\ra \yd$ is called 
{\bf smooth} at a point $p\in Y_1$, if $p$ has a
neighborhood in which $u$ is equidimensional, the 
reduced fiber
$F$ through $p$ is smooth and in which $F$ is not 
contained in, but is transversal to $D_{red}$ where
$$D:= \bd Y_1^\bd - u^\circ \bd \yd=
 \bd Y_1^\bd - \big(u^* \bd \yd - (L_u-u^*K_Y)\big)$$
in the sense that $D_{red}$ is
locally trivial with respect to $u$ in a neighborhood of $p$. 
If $u^\bd$ is a morphism and $u$ is a fibration, then
we call $u^\bd$ a {\bf fibration morphism} and we denote its
discriminant locus by
$\Delta(u^\bd),$ the image of the nonsmooth locus of $u^\bd$.
\end{defn}

We note that if $\bd \yd\geq 0$ and 
$Y_1^\bd$ and $Y^\bd$ are smooth orbifolds,
an orbifold morphism $u^\bd: Y_1^\bd\ra \yd$ must satisfy 
$$u^{-1}(\bd Y^{|\bd|})\subset \bd Y_1^{|\bd|}.$$
It is immediate that the discriminant locus $\Delta(u^\bd)$
of a fibration morphism $u^\bd: Y_1^\bd\ra \yd$ 
is always a proper algebraic subset outside of which the fibers
are differentiably isomorphic as ``differentiable''
orbifolds. Also, if $i:F\hra Y_1$ is the inclusion of a
fiber of $u$ outside $u^{-1}(\Delta(u^\bd))$, then the
induced orbifold structure on $f$ is just $i^* \bd Y_1^\bd$
and it has normal crossing support if $Y_1^\bd$ does.\\

Of course, the above condition defining an orbifold
morphism is really a local one and 
it is not necessary to consider all maps $g: Y\dashrightarrow Z$
but only a subset of maps each of which takes some
components of $\bd(\yd)$ as fibers, see Lemma~\ref{K-map} and 
Lemma~\ref{two} and its proof for further simplifications
in particular circumstances.
In fact, the notion of orbifold morphism for smooth
orbifold is very simple geometrically and has much simpler
characterizations, see Appendix A. However, the 
above definition is a more practical and direct approach 
to our main results.
\\

A particular class of orbifold morphisms consists of log-morphisms
between log-manifolds as we show below. 
In this case the characterization 
is particularly simple in 
sharp contrast to the general case.
Hence given two log-manifolds $X\m A$ and $Y\m B$, a map
$\fd: X\m A\dra Y\m B$ is an orbifold morphism if and only if
it is a log-morphism, namely $f^{-1}(B)\subset A$.
This comes from the fact that we actually have an
inclusion of the log-cotangent sheaves
in this case. 

\begin{lem}\label{ful} Let $f:X\ra Y$ be a morphism where $X^\bd$
and $\yd$ are log-manifolds. Then $\fd$ with respect to 
these orbifold structures is an orbifold morphism if and 
only if it is a log morphism.
\end{lem}

\noindent
{\bf Proof:} We only prove the ``if'' direction since
the other direction is not needed. The result then
follows from the next self-evident lemma and the fact 
that a saturated line subsheaf of a vector bundle on $X$ is 
locally free as it must be reflexive and is, 
outside a codimension two subset of $X$, a subbundle
(see for example \cite{MP}).  \BOX

\begin{lem}\label{fu} 
Let $L$ be a saturated line subsheaf of $\Omega_X^i$
and $X\m A$ a log-manifold. Then the saturation of $L$
in $\Omega_X^i(\log A)$ is $L(A')$ where $A'$ consists
of components $D$ of $A$ whose normal bundles 
$N_D$ satisfy $N_D^*\wedge L=0$
in  $\Omega_X^{i+1}(\log A)|_D$. Hence given a fibration
$\fd: X\m A\dra Y,$  $L_\fd:= L_f(A\cap f)$ is the 
saturation of $L_f$ in $\Omega_X^r(\log A)$ where
$r=\dim Y$.
\BOX
\end{lem}

\begin{defn}\label{adapt}
A fibration $f^\bd:X^\bd \ra Y$ is called {\bf admissible}
if there is a birational orbifold morphism 
$v^\bd: X^\bd\ra X^\bd_0$ 
such that every reduced divisor $E$ in $X$ with
$f(E)$ of codimension two or higher is exceptional with respect to $v$.
It is called {\bf very admissible} if further $X^\bd$ is smooth
and {\bf adapted} if furthermore $Y\m D(\fd)$
is smooth. If moreover $\fd:X^\bd \ra Y\m D(\fd)$ is a morphism,
then $\fd$ is called {\bf well adapted}.
An orbifold fibration is called {\bf good} if it is 
a morphism to its orbifold base. An orbifold $\yd$ is called
an {\bf admissible orbifold} if it is the orbifold 
base of an admissible fibration $\fd :\xd\ra Y$
from a smooth orbifold. It is called effectively admissible
or {\bf e-admissible} if $\bd \xd\cap E(f)$ is
effective.
\end{defn}

It is important to note that, by definition, {\em smooth
orbifolds are e-admissible.}
The following are key lemmas that allow for a 
good theory of birational geometry for smooth orbifolds.
Recall that the exceptional divisor $E(f)$ of a morphism
$f$ is the union of prime divisors whose images under $f$
are codimension two or higher.

\begin{lem}\label{ford} {\rm (1)} Let  $f:X \ra Y$ be a 
generically finite surjective morphism with
$Y^\bd=Y\m B$ a smooth orbifold. Suppose 
$X\m \bar A$ is a smooth orbifold for 
$\bar A= (f^* B)_{red}$.
Then $X$ can be endowed with a smooth (and standard)
orbifold structure $A$ supported on $\bar A$
(not unique in general) so that $A=f^*B$ outside
the exceptional divisor of $f$ and so that
$\fd:\xd\ra \yd$ is a morphism.
In particular, one can choose $\fd$ 
to be birational in the case $f$ is 
birational. 

{\rm (2)} Let  $\fd:\xd \ra Y$ be a fibration with $Y^\bd=Y\m B$
smooth for $B=D(\fd)$. Assume that $f^{|\bd|}$ is a log morphism. 
Then the multiplicities of $\xd$ along $E(f)$ can be increased if 
necessary so that the resulting 
orbifold $X^{\bd'}$, which gives rise via $f$ to the same 
orbifold base $B$, gives rise to an orbifold morphism
$f^{\bd'}: X^{\bd'}\ra \yd$ via $f$. If further 
$\fd$ is admissible, then $\xd$ and $X^{\bd'}$
are birationally equivalent as orbifolds.
\end{lem}

\noindent
{\bf Proof:} 
This follows directly
by considering the absolute case 
$f^{|\bd|}: X\m \bar A\ra Y\m B_{red}$,
which is a log-morphism 
and hence an orbifold morphism by Lemma~\ref{ful}.  
Here, we take $\bar A= (f^* B)_{red}$ for (2)
as well.
\BOX\\

We note that the choice of multiplicities can always be taken
to be standard and not $\infty$ on any component of
$E(f)$ that does not map into a divisor with 
multiplicity $\infty$.\\

 We generalize
the notion of orbifold morphisms in the special
case of maps between orbifold base of fibrations as follows.

\begin{defn} \label{d0}
Given a commutative diagram of morphisms
\begin{equation}\label{dm2}
\begin{CD}
X_1^\bd @>{v^\bd}>> X^\bd\\
@V{f_1^\bd}VV @VV{f^\bd}V\\
Y_1 @>>u> Y
\end{CD}
\end{equation}
where $f_1^\bd$ and $f^\bd$ are fibrations with $f^\bd$ admissible and
where $v^\bd$ gives a morphism of orbifolds $X_1^\bd$ and $X^\bd$.
Then we say that $u$ gives rise to a {\bf premorphism}
$$u^\bd: Y_1\m D(f_1)\ra Y\m D(f).$$
If $v^\bd$ and $u$ are birational, then we say that
$u^\bd$ is {\bf prebirational}. A composition consisting
of premorphisms and morphisms is called a {\bf promorphism} and
if the factors are birational, we say that it is 
{\bf probirational}.
\end{defn}

\begin{lem}\label{none} 
Given a fibration  $\fd: X^\bd\ra Y$ with $\yd=Y\m D(\fd)$,
we have for $m>0$ divisible by all the relevant
multiplicities that 
$$H^0(X, L_\fd^m)\hra H^0(Y, K_\yd^m)\ \ \text{and so }\ \ 
\kappa(Y,\fd)\leq \kappa(\yd).$$
Moreover, if $\fd$ is an admissible fibration,
then
$$H^0(X, L_\fd^m)= H^0(Y, K_\yd^m)\ \ \text{and so }\ \ 
\kappa(Y,\fd)=\kappa(\yd).$$
{\bf\em If $\xd$ is smooth, one can find} a commutative
diagram as in diagram~(\ref{dm2}) with 
$v^\bd$ and $u$ birational and onto
so that the {\bf\em birationally equivalent 
fibration} $f_1^\bd$
is {\bf\em well adapted and} hence {\bf\em good}. Also 
$$
\kappa(Y_1,f_1^\bd)= \kappa(Y, \fd).
$$
The same holds for $\xd$ admissible with $v^\bd$
weakened to prebirational.
\end{lem}

\noindent
{\bf Proof:} The proof is essentially the same as Lemma~\ref{B}
and Lemma~\ref{A} using Lemma~\ref{ford}. The proof of the last 
part is also contained in the proof of Proposition~\ref{MR}. It is
sufficient to observe here that an adapted fibration is well
adapted if one choose the multiplicities of the components of
the exceptional divisors of
$v$ as given in Definition~\ref{adapt}
above to be $\infty$ and that
$$1-{1\over m} +{1\over m}\,\big(\,1-{1\over n}\,\big)=1-{1\over mn}.
\ \ \ \BOX$$ \\

The cautious reader would note, in comparing this lemma with
Lemma~\ref{B}, that the birationally equivalent admissible fibration
$f'$ there cannot be chosen in general to be an orbifold morphism
to its orbifold base, the crucial difference being that
we are allowing orbifold structures on the total space of a
fibration in the above lemma. 
We record here the following obvious but useful fact.

\begin{lem}\label{A*}
Consider the commutative diagram~(\ref{dm2}) of morphism above
with $X^\bd$, $X_1^\bd$ (smooth) orbifolds, $f$, $f_1$ fibrations
and $u$, $v^\bd$ birational. 
If $\fd$ is admissible, then so is $f_1^\bd$. $\BOX$\\
\end{lem}

Note that the class of premorphism includes that of 
orbifold morphisms and that of adapted fibrations,
the orbifold analogue in some sense of
strictly rational maps (c.f., \cite{Ii}). 
In view of Lemma~\ref{ford}, the category of
admissible orbifolds and promorphism may serve in considering
the birational geometry of (not necessarily smooth) orbifolds.
The use of this notion is not indispensable in this paper
although it allows for a larger class for which
most of the theorems in this paper hold.  It would
be of much interest to have a more intrinsic and
direct characterization of a good class of maps 
that includes premorphisms for this purpose. It is
easily seen from the definition that a composition
of orbifold morphisms is an orbifold morphism. 
This means in 
particular that a (pro)morphism from $\yd$ remains so
when restricted to the general fibers of a fibration
$g:Y\ra Z$. This is easily verified to be true also
for a premorphism. We record these facts as follows.\\

\begin{lem} \label{premis}
Let $\fd: X^\bd\ra \yd$, $\gd:\yd\ra Z^\bd$
and $\hd=\gd\circ\fd$. If $\fd$ and $\gd$ are (pre)morphisms,
then so is $\hd$. If $\fd$ is a morphism and $\gd$ a 
promorphism, then $\hd$ is a promorphism.
If $g$ is a fibration and 
$f^\bd_z:X^\bd_z \ra Y^\bd_z$ the 
restriction of $\fd$ to the general fiber $X^\bd_z$
of $h$, then $f^\bd_z$ is admissible, a morphism,
a premorphism 
or a promorphism respectively if so is $\fd$. \BOX
\end{lem}

It is understood whenever we speak of 
a premorphism $u^\bd:Y_1^\bd\ra \yd$ from now on that we have a 
commutative diagram as above where $Y_1^\bd$ and $\yd$ are
the orbifold base of fibrations
(in practice, from smooth orbifolds or from log-manifolds)
$X_1$ and $X$ respectively with the latter assumed
to be admissible and where $v^\bd$ is a morphism.
We remark that an orbifold map is really properly understood by
the fact that the differential of the map gives a 
(linear and holomorphic) map of the 
orbifold tangent bundles (which can be understood as the
``tangent bundle'' of the corresponding V-manifold for 
example) or alternatively as done in the next section
by orbifold cotangent sheaves. But for our purpose here, 
it is sufficient to take this more indirect  but perhaps
more practical approach. \\

The connection between premorphisms, morphisms 
and adapted or good orbifold
fibrations can be seen from the following lemma.

\begin{lem} \label{na}
In the situation of definition~\ref{d0},
suppose $\fd$ is a good (or more strongly, a well adapted)
fibration. Then $u$ gives rise to an orbifold morphism
$$u^\bd: Y_1\m D(f_1^\bd) \ra Y\m D(\fd).$$ 
\end{lem}

\noindent
{\bf Proof:} Set  $Y_1^\bd=Y_1\m D(f_1^\bd)$,
$\yd=Y\m D(\fd)$ and $B=D(\fd)$. Let $V$ be a 
neighborhood of a point $p\in B_{red}\cap u(Y_1)$ and 
$g: Y\cap V\dra Z$ a meromorphic map that compose 
holomorphically with $u|_{u^{-1}(V)}$ which we still
denote by $u$ by abuse of notation. We need to show that
the $\QQ$-divisor $L_{g\circ u^\bd}$ arising from $du$ 
is effective on an open subset of the form $W=u^{-1}(U)$ 
for a small neighborhood
$U$ of $p$ in $V$ on which $L_\gd$ is trivial. 
Now a nowhere vanishing 
section of $L_\gd^l|_U$ gives rise to a section of 
$L_{g\circ \fd}^l|_{f^{-1}(U)}$ as $\fd$ is a morphism, 
which in turn gives rise to a section of 
$L_{g\circ f\circ v^\bd}^l=L^l_{g\circ u\circ f_1^\bd}$ above  
$W$  as $v^\bd$ is a morphism 
and finally to a section of
$L_{g\circ u^\bd}$ on $W$ by the 
Hartog extension theorem
(c.f. Proposition~\ref{MR}), where $l$ is chosen to divide all 
multiplicities so that the $l$-th tensor powers above makes
sense as line bundles. So the desired conclusion follows. \BOX\\

\begin{defn}\label{orf} Let $f^\bd:X^\bd \ra Y$ be an
admissible fibration with 
$X^\bd=X\m A$ a log-manifold. 
Given a fibration $\gd: \yd=Y\m D(\fd) \ra Z$, we define 
$L(Z,g^\bd; \fd)$ to be the saturation 
of $f^*g^*K_Z$ in $\Omega^s(X,\log A)$ where $s=\dim Z$
and we let
$$L(Z,g^\bd\circ \fd) = L(Z,g^\bd; \fd).$$
We set $\kappa'(Z,\gd;\fd) = \kappa(L(Z,g^\bd; \fd))$ 
and $\kappa'(Z, \gd)= \kappa(L(Z,g^\bd; {\rm id}))$. 
\end{defn}

A crucial point is the fact that 
$(Z,\gd;\fd)$ is independent of $\fd$ in the sense
that $\kappa'(Z,\gd;\fd)$ is and in fact we will
identify $\kappa'(Z,\gd;\fd)$ with $\kappa(Z,\gd).$
We now show this fact by first warming up with the 
following special case.

\begin{lem}\label{XYZ}
Assume that $\fd:X^\bd \ra Y$ is an adapted fibration with 
$X^\bd=X\m A$ a log-manifold.
Given a fibration $g:Y\ra Z$,
define $L_Y$ (resp. $L_X$) to be the saturation of 
$g^*K_Z$ in $\Omega(Y,\log D(\fd))$ (resp. $f^*g^*K_Z$ in
$\Omega(X,\log A)$). Then $L_X=L_{g\circ \fd}$ and 
$\kappa(L_X)=\kappa(L_Y)$.
\end{lem}

\noindent
{\bf Proof:} The first equality follows from 
Lemma~\ref{fu}. For the second,
let $s=\dim Z$ and $v^\bd: X^\bd\ra X^\bd_0$
be a birational morphism for which $E(f)$ is exceptional.
For any $l\in \NN,$
any section of $\Omega^s(Y,\log D(\fd))^{\otimes l}$
gives rise to a section of $\Omega^s(X,\log A)^{\otimes l}$
outside $E(f)$ which therefore extends over $E(f)$ via
Hartog extension on $X_0$. So 
$\kappa(L_X)\ge\kappa(L_Y)$. Conversely, there is an
open subset $U$ with codimension two complement in $Y$ 
where $f|_{f^{-1}(U)}$ is 
smooth as a log-morphism, i.e., where $f_*(T(X,-\log A))$
has maximal rank in $T(Y,-\log D(\fd))$, since 
$D(\fd)$ is reduced. Hence, $L_X| _{f^{-1}(U)}=f^*(L_Y|_U)$
so that $\kappa(L_X)\le\kappa(L_Y)$ via Hartog extension on $Y$.
$\BOX$\\

\begin{prop}\label{MR} Let $\fd: X^\bd\ra Y$, $g:Y\dra Z$ 
be fibrations and $B=\bd Y=D(f^\bd)$. Let $h^\bd=\gd\circ\fd$, 
where $\gd: Y\m B \dra Z$ is the orbifold fibration that arises. 
Then there is an inclusion
\begin{equation}\label{sat2}
{f}^*L_{\gd}^N\hra L_{h^\bd}^N\ 
\end{equation}
outside $f^{-1}(U)\cup E(f)$ where $U$ is an open subset of
$Y$ with codimension two complement
and $N$ is a positive multiple of all relevant multiplicities.
Moreover, this inclusion is an isomorphism on an open subset
of $X$ that surjects to the complement of a codimension two
subset of $Y$. Thus, there is an inclusion 
$$H^0(L_\hd^l)\hra H^0(L_\gd^l)$$
for all $l>0$ divisible by the relevant multiplicities and so
$$\kappa(Z, h^\bd)\leq \kappa(Z,\gd):=\kappa(L_\gd).$$
In addition, in the case $X^\bd$ is a log-manifold,
one can replace $ L_{h^\bd}^N$ by $L(Z,\gd;\fd)^N$ in (\ref{sat2})
without affecting any conclusion. Also,
if $\fd$ is admissible, then
$$H^0(L_\hd^l)= H^0(L_\gd^l)\ \ \text{and }\ \ 
\kappa(Z,\gd)=\kappa(Z, h^\bd):=\kappa(L_{h^\bd})$$
and (in the case $X^\bd$ is a log-manifold)
$$\kappa(Z,\gd)= \kappa'(Z,\gd;\fd):=\kappa(L(Z,h^\bd)).$$
In particular, if $\yd$ is a log-manifold
(for $f$ the identity morphism), we obtain 
$$\kappa(Z,\gd)=\kappa'(Z,\gd):=\kappa(L(Z,\gd)).$$
\end{prop}

\noindent
{\bf Proof:} By applying the Hartog
extension theorem, the last four statements follow directly from the
claim about (\ref{sat2}), which we now prove. (In fact, the last
two statements also follow directly from Lemma~\ref{fu}.)

Replacing $Z$ by a modification, we may assume that all prime
divisors in $Y$ not dominating $Z$ maps onto divisors in 
$Z$. By this, we may take an  open set $U\subset Y$
with  codimension two complement on which $g$ is
equidimensional (onto the new $Z$) and where 
$g$ is smooth when restricted to each prime divisor.

Now, outside 
the discriminant locus $\Delta(f)$ of $f$,
the first inclusion in (\ref{sat2}) is an 
isomorphism. So it remains to look at a divisorial
component of $\Delta(f)$.
Let $R$ be such a component. 
If $g(R)=Z$, then $L_g=g^*K_Z$ on a neighborhood
of $R\cap U$ in $U$ by definition and 
so $L_h=L_{h^\bd}=f^*L_{g^\bd}=h^*K_Z$ 
($=L(Z,h^\bd)$ in the case $X^\bd$ is a 
log-manifold) there. If not,
let $S=g(R)\in \dv'(Z)$. 
{}From now on, we will only consider our problem over
a neighborhood of $R\cap U$ in $U$ and outside $E(f)$
in $X$. There, we
have $g^*(S)=mR$ for some $m\in \NN$. Let
$f^*R=\sum m_iD_i$ for $D_i\in \dv'(X_1)$ 
and set $m(D_i)=m_i\,m(D_i\cap \bd X^\bd)$ and 
$m'=\min_i m(D_i)=m(R\cap \bd Y^\bd)$.
Now, by a direct local computation, we have
(compare \cite{Re}) that 
\begin{equation}\label{a}
(L_{h^\bd})=h^* K_{Z} + h^*S-(h^*S)_{red}+\bd X^\bd
   =h^* K_{Z} + \sum_i (m\cdot m_i -{m_i/m(D_i)})D_i
\end{equation}
as divisors and similarly that 
$(L_{g})=g^* K_{Z}+(m-1)R$.
Hence
\begin{equation}\label{b}
(f^*L_{g^\bd})=(f^*L_{g})+(1-{1/m'})f^*R=
h^* K_{Z} + \sum_i (m-{1/m'})m_iD_i.
\end{equation}
We thus see that 
$(f^*L_{g^\bd})\le (L_{h_1^\bd})$
as a $\QQ$-divisor above this neighborhood
of $R$ and outside $E(f)$ with equality achieved above
every point of $R\cap U$ (on $D_{i_0}$
where $m(D_{i_0})=m'$). This coupled with 
the natural isomorphism of $L_h$
with $f^*L_g$ outside $\Delta(f)$ in $U$
establishes our claim about (\ref{sat2}),
and so the lemma follows. $\BOX$\\

The following are straightforward generalizations of
Lemma~\ref{A} and Lemma~\ref{B}
to orbifold fibrations respectively.

\begin{lem}[Birational Invariance]\label{AA} Consider a 
commutative diagram of maps:
\begin{equation}\label{dm11}
\begin{CD}
X_1^\bd @>v^\bd>> X^\bd\\
@V{f_1^\bd}VV @VV{\fd}V\\
Y_1^\bd @>u^\bd>> Y^\bd\\
@V{g_1^\bd}VV @VV{\gd}V\\
Z_1 @>w>> Z
\end{CD}
\end{equation}
where $f_1^\bd,\ \fd$ are orbifold
fibrations with orbifold bases
$Y_1^\bd=Y_1\m D(f_1^\bd)$ and  $\yd=Y\m D(\fd)$, 
$g_1^\bd$ and $\gd$ are fibrations, $v^\bd$ is a 
morphism and $w$ is generically finite and surjective.
If $\fd$ is admissible (so that $u^\bd$ is a premorphism), 
then we have
$$H^0(L_\gd^l)\hra H^0(L_{g_1^\bd}^l)\ \ \text{and so}\ \ 
\kappa(Z_1,g_1^\bd)\ge \kappa(Z,g^\bd)$$
for all $l>0$ divisible by the relevant multiplicities.
If $u^\bd$ is a morphism,
then this remains true without the upper part of
the diagram. 
If moreover $u^\bd$ and thus $w$ are birational, then
we have (without the upper part of the diagram),
$$H^0(L_\gd^l)=H^0(L_{g_1^\bd}^l)\ \ \text{and so}\ \ 
\kappa(Z_1,g_1^\bd)= \kappa(Z,g^\bd).$$
If $v^\bd$ is birational and onto and $u$ is generically
finite (hence birational), then the same conclusions hold
(with the upper part of the diagram). Moreover, all
the results remain valid if $g$, $g_1$ and $w$ are
only assumed to be meromorphic.
\end{lem}

\noindent
{\bf Proof of Lemma~\ref{AA}}: The last part of the lemma
without the last statement and
without assuming that the upper part of the diagram exists,
follows the same way as that of Lemma~\ref{A}.
If we assume that 
$\fd$ is admissible, then 
the proof of the first part and the second last statement is
also the same as that of Lemma~\ref{A}
using  Lemma~\ref{A*} and Proposition~\ref{MR}. 
Hence the result follows.
$\BOX$\\

Note in particular that a prebirational morphism 
preserves orbifold Kodaira dimension 
(by taking $g$ and $g_1$ to be the identity map).

\begin{lemma}\label{BB}
Let the second column of diagram~\ref{dm11} 
be given as in the previous lemma with $X^\bd$ 
smooth. Namely, let 
$\fd: X^\bd\ra Y$, $g:Y\ra Z$ be fibrations with $X^\bd$ 
smooth, $\yd:=Y\m D(f^\bd)$.  
Then one can complete this
column to a diagram as in diagram~\ref{dm11} 
with $u^\bd,v^\bd$ birational and
onto so that the fibrations 
$f_1^\bd,g_1^\bd$ and $h_1^\bd=g_1^\bd\circ f_1^\bd$ 
are all well adapted
and all divisors in $X$ (resp. $Y$) not mapping onto $Y$ 
(resp. $Z$) push forward to divisors in $Y_1$ (resp. $Z_1$). 
$f_1^\bd$ and $h_1^\bd$ are birationally equivalent to
$\fd$ and $\hd=g^\bd\circ f^\bd$ resp., and if $\fd$ is
admissible, $g_1^\bd$ is also birationally equivalent to $\gd$.
Hence if $\fd$ is admissible, then there are birationally 
equivalent fibrations to $\fd$, $\gd$ and $\hd$ such that
equality is attained in Lemma~\ref{comp}. Moreover, there is an
open subset $U$ of $Y$ with  codimension
two complement on which $f$ and $f_1$ are isomorphic
under this correspondence.

If further $g^\bd\circ \fd: X^\bd \ra Z$ is admissible, then
$$\kappa(Z,\gd;\fd)= \kappa(Z\m D(\gd)).$$
\end{lemma}

\noindent
{\bf Proof of Lemma~\ref{BB}:} 
The last statement follows from Lemma~\ref{B}.
So we only need to construct $f_1$ and $g_1$. For this 
purpose, we may first assume that $\fd$ is well adapted
by blowing up and replacing it with a birationally
equivalent one using Lemma~\ref{none}. Then
we take a flattening of $f$ and let the
resulting flat fibration be $f':X'\ra Y'$ and the base
change map be $u':Y'\ra Y$. Then we take
a flattening of $g\circ u'$ and let the
resulting flat fibration be $g_0:Y_0\ra Z_0$ with
map $u_0:Y_0\ra Y'$. Then we base change $f'$ by
$u_0$ to give a flat fibration $f_0:X_0\ra Y_0$.
We now normalize $X_0$ and $Y_0$ and denote the
resulting equidimensional fibrations by the same
symbols $f_0, g_0$ and $g_0\circ f_0$. Let 
$u_1: Y_1\ra Y_0$
be a resolution of singularities of $Y_0$. Since $\fd$
is well adapted, the resulting birational morphism from
$Y_1$ to $Y$ is an orbifold birational morphism
by Lemma.~\ref{na}.
Then $g_1:=g_0\circ u_1$ is admissible. Now let
$f'':X''\ra Y_1$ be the base change of $f_0$ by $u_1$.
Let $v_1:X_1\ra X''$ be a resolution of singularities 
of $X''$ so that the map from $X_1$ to $X$ gives an
orbifold birational morphism. Then $f_1:=f''\circ v_1$
is admissible since $f''$ is equidimensional. Also
$g_1\circ f_1$ is admissible since it factors through
the equidimensional morphism $g_0\circ f_0$. 
Hence if we chosse the resolutions to be normal
crossing and the multiplicity of the exceptional
divisors on $X_1$ to be $\infty$, we see that
$f_1$, $g_1$ and $h_1$ are all well adapted. $\BOX$\\

We see in particular that 
$\kappa(Z,\gd)=\kappa(L_\gd)$ is a
birational invariant of orbifolds.
In this sense, we may speak of $L{(Z,\gd)}$
({\em as an equivalence class of sheaves})
as a well defined object for $\gd$ and its birational class.\\

We summarize the main result of this section in the log-case:
Given an adapted fibration $\fd: X^\bd\ra Y$ where
$X^\bd= X\m B$ is a log-manifold, there is
an injective correspondence between birational equivalence classes of 
orbifold fibrations from $\yd=Y\m D(\fd)$ and saturated subsheaves of 
$\Omega^i(X, \log B)$ which are trivial on the generic 
fibers of $f$ for $i= 1,\dots, \dim Y$. Furthermore,
the Kodaira dimension of such a subsheaf is the
Kodaira dimension of the orbifold $Z_1\m D(g_1^\bd)$
for an admissible birationally equivalent
fibration $g_1^\bd: Y^\bd_1\ra Z_1$,
which one can always attain
by flattening any orbifold fibration
corresponding to the subsheaf and then resolving the singularities
of the resulting family suitably.

\section{The birational geometry of projective orbifolds}

For simplicity of treatment, all orbifolds treated in this
chapter are standard. The general non-standard case can be
obtained by considering (in general) a morphism (or fibration)
from a standard orbifold $X^\bd$ to $\yd$ as the composition of a 
morphism (respectively fibration) from a smooth orbifold
$\te X^\bd$ with the inverse of a finite morphism from 
$\te X^\bd$ to $X^\bd$. Alternatively, if one is not
strictly interested on line sheaves of holomorphic p-forms,
on can simply regard the $L_\gd$ constructed in this section
for the nonstandard case as line sheaves on $\te Y$ which
are pullbacks of the $\QQ$-line sheaves on $Y$ constructed in
the last section and also denoted by $L_\gd$ there.  
All the results of this section 
are valid in the more general category this way.

In the projective category, we have the following basic
technical covering lemma of Kawamata 
(\cite[Theorem 17, 19]{Ka}, \cite[Theorem 1-1-1]{KMM})
for reducing a question on a smooth
orbifold to a corresponding question concerning ordinary 
manifold. In general, the technical lemma is part of 
some fundamental results in the theory of
log-schemes mentioned in \cite{LM1}, 
which we will discuss elsewhere.

\begin{lem}[Kawamata Branched Covering Theorem]\label{K}
Let $\gamma': Y'\ra Y$ be a finite surjective morphism from a 
normal variety $Y'$ to a smooth variety $Y$ and 
$D$ a reduced divisor on $Y$. Suppose 
$D(\ga')_{red}+ D$ is a simple normal crossing
divisor and $D= \sum_i D_i$, 
where $D_i\in \dv'(Y)$ for all $i$. Given $n_i\in \NN$,
there exist a finite Galois covering $\ga: \te Y\ra Y$
factoring through $\ga'$ with $\te Y$ smooth.  Moreover, 
$\ga^{-1}(D)\cup R$ is
simple normal crossing where $R=(\det (d\ga))_{red}$.
\end{lem}

\noindent
{\bf Proof:} The statement above is slightly more general
than that stated in \cite{Ka} but the proof is exactly
the same. $\BOX$\\

\begin{lem}\label{one}
Let $Y^\bd=Y\m B$ be a smooth orbifold. Then 
\begin{itemize}
\item[\rm (I)] There exist a finite Galois covering $\ga: \te Y \ra Y$
whose ramification divisor R on $\te Y$ is supported on a simple
normal crossing divisor containing $\ga^*(B_{red})$ and 
$\ga^*(B)$ is a Cartier divisor on $\te Y$.

\item[\rm (II)] There is a unique locally free subsheaf $\Omega_\yd$
of $\Omega_{\te Y}$ with determinant sheaf $\ga^*K_\yd$ such that
$\Omega_\yd$ is the natural subsheaf 
$\ga^*\Omega_Y\hra \Omega_{\te Y}$ outside $R_{red}$. 
In particular, $\Omega_\yd$ is invariant under the Galois 
group $G$ of $\ga$.
\end{itemize}
\end{lem}

{\bf Proof:} (I) follows from the previous lemma. 
(II) follows from (I) by the 
following observation: Each $p\in Y$ has a coordinate 
neighborhood $U$ with coordinate $(y_1,\dots,y_m)$
on which we may write $\ga=\delta\circ \alpha$ for
$$\delta: \Delta^r \twoheadrightarrow U\subset Y;\ \ 
(x_1,\dots,x_r) \mapsto (x_1^{M_1},\dots,x_r^{M_r}),$$
where $B|_U$ is supported on the coordinate axis,
$M_i=m(D_i, Y^\bd)$ or $1$ according to whether $m(D_i, Y^\bd)$ 
is finite or not and where we have arranged the indices  
so that $D_i|_U=(y_i)$ for all $i$ in
$B_{red}=\sum_i D_i$. We then define
$\Omega_\yd|_{\ga^{-1}(U)}=\alpha^*\Omega_{\Delta^r}(\log D)
\hra \Omega_{\te Y}|_U$ where 
$D=\sum\{(y_i)\ |\ m(D_i, Y^\bd)=\infty\}$. Hence,
we may conclude the proof by noting that this definition is
invariant under the Galois group of $f$ and that
$\det\Omega_\yd = \ga^*K_\yd$ by construction. $\BOX$\\

We call the covering $\ga$ so constructed a {\bf $\bd$-covering}
of $Y\m B$. It can also be written as a finite surjective 
orbifold ``\'etale'' morphism: $$\ga:  \te Y \ra \yd.$$
In general, a finite surjective morphism 
$\ga^\bd:  \te Y^\bd \ra \yd,$ is called {\bf \'etale} if
the induced structure from $\yd$ on $\te Y$, namely
$$
\ga^*(\bd \yd)-(\det(d\ga))
$$ 
coincides
with that of $\te Y^\bd$.

\begin{lem} \label{two}
Let $u: Y_1\ra Y$ be a generically finite surjective morphism with 
$Y_1^\bd$ and $Y^\bd$ smooth orbifolds. 
Let $\ga: \te Y\ra Y$ be a $\bd$-covering
of $Y\m B$ (as given in Lemma~\ref{one}).
Then there is a finite Galois covering $\ga_1$ as given by
Lemma~\ref{one} for the orbifold $Y_1^\bd$ 
where $u\circ \ga_1=\ga\circ \te u$ for a generically finite 
morphism $\te u: \te Y_1\ra \te Y$.  Furthermore, 
$u^\bd$ is a morphism if and only if there is a natural
inclusion $\iota : \te u^*\Omega_{Y^\bd} \hra \Omega_{Y_1^\bd}$
which is an isomorphism on $\ga_1^{-1}(W)$ for the open set $W$
on which $u$ is an isomorphism. In particular, if $u^\bd$ is orbifold
birational, then $\iota$ is an isomorphism outside the exceptional 
divisor of $\te u$.
\end{lem}

\noindent
{\bf Proof:} To construct $\ga_1$, we first take the normalization
$Y'_1$ of the fiber product $Y_1\times_Y \te Y$. Then Lemma~\ref{K}
guarantees the existence of our required $\ga$ factoring through
$\ga_1$. If we have an inclusion $\iota$ as above, then it is
easy to see that $u^\bd$ is a morphism and we leave this 
verification to the reader. Suppose now that $u^\bd$ is a morphism.
To see the inclusion $\iota$, we first observe that the
natural inclusion $j: \te u^*\Omega_{Y^\bd} \hra 
\Omega_{Y_1^{|\bd|}}=\Omega(Y_1, \log D)$ is
an isomorphism outside $D=\ga_1^{-1}((\bd Y_1)_{red})$.
Now $j$ factors through the natural inclusion 
$\Omega_{Y_1^\bd}\hra \Omega(Y_1, \log D)$ outside the
exceptional divisor $E$ of $\te u$ by virtue of 
the Hartog extension theorem and the condition 
$K_{Y_1^\bd}\ge u^* K_{Y_1^\bd}$, which on $\te Y_1$ is 
the condition 
\begin{equation}\label{in2}
\ga_1^*K_{Y_1^\bd}\ge \te u^*\ga^* K_{Y_1^\bd}.
\end{equation}
To show the desired inclusion $\iota$ on $E$, it is
equivalent to showing that the holomorphic map
$\te u^*\Omega_{Y^\bd} \hra \Omega_{Y_1^{|\bd|}}$ of
vector bundles factors through $\Omega_{Y_1^\bd}$.
Since $\te u^*\Omega_{Y^\bd}$ is trivial on $E$, we only 
have to show this for each vector $v$ in the vector bundle
$\Omega_{Y^\bd}$ restricted to $E$.
But this follows from the fact
that one can find a map $g:Y\ra C$ where $C$ is a curve
with $v\in L_{g\circ \te\ga}\in \Omega_{Y^\bd}$ and that 
$\te u^*L_{g\circ \te\ga}$ includes in $\Omega_{Y_1^\bd}$
as $u^\bd$ is a morphism.
If further $u^\bd$ is birational, then (\ref{in2}) is
an equality outside the exceptional divisor of $\te u$
forcing $\iota$ to be an equality there as well. $\BOX$\\

It follows in particular that 
given two coverings $\ga_1$ and $\ga_2$ both as 
prescribed by Lemma~\ref{one} for a smooth orbifold $\yd$,
we can create a third such covering $\ga$ factoring
through $\ga_i$ and on which the $\Omega^p_\yd$ 
is naturally isomorphic to the pull back of the
$\Omega^p_\yd$ defined by $\ga_i$ for each $i$.
It is in this sense that $\Omega^s_\yd$ is well
defined independent of the covering $\ga$ that one
can choose to satisfy  Lemma~\ref{one}.\\

The same proof as above gives us the following
more standard characterization of an orbifold morphism
for any map between manifolds with orbifold structures.

\begin{prop}\label{three}
Let $u: Y_1\ra Y$ be a  morphism with 
$Y_1^\bd=Y_1\m B_1$ and $Y^\bd= Y\m B$ smooth orbifolds 
such that $u^{-1}(B_{red})\subset D:=(B_1)_{red}$. Then
there exists a commutative
diagram
\begin{equation}\label{dm4}
\begin{CD}
\te Y_1 @>{\te u}>> \te Y\\
@V{\ga_1}VV @VV{\ga}V\\
Y_1^\bd @>u>> Y^\bd\\
\end{CD}
\end{equation}
where ${\ga_1}$ and $\ga$ are $\bd$-coverings
and $\te u$ is a morphism. Moreover, setting 
$\te D=\ga_1^{-1}(D)$, $u$ gives rise to
a morphism $u^\bd: Y_1^\bd\ra Y^\bd$  if and only if
the natural morphism
$j: \te u^*\Omega_{Y^\bd} \ra \Omega(Y_1, \log \te D)$
factors through the inclusion 
$\Omega_{Y_1^\bd}\hra \Omega(Y_1, \log \te D).\ \ \ \BOX$
\end{prop}

\begin{defn}  Given a fibration
$\gd: \yd=Y\m B \ra Z$ with $\ga:\te Y\ra Y$ as given in
Lemma~\ref{one}, we define 
$L_{(Z,g^\bd)}$ (or $L_\gd$ for short) to be the saturation 
of $\ga^*g^*K_Z$ in $\Omega^s_\yd$ where $s=\dim Z$. 
We set $$\kappa(Z,\gd)= \kappa(L_\gd)$$ and call it the Kodaira
dimension of the orbifold $(Z,g^\bd)$ defined by $\gd$.\\
\end{defn}

Note that we have just given a new definition of
$\kappa(Z,\gd)$ and of $L_\gd$ and that this
$L_\gd$ depends on $\gamma$. However,
our main object of interest $\kappa(Z,\gd)$ 
is the same as that of the old definition 
because of Lemma~\ref{B'},
which therefore allows us to keep this ambiguity in
$L_\gd$ without harm. Note also that the 
$L_\gd$ defined previoiusly
was only a $\QQ$-line sheaf but its pullback
to $\te Y$ is a line sheaf that we also denote
by $L_\gd$, by abuse of notation.

\begin{lem}[Birational Invariance]\label{A'} Consider the 
diagram of morphisms
\begin{equation}\label{dm1'}
\begin{CD}
\te Y_1 @>{\te u}>> \te Y\\
@V{\ga_1}VV @VV{\ga}V\\
Y_1^\bd @>u^\bd>> Y^\bd\\
@V{g_1^\bd}VV @VV{\gd}V\\
Z_1 @>w>> Z
\end{CD}
\end{equation}
where $g_1^\bd$ and $g$ are fibrations from smooth orbifolds
$Y_1^\bd$ and $\yd$ respectively. Assume that the lower square of the
diagram commutes and that $w$ is generically finite and surjective.  
Then there exist $\bd$-coverings $\ga_1$ for $Y_1^\bd$
and $\ga$ for $Y^\bd$ and a morphism
$\te u:\te Y_1\ra \te Y$ as given by Lemma~\ref{three}
making the above diagram commutative and
as a result,
$$\kappa(Z_1,g_1^\bd)\ge \kappa(Z,g^\bd).$$

If moreover $u^\bd$ and thus $w$ are birational, then
$$\kappa(Z_1,g_1^\bd)= \kappa(Z,g^\bd).$$
Again, all the results hold if $g$, $g_1$ and $w$ are
only assumed to be meromorphic.
\end{lem}

\noindent
{\bf Proof:} We restrict to the case where $\te u$ is 
surjective for simplicity. 
The existence of the upper diagram making the entire
diagram commutative is given by Lemma~\ref{three}. 
If $w$ is generically finite
and surjective, then we have 
\begin{equation}\label{in3}
L_{g_1^\bd}^m \supset \te u^*L_\gd^m,
\end{equation}
as sheaves by the commutativity of the above diagram. By the
surjectivity of $\te u$, we have an injection
$H^0(\te Y, L_\gd^m)\hra H^0(\te Y_1, \te u^*L_\gd^m)$ and so the first
result follows.  If $u^\bd$ is birational, then the 
sheaf inclusion above is an equality on
$\te Y_1\setminus E$ where $E$ is the 
exceptional divisor of $\te u$. Hence 
$$H^0(\te Y_1,L_{g_1^\bd}^{Nm}) 
\ge H^0(\te Y,L_{\gd}^{Nm}) \ge H^0(\te Y_1,L_{g_1^\bd}^{m}),$$
where $N=\deg \te u$ and the last inequality is obtained by 
the following procedure: We may assume (by taking a finite
covering if necessary) that $\te u$ is Galois for 
simplicity. Given a
section $s$ of $L_{g_1^\bd}^m$,
we first identify $s$ with a section 
of $\te u^*L_\gd^m$ on $\te Y_1\setminus E$.
Then we take the $\te u$-trace of $s$ via 
$$s\mapsto \prod \{\ g(s)\ |\ g\in G\ \}$$
where $G$ is the Galois group of $\te u$ and
apply the Hartog extension theorem on $\te Y$ as $v(E)$ is 
codimension two or higher in $\te Y$. (We remark here
that we can define the trace of sections of the pull back of a line sheaf
under any finite map similarly.) $\BOX$\\

\begin{lemma}\label{B'}
Any fibration  $\gd: \yd \ra Z$ has an admissible
birationally equivalent representative. Any admissible 
fibration $\gd :\yd\ra Z$ satisfies 
$$\kappa(Z\m D(\yd))=\kappa(Z ,\gd):=\kappa(\te Y, L_\gd).$$
\end{lemma}

\noindent
{\bf Proof:} The proof is similar to the one above and that
of Lemma~\ref{B} and thus left as an exercise. \BOX\\

This lemma shows that 
$\kappa(Z,\gd)= \kappa(L_\gd)$ is not only independent of 
the choice of the finite covering $\ga: \te Y\ra \yd$
used to define $L_\gd$ but it is also independent of the
birational model of the orbifold fibration $\gd$. As in
the previous section, we will thus think of $L_\gd$ 
({\em as an equivalence class of sheaves})
as a well defined object for $\gd$.

We summarize a main part of this section, in parallel
with the previous section, as follows. 
Given a smooth projective orbifold $\yd$,
there is an injective map from 
the birational equivalence classes of orbifold
fibrations and saturated line subsheaves of $\Omega^i_\yd $,
or equivalently, saturated subsheaves of $\Omega^i_\yd $
of rank $i$ for $i=1,\dots, r=\dim Y$. Furthermore, their
Kodaira dimensions are the Kodaira dimensions of the ``true''
orbifold bases of the orbifold fibrations.

\section{Bogomolov sheaves and general type fibrations}

We have established in the previous two sections 
some elementary facts concerning (orbifold) fibrations
from an admissible orbifold $Y^\bd$, in particular, the fact that
that their birationally equivalent classes are in one
to one correspondence with certain saturated line subsheaves
of exterior powers of the orbifold cotangent sheaf $\Omega_\yd$
(or $\Omega_{X^\bd}$ in the non-projective category, where 
$X^\bd$ is a log-manifold of which $Y^\bd$ is the orbifold base).
Hence such subsheaves are certainly worthy of study in 
birational geometry. A distinguished class of such is played
by such a subsheaf $L\hra \Omega^p_\yd$ 
(resp. $L\hra \Omega^p_{X^\bd}$) of {\bf general~type}, 
which we define by the condition $$\kappa(L)=p.$$
We will allow $p=0$ for which $\Omega_\yd^p={\mathcal O}_Y$.
When $Y^\bd$ is a projective manifold, such subsheaves were
studied by Bogomolov in \cite{Bo}. He also gave there an 
alternative characterization for them via a generalization
of the lemma of Castelnuovo-De Franchis. To state it, we
first give an important definition.

\begin{defn} A meromorphic fibration $f: X\dra Y$ is called
{\bf almost holomorphic} if $f$ is holomorphic above a Zariski
open subset of $Y$.
\end{defn}

\begin{theorem}[Bogomolov]\label{CD}
Let $L$ be a saturated line subsheaf of $\Omega^p_X$
with $X$ K\"ahler, or more generally of Fujiki class 
$\mathcal C$ (denoted by $X\in {\mathcal C}$). Then
\begin{itemize}
\item[\rm (I)] $\kappa(L)\le p$
\item[\rm (II)] If $\kappa(L)=p$, 
then the Iitaka fibration $I_L$ of $L$ defines 
an almost holomorphic fibration from $X$
to a projective base $B$ of dimension $p$ and
$I_L^*K_B$ saturates to $L$ in $\Omega_X^p$. 
In particular, $L=L_{I_L}$.
\end{itemize}
\end{theorem}

Bogomolov in \cite{Bo} only stated the theorem in 
the projective case and that, as far as we know, the
above generalization (to the K\"ahler case) is due to
Campana in \cite{Ca01}. But the proof of this more
general case is 
essentially the same and well known among experts
including, of course, Bogomolov (see \cite{ViVan}
for example).  In view of its importance, we 
indicate a derivation of a slightly more general
version.\\

\begin{prop}\label{logCD}
Let $L$ be a saturated line subsheaf of $\Omega^p(X,\log A)$
with $X$ K\"ahler, or more generally of Fujiki class 
$\mathcal C$ (denoted by $X\in {\mathcal C}$),
and $A$ a normal crossing divisor in $X$. Then
\begin{itemize}
\item[\rm (I)] $\kappa(L)\le p.$
\item[\rm (II)] If $\kappa(L)=p$, 
then the Iitaka fibration $I_L$ of $L$ defines 
an almost holomorphic fibration of general type from $X$
to a projective base $B$ of dimension $p$ and
$I_L^*K_B$ saturates to $L$ in $\Omega^p(X,\log A)$. 
In particular, $L=L_{I_L^\bd}$.
\end{itemize}
\end{prop}

\noindent
{\bf Proof:}
Assume $\kappa(L)=r$. We then have, for some 
$m$, a fibration $f: X \ra Y$ birationally equivalent 
to $I_L$  defined by $H^0(X, L^m)$ such that
$r=\dim Y$. By a well-known covering trick of taking
roots of sections (see \cite{BlG}), we can  
reduce the problem to the case $m=1$ as was done
by Bogomolov in \cite{Bo} (see also \cite{ViVan}).
We now assume $m=1$. We may then take $r+1$
(logarithmic) $p$-forms $\omega_0,\dots, \omega_s$ 
in $L$ such that the $r$ meromorphic functions 
$y_i=\omega_i/\omega_0$, $i=1,\dots,r$ are then pull backs of 
meromorphic functions on $Y$ which serve as coordinates for $Y$
generically. Since $X$ is K\"ahler (or of Fujiki class 
$\mathcal C$), logarithmic forms are closed by
Hodge theory (or by \cite{Nog}) and so 
$$0=d\omega_i=(dy_i)\wedge \omega_0$$ for all $i$.
Since the $dy_i$'s are pointwise linearly independent
over an open subset of $Y$, we see that $r\le p$.

When $r=p$, linear algebra shows that 
$dy_1\wedge\dots\wedge dy_r$ is 
generically a nontrivial multiple of $s_0$.
This forces $f^*K_Y\hra \Omega_X^p(X,\log A)$ to factor
through $L\hra \Omega_X^p(X,\log A)$ genercally so that
$L=L_{I_L^\bd}$. Hence $I_L^\bd:X\m A \dra Y$ 
is a meromorphic orbifold fibration of 
general type. Now Proposition~\ref{almo} 
concludes the proof.
$\BOX$\\

The same proof gives us the following generalization:

\begin{prop}\label{orbCD} Let $\yd$ be an orbifold with
a $\bd$-covering $\ga:X\ra \yd$ (respectively an admissible 
fibration $\fd$ from a log-manifold $\xd$ with 
$\yd=Y\m D(\fd)$). Let $L$ be a subsheaf of 
$\Omega^p_\yd$ (respectively $\Omega^p_\xd$) invariant
under the Galois group of $\ga$ (respectively trivial
on the general fibers of $f$). If $X\in {\mathcal C}$,
then
\begin{itemize}
\item[\rm (I)] $\kappa(L)\le p$
\item[\rm (II)] If $\kappa(L)=p$, 
then $\ga$ (resp. $\fd$) dominates
the Iitaka fibration $I_L: X\dra B$ of $L$ giving 
an almost holomorphic fibration from $Y$
to a projective base $B$ of dimension $p$ and
$I_L^*K_B$ saturates to $L$ in $\Omega_\yd^p$
(resp. $\Omega^p_\xd$). In particular, $L=L_{I_L^\bd}$.
\end{itemize}
\end{prop}

We remark that these Bogomolov type theorems
fail for more general compact complex manifolds.
In view of the above theorems, we can make the following 
definition (compare also \cite{Ca}).

\begin{defn} A meromorphic fibration $g^\bd:\yd\dra Z$
is said to be of {\bf general type} if 
$$\kappa(Z,\gd):=\kappa (L_\gd)=\dim Z.$$
Given an orbifold $\yd$ which is smooth and projective
(resp. which is the orbifold base of an admissible fibration 
from a log-manifold or from a smooth and projective
orbifold $X^\bd$), 
a saturated line subsheaf of $\Omega_\yd^p$ 
(resp. of $\Omega_{X^\bd}^p$) for $p\in \{1,\dots,\dim X\}$
corresponding to a general type fibration is called
a {\bf Bogomolov sheaf} (of dimension $p$ on $\yd$).
A meromorphic map from an orbifold is said to be 
of general type if its Stein factorization is so.
\end{defn}

Hence, if $X$ is of Fujiki class $\mathcal C$, the above
theorem of Bogomolov 
states that saturated line subsheaves $L$ of 
$\Omega_X^p$ with $\kappa(L)=p$ for $p=1,\dots, \dim X$
are precisely the Bogomolov sheaves on $X$, i.e., they
are the sheaves corresponding to meromorphic 
general type fibrations. And this continues to 
hold for admissible orbifolds by Proposition~\ref{orbCD}
above. On the other hand, we note that for non-admissible $\xd$, 
our notion of general type above should be used with care
since it can no longer be made meaningful as 
a birational invariant and indeed is not a
birational invariant in general. \\

The following basic lemma is a direct consequence of Lemma~\ref{AA}.

\begin{lem}\label{gt*}
Consider a fibration $g^\bd: \yd\dra Z$ and a
premorphism $\phi^\bd: Y_1^\bd\ra \yd$ such that $\psi=g\circ \phi$
is surjective. If $g^\bd$ is of general type,  
then so is $\psi^\bd$ (or more properly, the birational
class of its Stein factorization). $\BOX$\\
\end{lem}

We now establish an important technical lemma for this paper, 
Proposition~\ref{tech0}.

\begin{lem}\label{tech0'}
Consider the following commutative diagram 
\begin{equation}\label{dmT'}
\begin{CD}
X^\bd_t @>i^\bd>> X^\bd @. \\
@V{f_t^\bd}VV @VV{f^\bd\ \ \ \, \searrow\, h^\bd}V @.\\
Y_t @>j>> Y @>g>> Z\\
@. @VVuV @.\\
@. T @.
\end{CD} 
\end{equation}
where $\fd$, $g$, $h^\bd=g\circ \fd$ and $u$ are fibrations, 
$i$ and $j$ are the inclusion 
of the general fibers $Y_t:=u^{-1}(t)$ and $X_t:=f^{-1}(Y_t)$ over $T$
(by which we mean that $t$ is a representative point in
the complement of some countably many closed proper subsets of $T$)
and $i^\bd$ is an orbifold inclusion.
Let $h=g\circ f$. Assume that $g\circ j$ is generically finite.
Then we have with $p=\dim Z$ and
$q=\dim Y_t$ that
$$\kappa(Z, h^\bd)- (p-q) \leq \kappa(i^*L_{(Z,h^\bd)})
\le \kappa(L_{f_t^\bd})=\kappa(Y_t, f_t^\bd).$$
In particular, if $h^\bd$ is of general type, then so is $f_t^\bd$.

All the results remain valid if one only assumes that
$h$, $g$ $f$ and hence $f_t$ are meromorphic maps 
as long as $u\circ f$ is holomorphic. 
\end{lem}

\noindent
{\bf Proof:} We first ignore the last statement. We note that 
$X_t^\bd= X_t\m i^*A$ for $A=\bd X^\bd$, true for general $t\in T$. 
We have the following commutative diagram of sheaf morphisms
\[
\begin{array}{lccc}
i^*h^*K_Z\, \hra \, i^*L_{(Z,h)}\, \hra\, i^*\Omega_{X}^p\ \  \ra 
&\Omega^q_{X_t} \otimes \wedge^{p-q}N_{X_t}^\vee 
&\ra 
&\Omega^q_{X_t} \\
\ \ \ \ \ \| & \uparrow & &\uparrow\\
f_t^*j^*g^*K_Z\ \ \ \ \hra\ \ \ \ f_t^*j^*\Omega_Y^p\ \  \ \lra
&f_t^*K_{Y_t} \otimes \wedge^{p-q}N_{X_t}^\vee   &\ra
&f_t^*K_{Y_t}
\end{array}
\]
where the vertical arrows are inclusions,
$N^\vee_{X_t}={\mathcal O}_{X_t}^{\oplus(p-q)}$ 
is the conormal bundle of $X_t$ and where the two right arrows
on the very right are projections obtained by choosing
a section $s$ of $\wedge^{p-q}N^\vee_{X_t}$ so that the 
image of  $i^*h^*K_Z$ in $f_t^*K_{Y_t}$, or equivalently,
that the image of $j^*g^*K_Z$ in $K_{Y_t}$ is nontrivial
(such an $s$ exists for general $t$ by construction).
Hence the saturations of $f_t^*K_{Y_t}$ and of $i^*L_{(Z,h)}$
in $\Omega^q_{X_t}$ are the same, namely $L_{(Y_t, f_t)}$.
The above argument still works under the more general situation
given by the last statement of the lemma since we only 
need to establish the equality of $f_t^*K_{Y_t}$ and 
$i^*L_{(Z,h)}$ generically on $X_t$ and so we only need to
work on the open set $U$ on which $f$ holomorphic as
$U$ and therefore its restriction to $X_t$ has 
at most codimension two complement.  Hence,
the right hand side of the above inequality follows from 
the easily deduced fact that 
$i^*(A\cap h)=(i^*A)\cap h\circ i=(i^*A)\cap f_t.$

For the  rest of the inequality and hence the lemma, we 
generalize the proof of the easy addition law: 
We may assume that 
$\kappa(Z, h^\bd):=\kappa(L_\hd)\geq 0,$ since
otherwise the result is clear. Then 
$\kappa(i^*L_\hd)\geq 0$ as well and so we may
also assume that $\kappa(L_\hd)> p-q\geq 0$.
We may then take the Iitaka fibration 
of $L_\hd$, which we may take to be
$\mu:X\dra X'$ with $X'$ projective. Since
$L_\hd$ is trivial on the general fibers of $h$, $\mu$ 
factors through a meromorphic fibration $\gamma: Z\dra X'$. 
The image $X'_t$ of $X_t$ in $X'$,
which is also the image of $Y_t$ in $X'$,
is positive dimensional as $\dim X'>p-q$. 
Let the codimension of $X'_t$ in $X'$ be $r$. 
We then have by construction that 
$\kappa(i^*L_\hd)\geq \dim X'_t$ and $p-q\geq r$.
It follows that $$\kappa(L_\hd)=\dim X'=\dim X'_t
+r\leq \kappa(i^*L_\hd)+ p-q$$ as required.
$\BOX$\\

We observe that the above proof does not 
use the fact that $Z$ is smooth. This is because 
we can replace $K_Z$ above by $K_U$ for a smooth 
open subset of $Z$ since we only need to consider
rational sections of $K_U$ because the saturation
considered above in $\Omega_{X_t}^q$ is uniquely 
determined by their pullbacks. For the same reason,
nor do we need $f,\ g$ and $h$ to be holomorphic.

\begin{prop}\label{tech0}
Consider the following commutative diagram 
\begin{equation}\label{dmT2'}
\begin{CD}
X^\bd_t @>{i^\bd}>> X^\bd @>v^\bd>> X_0^\bd \\
@V{f_t^\bd}VV @VV{f^\bd}V @VV{f_0^\bd}V\\
V_t^\bd @>{k^\bd}>> V^\bd @>{u_0^\bd}>> Y_0^\bd\\
@. @VVeV @VV{g_0^\bd}V\\
@. T @. Z_0
\end{CD} 
\end{equation}
where the vertical arrows are fibrations with $f_0^\bd$ 
admissible, $v^\bd$ a morphism,
$\bd Y_0=D(f^\bd_0)$,
$\bd V=D(\fd)$ and $\bd V_t=D(f_t^\bd)$,
$i^\bd$ and $k^\bd$ are orbifold inclusions of the fibers
at the general point $t\in T$ and $u_0$ is a
morphism. Assume that
$g_0\circ u_0$ is surjective with Stein factorization
given by $g:V\ra Z$ 
and that $g_0\circ u_0\circ k$ is nontrivial with
Stein factorization $g_t: V_t\ra Z_t$.
Then we have, with 
$h_0=g_0\circ f_0,\ h=g\circ f,\ h_t=g_t\circ f_t$
and $d_c={\rm codim}_{Z_0}\ g_0(u_0(k(V_t)))$, that
$$\kappa(L_{g_0^\bd})=\kappa(Z_0, h_0^\bd)\leq \kappa(Z,h^\bd)
\leq\kappa(Z_t, h_t^\bd)+ d_c\leq \kappa(L_{g_t^\bd}) + d_c.$$
In particular, if further $g_0^\bd$ is a fibration of 
general type, then so is $g_t^\bd: V_t^\bd\ra Z_t$. 
All results remain
valid if $g_0$ is only meromorphic.
\end{prop}

\noindent
{\bf Proof:}  
The left most inequality is a direct
application of Lemma~\ref{AA} while the above lemma
gives the rest via the following construction. 
Set $\phi'=(e,g): V\dra T\times Z$ and
let $\phi: V\dra Y$ be its Stein factorization, where
$Y$ is a smooth model of the normalization of the 
image of $\phi'$ in the function field of $V$ and
so we may assume that the map from $Y$ to
$T\times Z$ is a morphism.
Let $u:Y\ra T$ and $w:Y\ra Z$ be the projections 
to the factors, necessarily surjective, and
$\phi_t=\phi|_{V_t}:V_t\twoheadrightarrow Y_t$. Then
$w$ restricted to the general fiber $Y_t$ of
$u$ is finite by construction
and $e=u\circ \phi$ is holomorphic. Hence the 
result follows from the above lemma and the
fact that 
$\kappa(L_{\phi_t})=\kappa(L_{g_t^\bd})
\geq\kappa(L_{h_t^\bd})=\kappa(L_{\phi\circ \fd})$
via the first inequality in Proposition~\ref{MR}. $\BOX$\\

\begin{cor}\label{tech} 
Let $\gd: \yd\dra Z$ and $w:Y\ra T$ be fibrations
such that $g(Y_t)$ is nontrivial for the
general fiber $Y_t$ of $w$. If $\gd$ is of
general type, then so is its restriction 
$g_t^\bd:Y_t^\bd\ra g(Y_t)$ (or, more properly,
the Stein factorization of $g_t^\bd$) to the 
general fiber $Y_t$.\\
\end{cor}

A meromorphic fibration of general type enjoys the 
following well-known but important property 
(c.f. \cite{Re} for example)
whose proof we delay to the next section.

\begin{prop} A meromorphic fibration of general type 
from a smooth orbifold is almost holomorphic.\\
\end{prop}

We conclude this section by recalling a very useful lemma of Fujita 
(\cite{Fta}, \cite{Mori})
in the theory of Kodaira dimension:

\begin{lem} \label{fta}
Let $f:X\ra Y$ be a fibration with
$Y$ projective and $L$ a line
bundle on $X$. Then the following are equivalent:
\begin{itemize}
\item $h^0(L^mf^* H^{-1})\neq 0$ for some $m\in \NN$ and
an ample line bundle $H$ on $Y$.
\item $\kappa(L)=\kappa(L|_F) + \dim Y$,
where $F$ is a general fiber of $f$.
\item The Iitaka model $I_L(X)$ dominates $Y$, i.e., 
$f$ factors through $I_L$.
\end{itemize}
We say that $L$ dominate $Y$ in these cases, meaning 
that $I_L$ dominates $f$ in the sense that the base of 
$I_L$ dominates that of $f$.
\end{lem}

We remark that the lemma as stated in the reference
cited is for the algebraic case, but it is easy to see
that the lemma works for any fibration as long as $Y$
is projective (even Moishezon if we take the $H$ above
to be big rather than ample).  We refer the reader to the
above references for the proof.  The following is 
a generalization of a well-known lemma of Kodaira
to the orbifold base of a fibration.

\begin{cor} Let $f: X\ra Y$ be a fibration. Then it
is of general type if and only if there exist a
big line bundle $H$ on $Y$ such that
$h^0(L_f^mf^* H^{-1})\neq 0$ for some $m\in \NN$.
We may take $H$ to be ample if $Y$ is projective.
\end{cor}

\noindent
{\bf Proof:} This follows from the proposition above
since $L_f$ is trivial on the generic fibers of $f$
and hence of zero Kodaira dimension there. $\BOX$

\section{Special orbifolds and special fibrations}

The notion of special in our sense was originally motivated
from our search of a purely algebro-geometric characterization
of an algebraic manifold, including the quasiprojective case,
which admits a generically surjective holomorphic map from 
$\CC^n$ for some $n$, see \cite{BL1, BL2}. Then we weakened
the notion to manifolds which are holomorphically connected
(generalization of rational connectedness) and even to
manifolds whose Kobayashi metrics vanish identically.
We have succeeded soon after in finding the 
characterization of $\kappa'_+=0$ in
\cite{Lu01}, which works for all the cases where we
have a characterization of any of the above properties even
though we do not know if these properties are equivalent
in general.  We did not use the terminology of special at
the beginning when we were working with $\kappa'_+$ and its
associated Iitaka fibration as the iteration of the 
orbifold Iitaka and still conjectural Mori's fibration.
But as soon as we realized that Campana has independently 
arrived at essentially the same notion in \cite{Ca01}  
from a broader perspective but restricted to manifolds,
we have adopted the very apt terminology he introduced, 
namely  ``special,'' for this notion. Although his 
original notion in \cite{Ca01} was conjecturally 
equivalent to ours in the case of manifolds, i.e.,
without orbifold structure, he has later adopted our
notion so that there is no ambiguity in this 
terminology (except for the easy misinterpretation
of special to mean not of general type, which had
a popular informal circulation at one time). 
Neverthless, his original notion
of special should still be very useful for many
of the current applications since the above mentioned
conjectural equivalence (which is only
supported by scanty evidence given in \cite{Lu01}
outside that given by Theorem~\ref{ghs}) should be deep 
if true.

The terminology of special is also used in geometry
to mean compact manifolds with special holonomy groups
such as Calabi-Yau or holomorphic symplectic manifolds
in complex geometry. But our usage of the term is in
the context of complex birational geometry where our 
notion is more general. 

\begin{defn} Let $\yd$ be an admissible orbifold.
$\yd$ is called {\bf special} 
if it does not admit a meromorphic
fibration of general type with positive 
dimensional base. 
A fibration $f^\bd:\yd\ra Z$ 
is called special
if its orbifold general fibers are special and
if $f$ is only assumed to be meromorphic then
$\fd$ is special if
it is prebirational (birational modifications of 
$\yd$ via premorphisms and of $Z$ by morphism)
to a special orbifold fibration.
\end{defn}

The above definition of a special fibration makes
sense since the general orbifold fiber of a fibration
from an admissible orbifold remains admissible, see
Lemma~\ref{premis}. \\

Notice that in general we require $\yd$ to be an admissible 
orbifold above to make sense of a meromorphic special
fibration. It is clear that the notion of special
above would be void of meaning as a birational
invariant otherwise as easy examples would show.
We mention again the key fact that {\em smooth orbifolds
are admissible.} In general however, admissibility
can be dropped to allow for a possibly wider range of 
results even though birational invariance may
no longer be present. \\

As a direct corollary of Bogomolov's theorem,
Proposition~\ref{logCD}, we have the following
alternative characterization in category $\mathcal C$.

\begin{prop}\label{Bo} Let $\yd$ be a log-manifold
or more generally a standard orbifold 
with a $\bd$-covering. If $Y\in {\mathcal C}$, then $\yd$
is special if and only if it does
not admit any Bogomolov sheaf. \BOX\\
\end{prop}

The following theorem is simply a restatement of 
a typical theorem
in \cite{BL1,BL2} in light of our definitions here. 


\begin{theorem} Let $X$ be an algebraic surface (not
necessarily compact) not birational to a K3 surface
nor (in the case $X$ is not projective) birational to $\PP^2$,
then $X$ admits a generically surjective holomorphic map
from $\CC^2$ if and only if $X$ is special. \BOX
\end{theorem}


This is an example relating to our original
motivation for this paper for which we restrict
to surfaces for simplicity. We have since been
working with the notion of identically vanishing
Kobayashi metric in connection with special for
which Campana later conjectures to be equivalent in
the K\"ahler category.  
In this connection, we would like to 
mention that for all known examples of compact
complex manifold $X$ whose Kobayashi pseudometric
vanishes identically (and we have a good understanding
of this at least for the case of surfaces), 
$X$ is a special manifold.
But the same is not true for the infinitesimal
pseudometric (see \cite{BL1}). \\

In preparation for a discussion of examples, we state two 
results about special orbifolds and fibrations which
are fundamental in the classification theory in terms of
special objects.

\begin{prop}\label{spd}
Let $\fd: X^\bd\ra Y^\bd$ be a surjective 
premorphism between admissible orbifolds.
If $X^\bd$ is special, then $\yd$ is special.
\end{prop}

\noindent
{\bf Proof:} This follows 
directly from Lemma~\ref{gt*}. $\BOX$\\

\begin{prop}\label{spgt}
Let $\xd$ be admissible, $\fd: X^\bd \dra Y$ a meromorphic 
special fibration and $h^\bd: X^\bd \dra Z$ a 
meromorphic fibration of general type.
Then $h^\bd=g\circ\fd$ for a 
meromorphic fibration $g:Y\dra Z$. The same
statement hold if we drop ``meromorphic'' from it.
In particular, if $f$ is holomorphic
and $\yd=Y\m D(\fd)$ is special,
then $X^\bd$ is special.

\end{prop}

\noindent
{\bf Proof:} This is a direct consequence of 
corollary~\ref{tech} and the first inequality
of Proposition~\ref{MR}. \BOX\\

This elementary result is the key
to connecting the main result of this paper with 
Mori's Program and resolves in particular 
the conjectures given
in section 6.1 of \cite{Ca02} which are not already
present in Mori's Program, see our discussion right after
Conjecture~\ref{Mo}.\\


A simple example of a special manifold is an abelian
variety because any section of a positive representation
of its cotangent bundle is constant. Similarly 
$(\CC^*)^n=\PP^n\m D$ where $D$ is the sum of $n+1$
hyperplanes in general position is a special oribifold.
Since these dominates $\PP^n$, it follows that all
unirational manifolds are special.\\

It is clear that a quasiprojective curve is special
if and only if it has non-negative Euler characteristic
so that $\PP^1$ and elliptic curves are the only 
projective special curves. Also an orbifold curve
is special if and only if it has non-negative 
orbifold Euler characteristic (see \cite{BL2}), 
or equivalently,
the orbifold anti-canonical $\QQ$-bundle is
semi-ample.\\

Another class of special orbifolds consists of rationally
chain connected manifolds and smooth rationally chain connected
orbifolds, or towers of them via fibration morphisms. 
An orbifold $\xd$ is called {\bf rationally chain
connected} if any general pair of points is contained in
a chain of rational curves belonging to an irreducible family 
$\Gamma$ of {\bf orbifold rational curves}, 
i.e, nontrivial morphism from $C=\PP^1$ to 
$X$ such that the induced orbifold anti-canonical 
$\QQ$-bundle on $C$ is ample. The fact that such an orbifold
$\xd$ is special is not difficult to see
since for each fibration $\fd$ from $X$ with positive dimensional
base, the general member in this covering family $\Gamma$ of 
special orbifold rational
curves is not contained in the fibers of $f$
so that Proposition~\ref{spgt} shows that $f$ is not of 
general type.\\

A more difficult example of a special projective orbifold
is a smooth projective orbifold with trivial canonical 
$\QQ$-bundle and even one
with zero Kodaira dimension (see Theorem~\ref{Sp}). The
latter can also be deduced from the
Miyaoka's generic semipositivity theorem generalized
to the orbifold category as given in \cite{LM, Lu0}
assuming the standard conjecture
in Mori's program that an orbifold with zero Kodaira
dimension has a ``good'' model with trivial canonical 
$\QQ$-bundle. We note 
that without orbifold structures, the first instance of
this theorem concerning manifold with zero Kodaira
dimension is due to Campana in \cite{Ca01, Ca02}
with a slightly different definition of special 
initially (see also \cite{Lu02}) but neverthless
a cornerstone result in this subject. It should 
be mentioned that these facts remains true for
compact K\"ahler manifold (see section~9) but fails
otherwise since Calabi has constructed in  \cite{C}
non-K\"ahler compact complex threefolds with trivial canonical bundles 
each dominating a general type curve.\\

It is not difficult to deduce from the above results
that a smooth projective orbifold $\yd$ with semiample 
anticanonical $\QQ$-bundle $K^\vee_\yd$
is special since one can construct a
dominating smooth orbifold with trivial 
canonical $\QQ$-bundle. In fact, we show in
Theorem~\ref{-K} that if $K^\vee_\yd$ is
nef (see \cite{DPS} for example concerning 
these objects), then $\yd$ is special, which is
a conjecture in the non-orbifold category
by Campana in relation to the 
fundamental group and the Shafarevich conjecture 
concerning the universal cover, 
see also \cite{DPS0, DPS1, DPS2}.
It follows that smooth {\bf Fano orbifolds}, i.e., those 
with ample anticanonical $\QQ$-bundle, are special.
Hence, a generic hypersurface 
in $\PP^n$ is special if and only if it is of 
degree not more than $n+1$, 
since such a hypersurface has an open subset 
with trivial logarithmic canonical bundle.
These latter results can also be deduced 
directly from the fact
that Fano manifolds are rationally connected,
see \cite{KMM,Ca3}.\\

To prevent the wrong impression that Kodaira 
dimension may be strongly connected to being
special, we emphasize that among the manifolds
of dimension $n$, special ones exists for
each Kodaira dimension less than $n$ but
only the Kodaira dimension zero class 
consists entirely of special ones. We illustrate
this with some surfaces of Kodaira dimension
one.

\begin{ex} \label{K3}
Consider the Kummer K3 surface $X$
constructed from the product of two elliptic
curves as groups by the blowup of the 16 
ordinary double points that results from its
quotient by the involution of its group structure.
We obtain an elliptic fibration $f:X\ra \PP^1$
with four singular fibers of type $\text{I}_0^*$ but 
otherwise holomorphically trivial. 
Therefore, $D(f)=0.$ If we
base change the elliptic fibration by a 
finite map $\ga:C\ra \PP^1$ of positive
degree with the discriminant condition 
$\Delta(\ga)\cap\Delta(f)=\emptyset$
and denote the resulting elliptic fibration by
$f'$, then $D(f')=0$ also. Such a base change
is easily constructed for any curve $C$. Hence if 
$C=\PP^1$ or if $C$ is an elliptic curve, 
then $X'$ is special
by Proposition~\ref{spgt}. Otherwise, $f'$ is
of general type and $X'$ is not special. It can
be deduced (directly from
Kodaira's canonical bundle formula for an elliptic
fibration (\cite{BPV}) for example)
that all such $X'$ constructed has Kodaira
dimension one. In particular, if the orbifold base
of the Iitaka fibration of $X'$ has zero Kodaira
dimension, $\kappa(X')=1$ even though the general
fibers have zero Kodaira dimension. This is in
contrast to the corresponding properties of special and that
of being rationally connected as given by 
Proposition~\ref{spgt} and Theorem~\ref{ghs} 
respectively.
\end{ex}

As for special fibrations, there is first the 
rationally connected fibrations, or more generally 
the orbifold Mori's fibration defined below. Since
the general fibers of an orbifold Iitaka fibration
have zero Kodaira dimension, the Iitaka fibration 
of a smooth (or more generally an admissible)
orbifold is a special fibration. These two fibrations
combine via their iteration in the projective category 
to form the canonical special fibration of general type by
assuming the standard conjectures in Mori's program discussed
in Conjecture~\ref{Mo} and its ensuing paragraphs below. This
has been our approach to the subject in \cite{Lu01}
before \cite{Ca} when we identified this canonical special
fibration as the Iitaka fibration of the top Bogomolov sheaf. 
The latter, however, brought the general
type case of the orbifold Cnm conjecture to the fore, which
we solve completely in this paper. To
include the compact complex category in this discussion, we
will need the Algebraic reduction: Recall that if $X$ is a
compact complex manifold, then there is birational
modification $v:X'\ra X$ and a fibraion 
$Alg_X:X'\ra Y=Alg(X)$ unique up to birational equivalence
where $Y$ is a projective manifold
and the functional field of $Y$ coincides with that of 
$X$ via $Alg^*_X$ (we identify $Alg_X$ with 
the meromorpic map $Alg_X\circ v$ by abuse of notation).
If $h:X'\ra Z$ is a fibration, we also define the fibration
$$\tilde a_{X/Z}=(h, Alg_X): X'\ra Z\times Alg(X).$$

It is understood that a projective model is chosen
whenever we speak of $Alg(\bullet)$.
We will show that the Algebraic reduction is a special
fibration in the next section.\\

We now give an important class of special fibrations
that we obtain, for now, by a generalization of a well-known
weak abundance conjecture in Mori's Program which we call 
Mori's orbifold ansatz:

\begin{conj}[Mori's oribifold ansatz]\label{Mo}
Let $\xd$ be a smooth
projective effective orbifold with $\kappa(\xd)=-\infty$. 
Then $\xd$ is birational to a smooth orbifold having
an (orbifold) rationally chain connected fibration, 
called {\bf (orbifold) Mori's fibration}, in
the sense that the general fibers are (orbifold)
rationally chain connected as defined above. 
\end{conj}

Coupled with the orbifold Iitaka fibration, the 
conjectures implies that there are exactly three
possibilities for a projective smooth effective orbifold $\xd$,
i.e., the trichotomy:
$$
\begin{cases}
\text{(1) $\xd$ is of general type, or}\\
\text{(2) Up to a birational modification,
$\xd$ has a special fibration,} \\
\text{$\ \ \ \ \ $ namely the Iitaka fibration of $\xd$,}\\ 
\text{$\ \ \ \ \ $ whose general fibers are 
positive dimensional and of zero Kodaira dimension, or}\\
\text{(3) $\xd$ has an almost holomorphic special fibration,}\\
\text{$\ \ \ \ \ $ namely the Mori's fibration of $\xd$,}\\
\text{$\ \ \ \ \ $ whose general fibers are 
positive dimensional and rationally connected.}
\end{cases}
$$

Case (1) corresponds to $\kappa(\xd)=\dim X$,
(2) to $\kappa(\xd)\geq 0$ and (3) to
$\kappa(\xd)=-\infty.$ A weaker but equivalent formulation
of (3) is as follows (c.f., \cite{Db, Kol1}):\\

(3)' $\xd$ is {\bf uniruled} in the sense that
there is a countable union of subvarieties of $X$
outside of which every point is contained in an
orbifold rational curve.\\

Now if $\xd$ is not of general type, then Conjecture~\ref{Mo}
implies that it has a special fibration 
$\fd:\xd\ra Y$, which we may assume to be well adapted
by modifying $X$. If $\yd=Y\m D(\fd)$ is not of 
general type, then it would have a specical fibration 
$\gd:\yd\ra Z$ which we may again assume to be well adapted
by modifying $X$ and $Y$. Note that by Proposition~\ref{spgt},
$h^\bd=\gd\circ\fd$ is then a special fibration with base
orbifold $Z^\bd=Z\m D(h^\bd)=Z\m D(\gd),$ 
which is smooth. We can then iterate
the process with $Z^\bd$ and so forth until we reach a 
special fibration of general type from some smooth birational
model of $\xd$. This produces, along with Proposition~\ref{almo},
a canonical almost holomorphic special fibration 
of general type from $\xd$. We summarize this as follows.\\

\noindent {\bf Fact:}
{\em Let $\xd$ be a smooth projective effective orbifold. Then Mori's 
ansatz implies that $\xd$ has a canonical special
fibration of general type, which is then necessarily almost
holomorphic.}\\

We will show that this statement is true without
assuming any conjecture and thereby resolving
in particular a part of Mori's orbifold ansatz, 
see Theorem~\ref{Main}, the remark in the next
section and the same in section~11.\\

We conclude this section with the following result
promised in the last section.

\begin{prop} \label{almo}
Let $f:X\dra Y$ be a meromorphic fibration
and $X^\bd=X\m A$ a smooth orbifold such that 
$\fd$ is a fibration 
of general type. Then $f$ is almost holomorphic.
\end{prop}

\noindent
{\bf Proof:} Let $v:X_1\ra X$ be the composition of 
a sequence of blowups $v_i$ with smooth centers which
resolves the indeterminancies of $f$. Let
$f_1=f\circ v$ and impose the infinity multiplicity
on the exceptional fibers of $v$ so that $v^\bd$ 
is birational and $f_1^\bd$ is of general type. 
Suppose $f$ is not almost holomorphic.
Then the exceptional divisor $V$ of some $v_i$, say
$v'$, is mapped surjectively
onto $Y$ and such that $f_1$ restricted to 
the fibers $\PP^l=V_t$ of $v'$ is nontrivial.
We impose the smooth orbifold structure $D=v_1^*A$ on 
$V$ which we easily verify to be special at least
when restricted to the general $V_t$ since
$V\m D_{red}$ is and since $D_{red}|_{V_t}$ can 
consists of no more than $l+1$ normal
crossing hyperplanes by construction
(one only need to verify this for point
blowups of normal crossing configurations). 
It is easy to see that the 
inclusion of $V\m D$ is a morphism to
$X^\bd_1$ so that $f':=f_1|_V$ is a fibration
of general type from $V$ onto $Y$. As $v'$ is a
special fibration from $V$, the pair $(v',f')$ then
clearly contradict the first part of 
Proposition~\ref{spgt} and concludes our
proof. \BOX\\

It is clear from the proof
that the smoothness assumption on $X^\bd$
is only necessary at the points of indeterminancy
of $f$. But the result is false otherwise.

\section{The refined Kodaira dimension}

We have deduced using Mori's ansatz  that any smooth projective
effective orbifold has a canonical special fibration of general type. 
It is clear from the construction that this fibration  
can be obtained by the Iitaka fibration
of a line sheaf associated to the fibration (by Bogomolov
theorem for example). Hence there is a natural Kodaira 
dimension associated to such an orbifold given the
validity of the conjecture. We define this Kodaira dimension
for any orbifold in the following way.

\begin{defn} \label{dgt}
Let $X^\bd$ be an orbifold. If $X^\bd$ is
not special, we define
$$
\kappa'_+(X^\bd)=\kappa'(X^\bd)
=\max\{\ \kappa(L_\fd)\ |\ \fd:X^\bd\dra Y\ \
\text{is of general type}\ \},
$$
which is necessarily positive. If  $X^\bd$ is
special, we define $\kappa'_+(X^\bd)=0$ and 
$$
\kappa'(X^\bd)=
\begin{cases}
-\infty & \text{if $\kappa(L_\fd)=-\infty$ for all 
nonconstant fibrations $\fd:X^\bd\dra Y$,}\\
0 & \text{otherwise.}
\end{cases}
$$
\end{defn}

It follows from this definition that $X^\bd$ is special
if and only if $\kappa'_+(X^\bd)=0.$
There are alternative definitions which obviate the
use of cases, especially in the projective (or $\mathcal C$)
category. Of course, whether $\xd$ is special or not, we have
$$
\kappa'_+(X^\bd)=\max\{\ \kappa(L_\fd)\ |\ \fd:X^\bd\dra Y\  
\text{is of general type}\ \},
$$
and in the case $\xd$ is a log-manifold and $X\in {\mathcal C}$,
we have
$$
\kappa'_+(X^\bd)=\max\{\ \kappa(L)\ |\ L\hra\Omega_\xd^p\ 
\text{ is a line subsheaf with}\ 
\kappa(L)=p,\ p=0,\dots,\dim X\ \},
$$
by Bogomolov Theorem, Proposition~\ref{logCD}.
This was our original definition by which we worked with 
the notion of special.\\

As for $\kappa'$, we only offer the following 
alternative definitions which we deduce for now
from Mori's orbifold ansatz. It is actually true
without any assumption by the subadditivity of
Kodaira dimension theorem that we obtain in the next section.
It brings out (and reduces
Mori's orbifold ansatz to) an important 
conjecture concerning rationally connected orbifold.

\begin{lem}\label{defk}
Let $\xd$ be a smooth projective effective orbifold. Assume Mori's
ansatz. Then 
\begin{align*}
\kappa'(X^\bd) &=\max_{m>0}\ \min \{\ \kappa(L_\fd)
\ |\ \fd:X^\bd\dra Y\  \text{a special fibration 
with}\ 0\leq \dim Y\leq m\ \} \\
&=
\begin{cases}
-\infty\ \ \ \ \ \ \ \ \ \  \ \ \ \ 
\text{if  $\kappa(L_\fd)=-\infty$ for all 
nonconstant fibrations $\fd:X^\bd\dra Y$.}\\
\min \{\ \kappa(L_\fd)\ |\ 
\fd:X^\bd\dra Y\  \text{a special fibration with}\ \dim Y>0\ \}
\ \ \text{otherwise.}
\end{cases}
\end{align*}
\end{lem}

\noindent
{\bf Proof:} The two expressions on the right hand side are
easily seen to be the same, which we denote by $\bar\kappa$.
We have seen that Mori's orbifold ansatz implies
the existence of a canonical special fibration of general type
$\fd:\xd\dra Y$. If $m=\dim Y$ is positive, then it is
clear by construction that $\bar\kappa=m=\kappa'(X^\bd).$ Hence
we assume that $\dim Y=0$, i.e. that $\xd$ is
special. This means we can take $\xd$ to be
a tower of (orbifold) Iitaka and Mori's fibrations over a point.
If there is an (orbifold) Iitaka fibration in this tower, then 
we have again by construction that
$\bar\kappa=0=\kappa'(X^\bd)$. Otherwise, this tower
consists entirely of orbifold rationally connected fibrations.
The result now follows from the following Lemma. \BOX\\

\begin{lem}\label{tower}
Suppose $\xd$ is a smooth orbifold having a tower of 
Mori's fibration over a point. Then 
$\kappa'(\xd)=-\infty$.  If further $\xd$ is a log manifold,
then $\kappa(L)=-\infty$ for all subsheaves $L$ of $\Omega^p_\xd$
with $p>0$.
\end{lem}

\noindent
{\bf Proof:} Since the Kodaira dimension of a fibration
whose orbifold base is an orbifold rational curve is
by definition $-\infty$, it follows by Proposition~\ref{tech0}
that any fibration from $\xd$ of nonnegative Kodaira dimension
must factor through each stage of the tower and hence must
have zero dimensional base. The last statement of the lemma
without the log structure on $X$ is a theorem of Campana in 
\cite{Ca} whose proof is later completed by Theorem~\ref{ghs} for
which we refer to \cite{Ca} for the very brief details. 
The log case is an
easy generalization for which we leave to the reader.
\BOX\\

It follows from the above proofs that Mori's
orbifold ansatz implies the following statement, which
we state as a conjecture (\cite{Ca} has
a similar one in the non-orbifold case).

\begin{conj} \label{rc} A smooth projective orbifold 
$\xd$ is birational to one with a tower of Mori's
fibration if and only if $$\kappa'(\xd)=-\infty.$$ 
\end{conj}

Conversely, this conjecture implies Mori's orbifold 
ansatz by 
assuming the subadditivity of Kodaira dimensions in the 
orbifold case for which we refer to the next section.

\begin{rem}\label{rm}
Conjecture~\ref{rc} implies Mori's orbifold ansatz,
i.e., the orbifold weak abundance conjecture as
stated in Conjecture~\ref{Mo}.
\end{rem}

\noindent
{\bf A proof of Remark~\ref{rm}:} Let $\xd$ be a smooth
projective orbifold with $\kappa(\xd)=-\infty$ and $\dim X=n$. 
We need to prove that $\xd$, up to a birational modification,
has a Mori's fibration with positive dimensional fibers.
We proceed by induction on $n$ so assume that the
result is true for such orbifolds of smaller dimension
(the result being trivially true in dimension one).
We may assume that there is a fibration $\fd:\xd\dra Y$
with $\dim Y>0$ and $\kappa(Y,\fd)\geq 0$.
For otherwise $\kappa'(\xd)=-\infty$ so that
$\xd$ is birational to a smooth projective
orbifold with a rationally connected fibration
by Conjecture~\ref{rc}.  We may assume that $f$ is
holomorphic and $\fd$ is well adapted by modifying 
$\xd$ and $Y$.
If the general fibers of $\fd$ has nonnegative
orbifold Kodaira dimension, then the subadditivity
assumption would give $\kappa(\xd)\geq 0$ which contradicts
our assumption. Hence, there is an infinite set of
fibers $F^\bd$ of $\fd$ above the complement of 
$\Delta(\fd)\cap \Delta(f)$ with $\kappa(F^\bd)=-\infty$. 
But this implies that these $F^\bd$'s admits Mori's 
fibrations by our induction hypothesis. Since these
fibrations are almost holomorphic, Proposition~\ref{AH} 
and the standard argument of Fujiki in \cite{Fujk}
shows that the general fibers of $f$ have Mori's fibrations 
and that these fibrations
patch up to give a Mori's fibration for $\xd$.
\BOX\\

Conjecture~\ref{rc} is known for projective threefold
because Mori's ansatz is solved by Miyaoka in \cite{Miya1}
in this case (see also \cite{MP}). It is also known 
for some (but not yet all) cases of orbifold surfaces.
We will address this conjecture for orbifold surfaces
with some remarks for higher dimensions in a 
forthcoming paper. We remark that this conjecture
implies a conjecture of Koll\'ar in \cite{Kol1} that
a manifold $X$ is  rationally connected if and only if
for all $l\in \NN$,
$$
H^0(\Omega_X^{\otimes l})=0.
$$

\begin{rem}\label{rm2}
Lemma~\ref{defk} and Proposition~\ref{spgt} shows that
$\kappa'_+(\xd)$ is precisely the dimension of the base
of the canonical special fibration of general type
constructed in the last section, assuming Mori's
orbifold ansatz. In fact, we show in Section 10 
(Theorem~\ref{M}) that if $\xd$ is not special,
then there is only one Bogomolov sheaf $L_\fd$ corresponding
to $\kappa'_+(\xd)$ in Definition~\ref{dgt}
and therefore only one general type 
fibration $\fd:\xd\dra Y$ up to birational
equivalence of $Y$ where
$\kappa'_+(\xd)=\dim Y$. We call this 
birational equivalence class of fibrations
the {\bf basic fibration (or reduction) of} $\xd$ 
and this line sheaf the {\bf top Bogomolov sheaf} 
of $\xd$.
\end{rem}

\section{Proper behavior in an analytic family}

We are interested here in the behavior of $\kappa'_+$, 
Bogomolov sheaves and general type fibrations in a
family. We approach this here following partly the
ideas of Campana in the construction of his relative 
general type reduction in \cite{Ca02} but generalizing 
them to the compact complex full orbifold category.
We actually only need to work out the required behavior in
a family in the projective category rather than 
in the category $\mathcal C$ Campana
chose to work out his relative general type reduction.
But for the sake of the subject
matter and for proper attribution, as he was the 
first to come up with the idea of the 
relative general type
reduction, we will work out his ideas in category 
$\mathcal C$ clarifying some necessary points in our 
more general full orbifold context. However, we
would like to mention that these ideas can be 
circumvented by more direct (and perhaps more
elementary) methods (using for example Fujiki's
theorem on the projectivity of the space of
divisors in a family, especially
\cite[Theorem 7]{Fujk2} and \cite{Fujk3}).\\

We start with the following proposition 
which follows from 
the upper semicontinuity 
of the dimension in an analytic family or simply by 
quoting \cite{LS} for example.

\begin{prop}\label{us}
Let $f:X\dra Y$ be a meromorphic fibration 
and $h: X\ra Z$ a fibration. Suppose $\xd$ is
a smooth or an admissible orbifold.
Then $f_z^\bd:X_z^\bd\dra Y_z=f(X_z)$
is a well defined meromorphic fibration for $z$
in a Zariski open subset $U$ of $Z$ and $\kappa(Y_z, f_z^\bd)$ 
is constant for $z$ in the complement $V$ of a countable number of
subvarieties in $X$ with $V\subset U$, i.e., for general $z\in Z$.
\end{prop}

\noindent
{\bf Proof:} Consider the meromorphic fibration
$\bar f=(h, f): X\dra \bar Y$ where 
$\bar Y=\bar f(X)\hra Z\times Y$. It is clear
that $\bar f_z=f_z$ is well defined as a meromorphic
map for $z$ in a Zariski open subset $U$ of $Z$ 
and hence so are $f(X_z)$ and $f_z^\bd$. By
blowing up $X$ and $\bar Y$, we may assume that 
$\bar f$ is an admissible morphism and $\bar Y$ is smooth
since this procedure does not change 
$\kappa(\bar Y_z, \bar f_z^\bd)$ 
for the general $z$ (by choosing the $\infty$ multiplicity
for the exceptional divisors of the blowups for example).
Let $k=\dim \bar Y -\dim Z$.
{}From the morphism (actually an inclusion modulo torsion) 
$\bar f^*\Omega_{\bar Y/Z}\ra \Omega_{X/Z}$ 
of the sheaf of relative differentials, we take double
dual of $\wedge^k$ of the saturation of the image of
the above morphism of sheaves to obtain a line subsheaf
$L'$ of $\wedge^k\Omega_{X/Z}^{\vee\vee}$ over $X$.
Let $i:X_z\hra X$ be the inclusion of 
the generic fiber of $h$ and let $()_z$
denote $i^*$ or the restriction to $X_z$
by abuse. By construction, 
$i^*\Omega_{X/Z}^{\vee\vee}=\Omega_{X_z}$ 
and $K_{Y_z}$ are locally free and
$L'_z$ is the saturation of $\bar f_z^*K_{Y_z}$
in $\Omega_{X_z}$, i.e., $L'_z=L_{\bar f_z}$. 
Let $L=L'(A\cap \bar f)$. Then 
$L_z=L_{\bar f_z^\bd}=L_{f_z^\bd}$
and so the result follows by looking at the 
sequence of images of natural maps 
$X_{U_m}=h^{-1}(U_m)\dra {\rm \bf Proj}(h_*(L^m))$
over $U_m$ where $U_m$ is a Zariski open subset of $Z$ on which
$h$ is smooth and $h_*(L^m)$ has constant rank.
\BOX\\


Now the basic idea of Campana in \cite{Ca02} is to look
at the family of cycles given by the fibers of fibrations.
For this we recall the theorem of Barlet (\cite{Ba}) that for a
compact complex space $X$, there exists a reduced 
complex space $B_q(X)$ and an analytic family 
$\{A_b\ |\ b\in B_q(X)\}$ of $q$-cycles on $X$ such 
that every analytic family $\{A_s\ |\ s\in S\}$ of
$q$-cycles on $X$ is induced by a unique morphism
$h:S\ra B_q(X)$, i.e., $A_s=A_{h(s)}$ for all $s\in S$.
Every analytic family $S$ of cycles of $X$ corresponds
to an analytic subset $i: A_S\hra S\times X$ defined by 
$A_S=\{(s,x)\ |\ x\in {\rm supp}(A_S)\}$, see 
\cite[Th\'eor\`eme 1]{Ba}, called the graph of $S$.
It is known (\cite{Fujk1}, for example) that
$B_q(X)$ has a countable base for its topology if $X$ is so.
A terse acount of cycle space theory can be found in \cite{CP}.\\

For the next lemma, 
suppose $f:X\dra Y$ is an almost holomorphic fibration
where $X$ is a compact complex manifold as usual.
Let $U$ be a Zariski open subset of $Y$ above which $f$ 
is defined and smooth. Then we have a map $k:U\ra B_q(X)$
where $q=\dim X-\dim Y$.

\begin{lem}\label{k_f}
$k$ extends to a birational map $k_f:Y\dra B_q(X)$
onto a (necessarily compact) component $S_f$ of $B_q(X)$.
Conversely, every compact component $S$ of $B_q(X)$ satisfying
the following condition {\rm (*)} gives rise to an almost
holomorphic (rational) map $f_S:X\dra S$ with $q=\dim X-\dim Y$.
\begin{itemize}
\item[\rm (*)] The general members of $S$
as $q$-cycles of $X$ are 
irreducible with multiplicity one and whose supports are disjoint
and fills up an open subset of $X$.
\end{itemize}
\end{lem}

\noindent
{\bf Proof:} 
For the second statement, consider the graph 
$i: A_S\hra S\times X$ of $S$. Let $p_1$ and $p_2$ be 
the projections of $S\times X$ onto the 
respective factors. Then $p_2\circ i$ is a map
onto $X$ defined and one-to-one above an open
subset of $U$ of $S$ by condition (*). Hence its inverse 
composed with $p_1$ is an almost holomorphic rational map
by the compactness of $S$ and we can identify this
map as $f_S$.

For the first statement, we can flatten $f$ via 
birational modifications $v:X'\ra X$ and $u:Y'\ra Y$
where by construction all fibers of $f': X'\ra Y'$ map
to $q$ cycles in $X$ via $v$. Hence we obtain a 
morphism $r: Y'\ra B_q(X)$. Consider now the 
component $S$ of $B_q(X)$ containing $r(Y')$ and
the map $h=(id,f)\circ i:A_S \dra S\times Y$ 
where $i: A_S\hra S\times X$ is the graph of $S$. 
By hypothesis, there exists
$s\in S$ such that $f$ (and hence $h$) is holomorphic
in a neighborhood of $A_S$ in $X$ (resp. in $A_S$)
and $h(A_s)$ is a point. So by upper semicontinuity
of the dimension of the fibers of the image of $h$
over $S$, we have that $h(A_s)$ is a point for $s$
in an open subset $U$ of $S$ (on which $h(A_s)$ is
defined). This forces $r(Y')$ to contain $U$.  \BOX\\

\begin{cor}
The set of almost holomorphic (and hence of general type)
fibrations from an orbifold $\xd$ is at most countable 
and is discrete in the compact open topology. \BOX
\end{cor}

We now generalize Lemma~\ref{k_f} to the relative situation.
Let $h: X\ra Z$ be a fibration.
Recall that by \cite{Ba, Fujk0}, there exists a 
relative Barlet space 
$B_q(X/Z)=\cup_{z\in Z}B_q(X_z)$ 
which is a reduced complex subspace of $B_q(X)$.
There is a natural map $h_*:B_q(X/Z)\ra Z$ which 
is holomorphic up to a normalization.
Furthermore, Fujiki in \cite{Fujk0,Fujk4} 
has shown using the result of
of Harvey-Shiffman (\cite{HS}) (itself a variation of a well
known theorem of Bishop used in \cite{Lieb} for the same 
purpose) that all connected (hence also irreducible) components 
of $B_q(X/Z)$ are compact if $X\in \mathcal C$ and more 
generally if $h\in {\mathcal C}/Z$. For our purpose, it is not
necessary for us to know the precise meaning of 
$h\in {\mathcal C}/Z$ which roughly means that
each fiber is in $\mathcal C$ where $X\in \mathcal C$
means that $X$ is holomorphically 
dominated by a K\"ahler manifold.
It is sufficient for us to know that $h:X\ra Z$ is in 
$\mathcal C/Z$ (in fact, K\"ahler) if it is a projective
morphism or if $X$ is K\"ahler. Let 
$Z_0=Z\smallsetminus \Delta(h)$ and $X_0=h^{-1}(Z_0)$.
We will need to consider
the subset $B_q^F(X_0/Z_0)$ of $B_q(X_0/Z_0)=h_*^{-1}(Z_0)$
defined to be the largest subset of cycles in $B_q(X_0/Z_0)$
whose intersection with $B_q(X_z)$ is either empty
or fills up $X_z$ for each $z\in Z_0$.

\begin{prop}\label{AH}
Let $AH_q(X_0/Z_0)$ be the subset of $B_q^F(X_0/Z_0)$ defined
by the union of $S_f$ (defined in Lemma~\ref{k_f}) 
over all almost holomorphic rational maps $f$ from $X_z$ 
as $z$ vary over $Z_0$. Assume that all 
components of $B_q(X/Z)$ are compact, which is
the case, in particular, if $X\in \mathcal C$
or more generally if $h\in {\mathcal C}/Z$.
Then there exist a Zariski closed subset $B_q^F(X/Z)$
of  $B_q(X/Z)$ containing $B_q^F(X_0/Z_0)$ as a
Zariski open subset and  
$B_q^F(X_0/Z_0)=B_q^F(X/Z) \cap B_q(X_0/Z_0)$.
Moreover, 
$AH_q(X_0/Z_0)$ consists of Zariski open subsets
of irreducible components of $B_q^F(X/Z)$.
\end{prop}

\noindent
{\bf Proof:} The first claim is clear by upper
semicontinuity of dimension in a family.
We may assume that the $f$ considered
above are all almost holomorphic fibrations by Stein
factorizing. Fix $z$ and an almost holomorphic 
fibration $f:X_z\dra Y_z$. Let $S$ be a component
of $B_q^F(X/Z)$ containing the general fibers of $f$.
By Lemma~\ref{k_f}, a component of $S\cap B_q(X_z)$
is $S_f$ and so satisfies (*). As $S$ is compact
by our assumption on $X$, the result follows since the
condition for a component of $S\cap B_q(X_{z'})$
to satisfy (*) is a Zariski open condition on the 
ramified cover of $Z$ defined by $S$ (again
by upper semicontinuity of dimension in a family).
\BOX\\

\begin{prop}\label{GT_q}
Let $\xd$ be a smooth or admissible orbifold 
and $h: X\ra Z$ a fibration with $Z$ projective or Moishezon.
Let $GT_q(\xd/Z)$ be the subset of $B_q(X/Z)$ defined
by the union of $S_f$ over all general type fibrations
$\fd$ from $X^\bd_z$ as $z$ vary over $Z$.
Then, above a countable intersection of Zariski 
open subset of $Z$, 
$GT_q(\xd/Z)$ coincides with 
a countable union of (Zariski opens of) 
compact complex analytic
subvarieties of $B_q(X/Z)$.
\end{prop}

\noindent
{\bf Proof:} We only work out the case $\xd$ is smooth
as the admissible case would follow by the basic lemmas
already established on premorphisms. We assume for
simplicity that $Z$ is projective.
Let $\Omega=\Omega^s_{X/Z}$ with $s+q=\dim X/Z$ and
$r: P\ra {\rm \bf Proj}\ \Omega$ 
be a resolution of singularity of ${\rm \bf Proj}\ \Omega$
with $p: P\ra X$ the projection composed with
$r$. Let ${\mathbb I}=r^*{\mathcal O}_\Omega(1)$ 
and ${\bar h}=h\circ p$. 
We may choose a smooth projective $Alg(P)$
and $r$ so that $p^{-1}(D_{red})$
is normal crossing and so that 
the algebraic reduction $Alg_P:P\ra Alg(P)$
of $P$ is a morphism by Hironaka's theorem.
We focus our attention on the open set 
$U_0:=Z\smallsetminus (\Delta(\hd)\cup \Delta(\bar\hd)),$ 
where $p$ has multiplicity one.
Note that $P_z^\bd$ and $X^\bd_z$ are smooth for $z\in U_0$.
So Lemma~\ref{ford}, or its proof, shows that we can
find an orbifold structure on $P$ (by choosing the $\infty$ 
multiplicity on the exceptional fibers of $r$ and
choosing $\bd P^\bd$ to be $p^*D$ over $U_0$ otherwise for 
example) so that $p^\bd:P^\bd\ra \xd$ is a well adapted 
special fibration over $U_0$. This means in particular
that there is a Zariski open subset $U$ of $Z$ contained 
in $U_0$ where $p^\bd_z: P_z^\bd\ra X^\bd_z$
is a well adapted special fibration for all $z\in U$.
By blowing up further (again choosing the
$\infty$ multiplicity on the exceptional fibers if one
prefers), we may assume that  $Alg_P^\bd:P^\bd\ra Alg(P)$
and $\bar a^\bd:P^\bd\ra \bar Y$ are well adapted 
fibrations (over $U$) where 
$r':\bar Y\ra \tilde a_{P/Z}(P)$ 
is a resolution of singularity of 
$\tilde a_{P/Z}(P)\subset Z\times Alg(P)$ 
and $\bar a$ is defined by 
$\te a_{P/Z}= r'\circ \bar a$ (see the 
definition just before Conjecture~\ref{Mo}). 
Hence by shrinking $U$ if necessary, we may
assume that $p^\bd_z: P_z^\bd\ra X^\bd_z$
and $\bar a_z^\bd:P^\bd_z\ra \bar Y_z=(r')^{-1}(p_1^{-1}(z))$ 
are well adapted fibrations for all $z\in U$ and that
$(p_1\circ r')^\bd$ is smooth over $U$.

Now fix a point $z\in U$ and let $i:P_z\hra P$ be
the inclusion of the fiber over $z$.
As $p^\bd_z: P_z^\bd\ra X^\bd_z$
is a well adapted special fibration, it follows from 
Proposition~\ref{spgt} that every general type
fibration ${\bar f}^\bd: P_z^\bd\dra Y_z$
gives rise to a general type fibration 
$f^\bd: X_z^\bd\dra Y_z$ such that 
${\bar f}=f\circ p$ and conversely. Now fix
such a $\bar f^\bd$ and its associated $f^\bd$.
Let $L=L_\fd$. Then we have an inclusion of line sheaves
$L\hra L_f(D_z)$. Since $p_z^* L_f$ is a subsheaf of 
$\mathbb I_z$, we see that $p_z^*L$ is a subsheaf of 
the line sheaf $J_z$ where $J={\mathbb I}(p^*D)$.
Also there exist $l\in \NN$ such that $\fd$ 
(resp. $\bar\fd$) is defined by the linear system 
of $L^l$ (resp. of $p_z^*L^l$) by definition.
Let $V_l={\bar h}_*(J^l)$ and $H_l$ be an ample divisor
on $Z$ such that $V(l):=V_l(H)=\bar h_*(J^l(\bar h^*H))$
is spanned. Let $U_l$ be a Zariski open subset of $X$
contained in $U_0\smallsetminus (H_l)_{red}$ over which
$V(l)$ has constant rank and where the base change
map $V(l)/{\mathcal I}_z V(l)\ra H^0(J_z^l(\bar h^*H))$
is bijective (see \cite[Theorem I.85]{BPV}). 

Suppose now and furthermore that $z\in U_l$.
Then the natural restriction map 
$$\eps:  H^0(J^l(\bar h^*H))=H^0(V(l))\ra  
V(l)/{\mathcal I}_z V(l)=H^0(J_z^l(\bar h^*H))$$
is surjective.  Consider the subspace $V_z=H^0(p_z^*L^l)$
of $H^0(J_z^l)=H^0(J_z^l(\bar h^*H))$ and let 
$V=\eps^{-1}(V_z).$ Then the image of the natural
map $f'_z: P_z\dra \PP(V_z^\vee)\subset \PP(V^\vee)$
is by definition birational to $Y_z$ so that this
map to its image is birational to the almost holomorphic 
fibration $\bar f$ given by $\bar \fd$. 
On the other hand, the natural map
$f':P\dra \PP(V^\vee)$ factors through $Alg_P$ and
hence through the fibration given by $\bar a:P\ra \bar Y$,
i.e., there is a map $f'': \bar Y\dra \PP(V^\vee)$ with
$f'=f''\circ \bar a$.
Since $f''\circ \bar a\circ i=f'_z$, there is a map 
$e: \bar Y_z=\bar a(i(P_z))\dra Y_z$ 
birationally equivalent $f''|_{\bar Y_z}$ with
$e\circ \bar a\circ i=\bar f$. Note that $e$ is an almost
holomorphic fibration onto $Y_z$ since $\bar f$ is. Hence
the associated orbifold fibration 
$e^\bd: \bar Y_z^\bd=\bar Y_z\m D(\bar a_z^\bd)\dra Y_z$
is of general type since $\bar \fd$ is and since
$\bar a_z^\bd=\bar a\circ i^\bd$ is admissible.
Note that the key to the argument of this paragraph is 
the inclusion 
$\PP(V_z^\vee)\subset \PP(V^\vee)$ induced by $\eps$.

Finally let $z\in W=\cap_l U_l$. We have just shown
that a general type fibration $\fd$ from $X_z^\bd$
gives rise to a general type fibration $e^\bd$
from $\bar Y_z^\bd$ and conversely and that the
correspondence is obtained on the level of cycles
in the form $(p^*)^{-1}\bar a^*S_e=S_f$. Here the
map $\bar a^*$ can be defined ad-hoc on a component
$S'$ of $B_q^F(\bar Y/Z)$ intersecting 
$(p_1\circ r')_*^{-1}(W)$ by associating the
cycle $k[\bar a^{-1}({\rm supp} A_s)]$ to each $A_s$ 
for the general $s$ in $S'$
where $k$ is the multiplicity of $A_s$. We note that
$B_q^F(\bar Y/Z)$ is well defined with compact
components by
Proposition~\ref{AH} since $\bar Y$ is projective.
As before, $\bar a^*$ extends to a
meromorphic map from $S'$ birationally into $B_q(P/Z)$ 
(via the usual flattening argument).
Using $e\circ \bar a\circ i=f\circ p_z$, we see easily
that $\bar a^*(S')=p^*(S)$ for some compact subvariety
$S$ in $B_q(X/Z)$ bimeromorphically equivalent to 
$S'$ and on which $p^*$ is bimeromorphic. Hence, 
the map $(p^*)^{-1}\bar a^*$ 
is bimeromorphic on 
the relevant components of $B_q(\bar Y/Z)$
and so preserves compactness and the result follows. 
\BOX\\

The question naturally arises as to whether the above 
proposition is true for $Z\in {\mathcal C}$ or more 
generally for $Z$ arbitrary. The latter
seems farfetched while the former should be
possible as it is at least true for $Z$ Moishezon
and also true for $X\in \mathcal C$ for which we
know by the definition of $\mathcal C$ 
that $Z\in \mathcal C$ in this
case. By combining Propositions \ref{us}, \ref{AH} 
and \ref{GT_q}, we obtain the following corollary:

\begin{cor}
Let $h: X\ra Z$ be a fibration with 
$h\in \mathcal C/Z$ or with $Z$ Moishezon. 
Then $\kappa'_+(X_z)$
is constant for general $z\in Z$. \BOX
\end{cor}


Remark~\ref{rm2} now allows us to deduce the following
corollary, which gives a canonical almost holomorphic fibration 
$$ bc_\hd: \xd\dra bc(\hd)$$
associated to a fibration $\hd:\xd\ra Z$
known as the {\bf relative  basic reduction} of $\hd$,
from Theorem~\ref{M} and the above proposition and its proof.

\begin{prop}\label{GTR}
Let $h:X\ra Z$ be a fibration such that either $h\in \mathcal C/Z$
or $Z$ is Moishezon. Let $\xd$ be a smooth or admissible 
effective orbifold. Then there exist a fibration 
$Y\ra Z$ and an almost holomorphic special fibration $\fd:\xd\dra Y$
such that $f_z^\bd: X_z^\bd\dra Y_z$ is the basic
fibration of $X_z^\bd$ for the general $z$ in $Z$. 
Furthermore, $Y$ is dominated by a projective fibration
over $Z$.
\end{prop}

\noindent
{\bf Proof:} By Proposition~\ref{GT_q}, there exists a 
component $Y'$ of $B_q(X/Z)$ whose fiber $Y'_z$
at a general point $z\in Z$ consists of component of 
$GT_q(\xd/Z)_z$, each corresponding to a general type
fibration from $X_z^\bd$ up to birational modification
of the base. Let $q=\dim X_z-\kappa'_+(X_z)$
which is well defined by the above corollary. Then
Theorem~\ref{M} says that there is a unique general type
fibration from $X_z^\bd$ up to birational modification
of the base, which is called the basic fibration of
$X_z^\bd$. This means that the natural projection 
from $Y'$ to $Z$ is connected and its general fiber
$Y'_z\subset B_q(X_z)$ gives, via its graph, 
the basic fibration of $X_z^\bd$. Hence by considering
the graph of $Y'$ and letting $Y$ be a desingularization
of $Y'$ gives us the result.  \BOX\\

Hence, the basic reduction of $\hd$ restricted to
the general fibers $F$ of $h$ are the basic reductions
of $F^\bd$. Coupled with Theorem~\ref{M}, also proved 
in \cite{Ca02} for less general settings, 
this is also a result in \cite{Ca02} for
the case $X\in {\mathcal C}$ without orbifold structures.\\

Another important corollary, also with the help of 
Theorem~\ref{M}, is the specialness of the algebraic
reduction already mentioned and promised:

\begin{theorem}\label{Alg} 
Let $\xd$ be a smooth or admissible effective orbifold and
$$h=Alg_X: X\dra Z=Alg(X),$$  
the algebraic reduction of $X$.
Then $\hd:\xd\dra Z$
is a special meromorphic orbifold fibration. 
\end{theorem}

\noindent
{\bf Proof:} By blowing up, we may assume that $h$ is
a morphism. The fact that $h$ is a fibration follows
from the definition and Stein factorization 
(\cite[Proposition~3.4]{Ue}). If $\hd$ is not 
special, then $Y:=bc(\hd)$ would have dimension 
greater than that of $Z=Alg(X)$ by the last proposition.
As $Y$ is dominated by a projective fibration
over $Z$ and $Z$ is projective, $Y$ is Moishezon by
\cite[Corollary 3.9]{Ue} and so its algebraic
dimension ${\rm a}(Y)$ is strictly greater than
that of $Z$. But $Y$ is dominated by $X$ (via $bc_\hd$)
so that $\dim Z={\rm  a}(X)\ge {\rm a}(Y)> \dim Z$. This is 
a contradiction and hence $\hd$ is special. \BOX

\section{The orbifold Cnm conjecture}

We derive in this section the most important ingredient for our
main theorems. This is obtained from the ``positivity'' 
theorems of Kawamata and Viehweg, which were extensively
developed by them with an important
contribution by Koll\'ar and some initial
contributions by Fujita and Ueno among others (see
the article of Mori \cite{Mori} for an excellent thorough
account of the theory and its role in classification).
The key to our original
approach, beside that given in \cite{Lu02} (though already
known to the author at the time), is a coupling of
(1) Viehweg's famous (very deep but ingenious)
proof that weak positivity (see below for
the definition) of the direct images of dualizing sheaves
implies the same for tensor powers of dualizing sheaves
with (2) the base change property for dualizing sheaves
(both of which are essential ingredients in all the proofs of 
Kawamata-Viehweg type theorems on additivity of Kodaira 
dimension). The author is extremely indepted to Eckart
Viehweg for a chance to discuss the proof 
and for providing an independent proof
by him. He would also like to thank G. Dethloff
for first inviting him to Brest where these
ingredients were dicussed and given in an algebraic geometry
seminar talk as well as the organizers of Oberwolfach 
Seminar on Fundamental groups in Geometry (September 2002)
where he gave a seminar talk on this.\\

We restrict to the K\"ahler category in this section
unless specified otherwise. The fact that the theorems
used in this section
treated in the literature in the projective category
are actually valid in the K\"ahler case 
(even in Fujiki class $\mathcal C$) 
is believed to be well known among the experts but we
first learned this from Viehweg, see 
related comments and discussions in \cite{Vi2, Kol, Na0}.
For our purpose however, there is no loss of generality
by restricting to the projective category.
Recall that an orbifold is called {\bf effective} 
if its boundary divisor is.
We first state an orbifold generalization of a well known
conjecture of Iitaka concerning the Kodaira dimensions
of an algebraic fiber space.

\begin{conj}[Orbifold Iitaka's Cnm Conjecture]\label{Iit}
Let $\yd$ be a smooth or an e-admissible effective 
orbifold. Let $\gd:\yd\ra Z$ be a fibration and $Y_z^\bd$ its
orbifold general fiber.  Then
$$\kappa(\yd)\ge \kappa(Y_z^\bd)+\kappa(Z, \gd).$$
More generally, if $\xd$ is an e-admissible orbifold,
$f:\xd\ra Y$ a morphism with effective orbifold base $\yd$,
$g:Y\dra Z$ meromorphic and $h=g\circ f$ almost holomorphic.
Then we have with $f_z=f|_{Y_z}$ that
$$\kappa(Y, \fd)\ge \kappa(Y_z, f_z^\bd) + \kappa(Z, \hd).$$
\end{conj}

The author thanks F. Bogomolov for calling this conjecture to
the author's attention in connection to the canonical
fibration. Various forms of this type of orbifold Iitaka
conjecture can be found in \cite{Ca01, Ca02, Lu02}. However,
the first theorem of this type, with orbifold structure both
on the total space and on the base, precedes these by two
decades and is due to Kawamata in \cite{Ka}.\\

Without the orbifold $\del$, 
the Cnm conjecture of Iitaka is known in the 
case of general type $Z$ by theorems of Viehweg and 
Kawamata (see \cite{En,Ka, Vi,Vie, Mori})
and in the case of general type $Y_z$
by a theorem of Koll\'ar (\cite{Kol}).
We prove here the orbifold analog of their theorems.
We begin with a special case which is also given
in \cite{Lu02} and in \cite{Ca02}.

\begin{theorem}\label{Cnm} Let $f:X\ra Y$ and $g:Y\ra Z$ 
be fibrations with $g$ of general type.
Then $\kappa(Y^\del)=\kappa(Y_z, f_z)+ dim(Z)$ where
$Y_z$ is the general fiber and $f_z=f|_{Y_z}$.
\end{theorem}

This is  a direct corollary via Lemma~\ref{fta}
of the following theorem, which generalizes slightly 
the results of the people
mentioned above.

\begin{theorem}\label{wp} 
Assume that $f: X\ra Y$ is a fibration. Set 
$K_{X/Y^\del}=K_X\otimes f^*K_{Y^\del}^{-1}$ as $\QQ$-line 
bundles. 
Let $m$ be a positive integer such that $mD(f)$ is a Cartier divisor
on $Y$. Then there is an effective divisor $E$ 
with $E_{red}\subset E(f)$ such that 
$f_{*}(K_{X/Y^\del}^{\otimes m}(E))$ is weakly 
positive in the sense of \cite{Vi}.
\end{theorem}

Without the orbifold $\del$, this is the Kawamata-Viehweg's theorem.
We will not define weak positivity except to note the following
two properties (see \cite{En}).

\begin{itemize}
\item[(P1)] A torsion free coherent sheaf is weakly positive if it 
contains a weakly positive subsheaf of the same rank.
\item[(P2)] Let $\Upsilon$ be a torsion free coherent sheaf on $Y$.
Let $v: Y'\ra Y$ be a flat base change. If $v^*\Upsilon$
is weakly positive, then so is $\Upsilon$.
\end{itemize}

\noindent
{\bf Proof of Theorem~\ref{wp}}:
We denote $D(f)$ by $D^f$ for simplicity.
Recall that $D=D^f_{red}=\sum_i D_i$ is a 
simple normal crossing divisor on $Y$. 
Write $f^*D_i=\sum_j m_{ij}D_{ij}$ and let 
$M_i$ be the least common multiple of 
$\{\ m_{ij}\ |\ f(D_{ij})=D_i\ \}$. 
The Kawamata branched covering trick (\cite{Ka}) 
guarantees the existence of a finite Galois cover 
$v: Y'\ra Y$ with $Y'$ smooth such that $v^*D_i=M_iD'_i$ 
for some effective divisor. We let 
$u:X'\ra X$ be the composition $u=u_1\circ d$, where
   $u_1:X_1\ra X$ is the base change of $X$ by $v$
composed with the normalization map, and
$d:X'\ra X_1$ is a desingularization which 
is an isomorphism above the smooth locus of $X_1$.
We also denote by $f_1:X_1\ra Y'$ and $f'=f_1\circ d:X'\ra Y'$ the 
resulting pullbacks of $f$ by $v$. Note that if $f(D_{ij})=D_i$,
then $u^*(m_{ij}D_{ij})=M_iD'_{ij}$ 
for some effective divisor $D'_{ij}$ by construction.
Hence, applying $u^*$ to the $\QQ$-divisor $f^*D^f$ renders
it into an effective (integral) divisor on $X'$.

Now, a local computation shows the
following natural 
inclusion of torsion free coherent sheaves of
the same rank outside a codimension two subset of $Y'$:  
$$
(f')_{*}\big(K_{X'/Y'}^{\otimes m}\big)
\hra v^{*}f_{*}\big(K_{X/Y}(-(df))^{\otimes m}\big)
\hra v^{*}f_{*}\big(K_{X/Y}(-f^*D^f)^{\otimes m}\big)
=v^{*}f_{*}\big(K_{X/Y^\del}^{\otimes m}\big).
\ \ \ \ \ (*)
$$
A round about way to see this is as follows:
This fact without the $(-f^*D^f)$ is well-known
(e.g., \cite{Ka}, \cite[Lemma 13]{En}). Hence, we
only need to check $(*)$ above an analytic neighborhood
$V$ in $Y$ of a generic point $p\in D'_i=v^{-1}D_i$. Let $D'_{ij}$
be a component of $f_1^{-1}D'_i$ dominating $D_{ij}$ and let
$U$ be a small neighborhood in $X_1$
of a generic point $q$ on $D'_{ij}$ above $p$. These neighborhoods
as $j$ varies cover $f_1^{-1}(V)$ outside
a codimension two subset if we include also the open subset
$f_1^{-1}(V\setminus D'_i)$, above which $K_{X_1/Y'}$ and 
$u_1^*K_{X/Y}$ are naturally identical. 
Since $d_*K_{X'/Y'}=K_{X_1/Y}$ on $X_1$,
see the proof of Lemma~13 in \cite{En}, and since $X_1$ is
normal, it is sufficient to show that $u_1^*K_{X/Y}(-u_1^*(df))$ is naturally
identified with $K_{X_1/Y'}$ on $U$, where 
$(df)=\sum f^*(D_i)-f^*(D_i)_{red}\geq f^*D^f$ with equality
attained everywhere on $V$ as given in 
Lemma~\ref{df}. This identification is
easily seen from the fact that, on a neighborhood of $q$, we
can write $u=u_1=u_{ij}\circ u'_{ij}$ where $u_{ij}:X_{ij}\ra X$ is
obtained from the  base change by an $m_{ij}$-cyclic cover $Y_{ij}$ 
of $Y$ branched precisely on $D_i$ and $u'_{ij}$ is obtained
from the (remainder) cyclic cover of $Y'$ over $Y_{ij}$.
Let $f_{ij}:X_{ij}\ra Y_{ij}$ be the map corresponding to $f$.
Then by construction, $K_{X_{ij}}=u_{ij}^*K_X$ and $f_{ij}$ is
smooth on a neighborhood $W$ of $u'_{ij}(q)$. The first implies 
$K_{X_{ij}/Y_{ij}}=u_{ij}^*K_{X/Y}(-(df))$ on $W$ while the
second implies that $K_{X_1/Y'}=(u'_{ij})^*K_{X_{ij}/Y_{ij}}$ 
on $U\cap u_{ij}^{\prime -1}(W)$
as required for $(*)$. Since $K_{X_1/Y'}\hookrightarrow 
u_1^* K_{X/Y}$ in $X_1$ by Lemma~13 of \cite{En} and we
have just shown that this inclusion factors through
$$K_{X_1/Y'}^m=u_1^* (K_{X/Y}(-(df)))^m\hookrightarrow
u_1^*K_{X/Y^\del}^{\otimes m}\ \ \ \ (\dag)$$
outside a codimension
two subset above $V$ where equality
is achieved at every point over $V$, it follows that 
$K_{X_1/Y'}^m\hookrightarrow
u_1^*K_{X/Y^\del}^{\otimes m}(E)$ on $X_1$ 
with equality achieved everywhere on $V$ for some
effective divisor $E$ supported on the exceptional divisor 
of $f$ where $E$ dominates $f^*D^f$. 
Combining this with the theorem of Kawamata-Viehweg
via (P1) and (P2) now gives us theorem~\ref{wp}
and hence theorem~\ref{Cnm}. \BOX\\

By taking $X=Y$ in Theorem~\ref{Cnm}, we get the following 
corollaries.

\begin{theorem}\label{Cnm+} Let $g:Y\ra Z$ 
be a fibration of general type.
Then $$\kappa(Y)= \kappa(Y_z)+ dim(Z),$$ where $Y_z$ is the 
general fiber of $g$. In particular, if 
$Y_z$ is of general type, then $Y$ is also. \BOX
\end{theorem}

\begin{theorem}\label{Sp} A projective 
(or K\"ahler) manifold $X$ is special
if $\kappa(X)=0$. \BOX
\end{theorem}

\begin{rem}
Although we had written up the above theorems at around
the same time in \cite{Ca02} and in \cite{Lu02}
independently, Campana in \cite{Ca01}
had at least claimed similar but slightly 
simpler theorems and therefore they should be credited to him.\\
\end{rem}


The following is the key technical result of this section. We 
say that an orbifold fibration $\fd: \xd\ra Y$
is effective if the horizontal orbifold divisors are,
i.e., $A-(A\cap f)$ is effective for $A=\bd \xd$.

\begin{theorem}\label{wpd} 
Assume that $\fd: \xd\ra Y$ is an effective fibration from a smooth
or e-admissible orbifold. 
Set $K_{\xd/Y^\del}=K_\xd\otimes f^*K_{Y^\del}^{-1}$ 
as $\QQ$-line bundles. 
Let $m$ be a positive integer such that $mD(\fd)$ and $m\bd\xd$ 
are Cartier divisors. Then there is an effective divisor $E$ 
with $E_{red}\subset E(f)$ such that 
$f_{*}(K_{\xd/Y^\del}^{\otimes m}(E))$ 
is weakly positive.
\end{theorem}


\noindent
{\bf Proof:} We reexamine 
the proof of Theorem~\ref{wp}. We observe via $(*)$
and Lemma~\ref{df} that given any vertical orbifold
structure on $X$, meaning that no component of $\bd \xd$
surjects to $Y$, we have the following injection 
outside a codimension two subset of $Y'$
of torsion free coherent sheaves of
the same rank on $Y'$: 
$$
(f')_{*}\big(K_{X'/Y'}^{\otimes m}\big)
\hra v^{*}f_{*}\big(K_{X/Y}(-(df))^{\otimes m}\big)
\hra v^{*}f_{*}\big(K_{\xd/\yd}^{\otimes m}\big),\ \ \ \ (**)
$$
where $\yd$ is the orbifold base of $\xd$. 
Now the key point to obtaining this is the factorization
through the inclusion $(\dag)$. This is unaffected if we
add an admissible (or smooth) horizontal 
(i.e., nonvertical) orbifold divisor 
to $X$ which pulls back only to a horizontal orbifold 
divisor on $X'$ and push forward to the original horizontal
divisor. Hence $(**)$ holds for any admissible $\xd$ after
adding in the horizontal orbifold structures to $X$ and
$X'$ in $(**)$. So
the theorem follows directly from the following combination
of theorems of Kawamta and Viehweg.
\BOX

\begin{prop}\label{horizon}
Assume that $\fd: \xd\ra Y$ is a fibration where $\xd=X\m A$ is
a smooth effective orbifold (with $A$ horizontal). 
Then $f_{*}(K_{\xd/Y}^{\otimes m})$ 
is weakly positive for any $m\in \NN$ such that
$mA$ is Cartier.
\end{prop}

We should remark that very similar theorems that are 
stronger than ours but in more restrictive situations
can be found in \cite{Na0}, \cite{Vi5} and \cite{Ka01}.
Since this is the key point in obtaining the canonical
fibration, we give two proofs of this theorem.
The first is an easier version given to us by Viehweg 
and depends on an easy computation of taking 
integral parts. 
The second is longer and more involved notationally
but more direct conceptually (at least to us) as it reduces
the proof to the same trick
of Viehweg \cite[Corollary~5.2]{Vi} directly by the base 
change property for dualizing sheaves without 
any computation. The second was our original
approach to the problem but we are completely 
indebted to E. Viehweg for the first approach.\\

\noindent
{\bf Proof:} The proof reduces, via the argument of 
\cite[(5.3)]{Vi}, to the following proposition,
which is an orbifold analogue of Corollary~5.2 of 
\cite{Vi}. \BOX

\begin{prop}\label{ho}
Let $\fd: \xd\ra Y$ be a fibration (or a surjective 
morphism) where $\xd=X\m A_0$ is
a smooth effective orbifold (with $A_0$ horizontal)
and $H$ an ample divisor on $Y$ such
that for $m$ as above and some $N\in \NN$, we have
\begin{itemize}
\item[\rm (++)] the sheaf $\hat S^N(f_*(K_{\xd/Y}(H)^{\otimes m}))$
is generated over an open set $U$ by global sections.
\end{itemize}
Then $f_*(K_{\xd/Y}(f^*H)^{\otimes m})(-H)$ is weakly positive.
\end{prop}

We remark that (++) implies that $M^N$ is generated over $f^{-1}(U)$
by global sections where 
$$M = {\rm Image\ of\,}\{f^*f_*(K_{\xd/Y}(f^*H)^{\otimes m})
\ra K_{\xd/Y}(f^*H)^{\otimes m}\},$$
which is actually the condition used in the proof 
(so that it is not necessary to know the exact
meaning of (++)).
Also our proof depends on a strengthening of 
Lemma~5.1 of \cite{Vi} where we replace $T$ by a 
log-manifold with boundary (horizontal with respect
to $g$ and) having no component in common with that of  
$E$ in the notation there. This strengthening is
a {\em direct consequence of Kawamata's Theorem}~32 in 
\cite{Ka} applied to the proof of this lemma (see also the proof of 
Lemma~2 in \cite{Ma}) so we assume that Lemma~5.1 of \cite{Vi} 
has been so strengthened. \\

\noindent
{\bf First Proof of Proposition~\ref{ho} (Viehweg):} \\
Set $K=K_{X/Y}$
and $A_0=\bd\xd=C+D$ where $C$ is reduced and $D$
a $\QQ$-divisor with coefficients 
in $(0,1)$. Let $A=f^*H$
and $L_0:=K+C+D+A=K_{\xd/Y}(f^*H)$ where we
now identify $K$ with a divisor and write
tensor products additively by abuse of notation. Let 
$L=L_{red}=K+C+D_{red}+A$. Then $L$ and $mL_0$ are
Cartier. 

By blowing up $X$, we may assume that
the the base scheme $E$ of $mL_0$ relative
to $f$ is a divisor, or what is the same
that $M$ as defined above is locally free
and $L^m_0=M(E)$. 
This is possible because the situation is compatible with 
blowing ups. In fact, if $\tau: X'\ra X$ is
a good blowing up, then $\tau^*mL_0$ is contained in
the Cartier divisor
$m(K_{X'/Y}+\tau^\circ(C+D)+A)$ by definition of
$\tau^\circ$ (see Definition~\ref{defmo}) and 
they have the same relative global sections,
i.e., direct images with respect to $f\circ \tau$,
as $\tau^\bd$ is an orbifold birational morphism 
(and hence so is its restriction to open subsets).
Note that $\tau^\circ(C+D)$ is reduced in 
$\tau^{-1}(C)$ and not reduced otherwise so that
we again have a decomposition of $\tau^\circ(C+D)$
in the form $C+D$ (reduced plus nonreduced)
as the boundary divisor for the
new $X$ (replacing $X$ by $X'$). By the same token,
we may assume that $E+mC+mD$ has normal crossing support,
where the base scheme $E$ of $mL_0$ relative
to $f$ is a divisor.

Note that the statement of the proposition concerns
solely $f_*(mL_0)=f_*(K_{\xd/Y}(f^*H)^{\otimes m})$, which is
unchanged if we subtract off part of the relative base
divisor $E$ from $mL_0$. More precisely, if we let $mL'$ be
the resulting Cartier divisor after the subtraction, then
we still have $M\hra mL'\hra mL_0$ as line sheaves and so
we see that $f_*(mL')=f_*(mL_0)$ from the fact that 
$f_*M=f_*(mL_0)$  by taking the direct images.
Hence we may assume that $C$ has
no common component with $E$ by the following 
procedure: For a component $D'$ of $E$ lying 
in $C$, we delete this component in $C$ and
if its multiplicity $s$ in $E$ is less than
$m$, we further replace $D$ by $D+sD'/m$. We
can also do the same to ensure that $D$ has
no common component with $E$ but it is not 
needed in the argument. Note that $C$, $D$, $E$, $L_0$ and
$L$ are all changed by this but the point is that
the concluding statement of the proposition is unaffected
by this change.

Now $mL=M(R)$ where $R=E+m(D_{red}-D)$ is an effective 
Cartier divisor with normal crossing support and $M^N$ is
generated by global sections on $f^{-1}(U)$ by hypothesis.
Hence we can apply Lemma~5.1 as strengthened above to 
obtain the weak positivity of $f_*(K+C+L^{(m-1)})$
where 
$$
L^{(m-1)}=(m-1)(K+C+D_{red}+A) - [(m-1)R/m]
$$
and the notation $[G]$ means the Cartier divisor obtained by
taking the integral part of the
coefficients for each component of a $\QQ$-divisor $G$.
Since 
$$ 
E + (m-1)(D_{red}-D) \ge (m-1)R/m \ge (m-1)(D_{red}-D) 
$$ 
and since $-[-(m-1)D]=mD-[D]=mD,$ we have 
$$m(K+C+A)+mD-E-A\le K+C+L^{(m-1)}\le m(K+C+A)+mD-A$$
by taking $[\cdot]$. By the definition of $E$, the
resulting direct images
$$f_*(m(K+C+A)+mD-E-A)\hra f_*( K+C+L^{(m-1)})
\hra f_*(m(K+C+A)+mD-A)$$
are therefore isomorphic, giving the weak
positivity of $f_*(mL_0)H^{-1}$ as required. \BOX\\

\noindent
{\bf Second Proof of Proposition~\ref{ho}:}
Let $\ga:X'\ra X$ be a fibration so that $\ga^*A$
and thus $\ga^*K$ where $K=K_{\xd/Y}$ are restrictions of
Cartier divisors on an open set $V$ of $X'$ such that 
$X\smallsetminus \ga(V)$ is contained in the complement
of a codimension two subset of $X$.
Such a fibration always exists from the problem in
practice but can be constructed in general by the Kawamata 
Galois covering construction of section~4 and twisting the
Galois group by a generically free action on a variety
such as an elliptic curve or an abelian variety (resulting
in an abelian fibration). Let $f'=f\circ \ga$. Let $K'$ be
a Cartier divisor on $X'$ that is identical with
$\ga^*K$ on $V$ and such that $D:=m(K'-\ga^*K)$ is an
effective (and Cartier by construction) divisor. Let $H_0$ be 
an ample divisor on $Y$ with $H=f^* H_0$ and $H'=\ga^*H.$
Note then for any integer $l$ that 
$\ga_*\ga^*(mK+lH)=(mK+lH)$ is an invertible
subsheaf of the torsion free line sheaf $\ga_*(mK'+lH')$
and that these sheaves are identical outside a codimension
two subset of $X$. Hence they give the same result
after applying $f_*$ (for the relevant properties 
such as weak positivity) and so we identify them by
abuse of notation.
As in \cite[Corollary~5.2]{Vi}, we only need to show 
that if  $S^N(f_*(m(K+H)))=S^N(f'_*(m(K'+H')))$ 
is generated by global
sections over an open set $U$ for some $N$, then 
$f_*(mK+(m-1)H)=f'_*(mK'+(m-1)H')$ is weakly positive. The rest
of the argument is the same as that of \cite{Vi} to give
the weak positivity of  $f_*(mK_{\xd/Y})$.

For this purpose,
let $L=K+H$ and $L'=K'+H'$. We can regard $L'$ as
a line bundle on $X'$. Now $L$ is only a line bundle
outside a codimension one subset where $L'=\ga^*L$.
Still, $L$ can be completed
to a singular space $\bar L$ with a vector space structure 
for each fiber over
$X$ (this is clear if the orbifold is standard since
one can use the ``singular'' canonical bundle of the 
corresponding V-manifold but in general the reader 
should have no problem in defining it via $L'$). 

Let $M$ be the image of $f^*(f_*(mL))$ in the invertible
sheaf $mL$. As in Lemma~5.1 of \cite{Vi}, we may
assume that $M$ is spanned over $U$ so that after
suitably blowing up $X$ outside $U$, we may assume
that the effective divisors $mL-M$ and $M$ together
has normal crossing support. We may assume that the 
pull back of these divisors to $X'$  and $D$ are jointly 
normal crossing by blowing up. Then the sum of these
divisors add up to $mL'$ and so the standard construction
of taking $m$-th root apply to obtain a subvariety in 
$L'$ which we denote by $T'$. By construction, $T'$ is
the inverse image of a subvariety $T$ in $\bar L$. Let
$u:Z\ra X$ be the finite map obtained from the normalization
$Z$ of $T$ and $v:Z'\ra X'$ the finite map obtained from
the base change via $u$ with the fibration morphism $g:T'\ra T$.
Let $\te Z$ be a resolution of singularity of $Z$ and 
$\te u:\te Z\ra X$ be the resulting morphism to $X$.
{}From the base change property for dualizing sheaves we get 
$g_*v^!L'=u^!\ga_*L'=u^!(K_{X/Y}+\bar A+H)=K_{Z/Y}+u^*(\bar A+H)$
where $\bar A$ is the reduced part of $A$.
Using the proof of Lemma~5.1 in \cite{Vi}, 
we see that the problem now
reduces to establishing the weak positivity of 
$f_*u_*(K_{Z/Y}+u^*\bar A)$ or, what is equivalent,
of $f_*\te u_*(K_{\te Z/Y}+\te u^*\bar A)$ which we know from
\cite[Theorem~30]{Ka}. \BOX\\

Hence, with notation as in Theorem~\ref{wpd}, we have shown
that $f_*(mK_{\xd/\yd}+E)$ is weakly positive.
So, by the definition of weak positivity and by the fact
that $K_\yd$ is big, given an ample $H$ on $Z$ and $A=f^*H$, 
there exist an $N>0$ with $mNK_{\xd/\yd}+A+NE$ 
and $mNK_\xd -2A$ effective. This shows that
$N(mK_\xd+E)-A$ is effective and so Lemma~\ref{fta} applies
to give
$$
\kappa(K_\xd)=\kappa(mK_\xd+E)=\kappa(mK_\xd|_F)+\dim Z
=\kappa(F^\bd)+\dim Z
$$
where $F^\bd$ is the orbifold general fiber of $f^\bd$
and the first equality follows by assuming $\fd$ is
well adapted, which we can do without loss of generality.
Thus we obtain as before the following corollaries by looking
at the relevant birationally equivalent well-adapted fibrations:

\begin{theorem}\label{Cnm+d} Let $\gd:\yd\ra Z$ 
be a fibration of general type with $\yd$ smooth
or more generally e-admissible effective.
Then $$\kappa(\yd)\geq dim(Z)+\kappa(Y_z^\bd),$$ 
where $Y_z^\bd$ is the 
general fiber of $\gd$.
Hence, Conjecture~\ref{Iit} is true if the $\gd$ or the $\hd$
there is of general type. In particular, if 
$Y_z^\bd$ is of general type, then $\yd$ is also. \BOX\\
\end{theorem}

\begin{theorem}\label{Sp+d} A smooth effective orbifold
or an e-admissible effective orbifold 
$\xd$ is special if $\kappa(\xd)=0$. \BOX
\end{theorem}

We have an analog of a result of N. Nakayama
\cite[Theoremm~7.3.5]{ Na0} on the subadditivity of 
sectional Kodaira dimensions:

\begin{theorem}\label{NN}
Let $\fd:\xd\ra \yd$ be a surjective morphism of
smooth effective orbifolds with connected fibers, 
i.e., a fibration morphism, and let
$F^\bd$ be a general fiber. Then 
$$
\ks(K_\xd-f^*K_\yd+f^*D)\geq \ks(K_{F^\bd}) + \ks(D), 
$$
for any $\QQ$-divisor $D$ of $Y$. In particular,
$$
\ks(K_\xd)\geq \ks(K_{F^\bd}) + \ks(\yd).
$$
\end{theorem}

\noindent
{\bf Proof:} We may assume that $\fd$ is admissible. Lemma~\ref{df} then
allows us to translate Theorem~7.3.3 and its corollary in 
\cite{Na0} into the orbifold context by adding in $f$-exceptional
divisors. Then the proof of Theorem~7.3.5 of \cite{Na0} applies verbatim, 
noting that the $f$-exceptional divisors have no effect 
as $f$ is admissible, to give the result above.
\BOX

We also refer the reader to Definition~6.2.7 in \cite{Na0} 
for the definition of the sectional Kodaira dimension 
$\ks$. For our purpose, it is
sufficient for us to know that $\ks(K_\xd)\geq 0$ if
and only if $K_\xd$ is pseudoeffective and that 
$\kappa(L)\le \ks(L)\le \dim X$ for any $\QQ$-divisor $L$ on $X$.
Finally N. Nakayama has a generalization of the additivity
of Kodaira dimension theorem of Koll\'ar to the compact 
complex category, see \cite[Theorem~5.7]{Na1}. 
In light of our propositions above
we have the following corollary of this theorem,
which we will not use but many of our main results
can also be derived directly from it.

\begin{theorem}[Orbifold Kollar's additivity theorem]\label{kollar}
Conjecture~\ref{Iit} is true if $Y_z^\bd$ or
$(Y_z, f_z^\bd)$ is of general type with $X$ and $Y$
allowed to be compact complex manifolds. \BOX
\end{theorem}

\section{The canonical fibration}

The following theorem resolves some conjectures
of Campana on unramified coverings and on the
structure of his general type reduction in 
\cite{Ca02}. It
also answers in the affirmative a vital question 
raised by Bogomolov concerning Bogomolov sheaves.

\begin{theorem}\label{M}
Let $\xd=X\m A$ be a smooth or an e-admissible effective orbifold.
Suppose that $\xd$ is not special. Then there is one and
only one Bogomolov $\QQ$-sheaf $L\hra \Omega_X^p(A)$ (where
the inclusion is as $\QQ$-sheaves) such that 
$\kappa(L)=\kappa'_+(\xd)=p>0$. That is, there is a 
meromorpic fibration of general type 
$\fd:\xd\dra Y$ which is unique up to
birational modification of $Y$ with the 
property that $\dim Y=\kappa'_+(\xd)$.
Furthermore, given any meromorphic fibration of
general type $\hd:\xd\dra Z $, there exist
a meromorphic map $g:Y\dra Z$ such that
$\hd=g\circ \fd$.
\end{theorem}

\noindent
{\bf Proof:} Suppose $\fd:\xd\dra Y$ and
$\hd:\xd\dra Z$ are
meromorpic fibrations of general type with
$\dim Y=\kappa'_+(\xd)$.
It is sufficient to show that $h$ factors
through $f$. We may assume that $f$ and
$h$ are morphisms by blowing up.
Consider the product map
$\psi=(f, h): X\ra Y\times Z$. Let $\psi:X\ra \bar Y$
be obtained from the Stein factorization of $(f, h)$ 
by blowing up so that it is a morphism with $\bar Y$
smooth. We have natural morphisms
$p:\bar X\ra Y$ and $q:\bar X\ra Z$ which are
obtained from the projections to the
two factors of $X\times Y$. It is easily 
seen that these are fibrations since $f$ and
$h$ are. It follows that the general fiber 
$F=\psi^{-1}(y,z)= f^{-1}(y)\cap h^{-1}(z)$
is well defined and smooth if nonempty. 
Let $d_c=\dim X-\dim Y$.
If $\dim F\neq d_c$, then $\dim \bar Y> \dim Y$
and so $p$ is a fibrations with positive dimensional
general fibers. The general orbifold 
fibers of this fibration
must be of general type by Corollary~\ref{tech}
and so we have that $\bar Y^\bd$ must be of general
type by Theorem~\ref{Cnm+d}. 
This means that $\psi$ itself is of general 
type as it is admissible. This is a contradiction 
to the definition of $\kappa'_+(\xd)$. Hence 
$\dim F= d_c$, which means that 
$f^{-1}(y)\cap h^{-1}(z)=f^{-1}(y)$ if the
intersection is nonempty. This, in turn, means
that $h$ factors through $f$
\BOX\\

It follows from this theorem that there is no ambiguity
in calling the $\fd$ above as the top (or the maximum)
fibration of general type from $X$ (up to a modification 
on the base), which we denote by $bc_\xd^{}$. The following
was a conjecture of Campana in the non-orbifold case.

\begin{cor}\label{cover} Let $\xd=X\m A$ be an 
e-admissible effective orbifold.
Let $\pi: \bar X\ra X$ be a finite \'etale
covering and $\bar X^\bd= \bar X\m \pi^*A$. Then 
$$bc_{\bar X^\bd}^{}=bc_\xd^{}\circ \pi\ \ \ \ \ 
\mbox{and so}\ \ \ \ \ 
\kappa'_+(\bar X^\bd)=\kappa'_+(\xd).$$
In particular, $\bar X^\bd$ is special if
$\xd$ is special.
\end{cor}

\noindent
{\bf Proof:} Clearly we may assume that $\pi$ is
a Galois covering. Let $\bar f=bc_{\bar X^\bd}^{}$. 
By uniqueness, $\bar f$ 
commutes with the covering group $G$ 
and hence gives rise to a morphism $f:X\ra Y$ and
$\fd$ is of general type by construction. Since 
$\pi^\bd$ composed with any general type fibration 
$f_1^\bd$ from $\xd$ is of general type, Theorem~\ref{M} 
says that $f_1^\bd$ factors through $\fd$ so that
$\fd=bc_\xd^{}$.

An easier argument is to look at the problem on the
level of the top Bogomolov $\QQ$-sheaf $\bar L$ of
$\bar X$ which is invariant under $G$ by its uniqueness. 
Hence, this sheaf descends to a Bogomolov $\QQ$ sheaf
by the fact that we can take trace of sections, i.e., the
Kodaira dimension of $\bar L$ is unchanged after 
taking its quotient by $G$ by the same argument
as that for Lemma~\ref{A'}. \BOX\\

The next theorem is the main theorem of this paper.

\begin{theorem}\label{Main} 
Let $X^\bd$ be an e-admissible effective orbifold. Then,
up to modifications on the base, there is a unique 
(almost holomorphic) fibration of general type 
$$\hd={b \,\! c}_{\!\xd}^{}: X^\bd\dra Z=bc(\xd)$$
whose general fibers are 
special. We have $\kappa'_+(\xd)=\dim bc(\xd)$.
Moreover, the general fibers of $\hd$ 
are the maximal special subvarieties of $\xd$
through the general points of $X$.  
\end{theorem}

\noindent
{\bf Proof:} We have already seen in the previous 
theorem that there is a maximum fibration of
general type $bc_\xd^{}$ from $\xd$. It remains 
to show that this fibration is special since the
other claims follows from Proposition~\ref{spgt}
(the last statement follows from the same and
Proposition~\ref{spd} using 
a component of the Barlet space of $bc(\xd)$
whose general points corresponds to images
of maximal special subvarieties of $\xd$
through some general points of $X$).

We will assume that $\hd= bc_\xd^{}:\xd\ra Z$ is a
well adapted fibration by blowing up. Let
$bc_\hd^{}:\xd\ra Y$ 
be the relative basic
reduction (Proposition~\ref{GTR})
of  $\hd=bc_\xd^{}$, which we
again assume to be admissible by blowing up.
We have a canonical projection $p:Y\ra Z$.
If $p$ is not birational, then the general
fiber $Y_z$ of $p$ has positive dimension and
so we have again by Theorem~\ref{Cnm+d} 
that $\yd$ is of general type. But this
contradicts the definition of $bc_\xd^{}$.
Hence $p$ is birational, which means that
the general orbifold fibers of  $bc_\xd^{}$
are special by Proposition~\ref{GTR}.
\BOX\\

This resolves in particular a conjecture of Campana 
in \cite{Ca02}.

\begin{cor}  Let $X\in {\mathcal C}$.
The maximal special fibration $c_X$ defined in  
\cite{Ca02} is the same as 
${b \,\! c}_{\! X}^{}$ up to birational equivalence
and hence is of general type. \BOX
\end{cor}

Another easy corollary is the following theorem, which
gives us in principle many possible special fibrations
connected to a compact complex orbifold.

\begin{theorem}\label{M'}
Given a fibration $h:X\ra Z$ with either $h\in \mathcal C$
or $Z$ Moishezon. If $\xd$ is an e-admissible effective
orbifold, in particular if $\xd$ is smooth and effective, then
$bc_\hd$ is a special fibration. \BOX
\end{theorem}

\section{Remarks on the Albanese, Abundance and Cnm}

There has been quite a number of focus on the 
structure of the Albanese map of a projective or
K\"ahler manifold with nef anticanonical bundle, 
see for example \cite{DPS}. The best results so
far in the projective category are that of Qi Zhang in \cite{ZQ}
and a slight generalization by Demailly-Peternell-Schneider
in \cite{DPS}. These results are essentially reformulations
of a result of Miyaoka in \cite{Miya2} that requires positive
characteristic arguments, see also \cite{Miya0}. 
Although it is not difficult to deduce from their proofs
the following result in the case $X_0$ has no orbifold 
structure (from, for example, the proof of Claim~20 in \cite{DPS}
following the constructions given in the proof below), 
we give an alternative proof below without involving
positive characteristic and for the full orbifold case.
In particular, we solve completely the conjectures of
\cite{DPS} and of \cite{Ca02} concerning the Albanese
map of a projective manifold with nef anticanonical
bundle, and give a generalization to the quasiprojective
case of these results including that of \cite{ZQ}.

\begin{theorem}\label{-K}
Let $X_0^\bd$ be a smooth projective and effective orbifold 
such that $K_{X_0^\bd}^\vee$ is a nef $\QQ$-line bundle or 
equivalently $-K_{X_0^\bd}$ is a nef $\QQ$-divisor.
Suppose $f^\bd_0:X^\bd_0\dra Y_0$
is a meromorphic fibration of general type
with $Y_0$ smooth projective.
Then $\kappa(Y_0, \fd)\le 0$. 
In particular, $X_0^\bd$ is special. 
\end{theorem}

\noindent
{\bf Proof:} By Proposition~\ref{almo},
$f_0$ is almost holomorphic.
By flattening and then resolving using the
$-\infty$ multiplicity on the exceptional
fiber, we obtain by Lemma~\ref{none} 
a well adapted birationally equivalent
fibration $\fd:\xd\ra Y$ via birational
morphisms $v^\bd:\xd\ra X^\bd_0$ and 
$u:Y\ra Y_0$. Then we may write
$$K_\xd = K_X+ A = K + E$$
where $A=\bd\xd$ has normal crossing support, 
$K=v^*K_{X_0^\bd}$ and $E$ 
is an effective divisor 
supported on the exceptional divisor of $v$. 
We note that $-K$ is nef 
and $E$  maps to a codimension two
subset $S$ of $X_0$. 
Assume that $\kappa(Y_0, \fd)> 0$. 
By blowing up $Y$ further if 
necessary (and hence $X$ as well), 
we may assume that the Iitaka fibration
$g$ of $\yd:=Y\m D(\fd)$ is a morphism from $Y$ onto a 
projective manifold $Z$ of dimension $\kappa(Y_0, \fd)> 0$.
Let $H_0$ be an ample $\QQ$-divisor on $Z$ and $H=g^*H_0$.
It follows from the standard lemma on Iitaka fibrations
that we may choose $H_0$ such that 
$K_\yd-2H$ is linearly equivalent to an
effective $\QQ$ divisor, which we denote by
$2H\ \lesssim\ K_\yd$.
Let $A'$ be an ample $\QQ$-divisor on $X$. 
Now, there exist an $\eps>0$ such that 
$$E-f^*H+ \eps A'$$
is not pseudoeffective (To see this, restrict to a 
general complete
intersection curve $C$ on $Y_0$ where $f_0$ and 
$g\circ u^{-1}$ becomes 
holomorphic and then to a general complete 
intersection curve on $X_0$ above $C$ avoiding $S$). 
As $-K+\eps A'$ is ample, it is linear equivalent
to an effective $\QQ$-divisor $D$ such that 
$D\leq D_{red}$ and that $D_{red}\cup A_{red}$
has normal crossing support. It follows that 
$X^{\bd_1}=X\m (A+ D)$ is a smooth orbifold
such that $f^{\bd_1}:X^{\bd_1}\ra Y$ is 
still a fibration of general type dominating 
$\yd=Y\m D(\fd)$.  Since 
$K_{X^{\bd_1}}$ is linearly equivalent to $E+\eps A'$,
which restricts to a big
divisor on the general fibers $F$ of $f$ 
by construction, we have
$$
\ks(K_{X^{\bd_1}} - f^*K_\yd + f^*H)\geq 
\dim F + \ks(H)\ge \kappa(H) = \kappa(Y_0, \fd)> 0
$$
by Theorem~\ref{NN}. 
Hence $K_{X^{\bd_1}} - f^*K_\yd + f^*H$
is pseudoeffective and linearly equivalent to
$$
E - f^*K_\yd + f^*H + \eps A'\ 
\lesssim\ E-f^*H+ \eps A',
$$
forcing the right hand side to be 
peudoeffective also. This is a contradiction. \BOX\\

For application to the Albanese map, we need the following.

\begin{lem} 
Let $\xd=X\m A$ be an orbifold $\fd:\xd\ra Y$ a fibration
with orbifold base $\yd$. Let $K=K_\yd=K_Y+B$ where $B=D(\fd)$.
Let $E=f^*K-L_\fd\cap f^*K$, which is an effective 
$\QQ$-divisor supported on $E(f)$. Suppose there is
an $\eps>0$ such that $L_\fd(-\eps E)$ is linearly 
equivalent to an effective $\QQ$-divisor. Then
$$ \kappa(Y,\fd)=\kappa(\yd).$$
This holds for $A+(1-\eps) E(f)\ge A\cap0$ and
$\kappa(K_Y(B\cap0))\ge 0$ (e.g., $A\ge 0$ 
and $\kappa(K_Y)\ge 0$).
\end{lem}

\noindent
{\bf Proof:} We know 
$$\kappa(Y,\fd)=\kappa(L_\fd)\le \kappa(\yd):=\kappa(K)$$
by the Hartog extension theorem (see Lemma~\ref{none}).  
But the reverse inequality follows directly from our hypothesis,
which says that $f^*K-L_\fd\le E \lesssim \frac{1}{\eps}L_\fd$.
Since $L_f-f^*K_Y\ge E(f)$,
the last statement of the proposition guarantees that 
$\eps E(f)\le L_\fd-f^*K_Y(B\cap0)$, which is sufficient
for our purpose because $E\le E(f)$ and $K_Y(B\cap0)$ is
linearly equivalent to an effective $\QQ$-divisor. \BOX\\

The following is a simple consequence of the 
structure theorem of Kawamata \cite{Ka} (and Ueno \cite{Ue})
for subvarieties and finite covers of semi-abelian varieties
(used in \cite{DL} for example), see
\cite{Ca02} for the non-log case.  We refer the
reader to the original source \cite{Ka}
(see also \cite{KV} of Kawamata-Viehweg) for the details.

\begin{theorem}\label{Alb} 
Let $\xd=(X,D)$ be an algebraic log-manifold
with $V=X\smallsetminus D$. Suppose $\xd$ is special. 
Then the quasi-Albanese map of $(X,D)$, 
$$
\alpha=Alb_{(X,D)}: V\ra Alb(V),
$$ 
is an algebraic fibration (i.e., a dominant
algebraic morphism with connected general fibers) 
with $D(\alpha)=0$ (i.e., $\alpha$ imposes
no branched orbifold or log structure on its base
$Alb(V),$ which is a semi-abelian variety). 
In particular, $\alpha$ surjects to the complement of
a set of codimension two or higher in $Alb(V)$.
Hence, the log-irregularity $\bar q$ satisfies
$$\bar q(V):=\dim H^0(\Omega_X(\log D))
=\dim Alb(V)\le \dim X$$ and equality
implies that $\alpha$ is a birational morphism. \BOX
\end{theorem}

So the above theorem is a straightforward
but broad generalization of the structure theorem of 
Kawamata and Viehweg concerning varieties
with Kodaira dimension zero, which are well connected to 
the terminology of being special in the literature.\\

We now note the following result of Mihai Paun
in \cite{P} derived using a strong differential 
geometric result of Cheeger and Colding (\cite{CC}),
a result of Qi Zhang (\cite{ZQ}) concerning the
Albanese map that employs characteristic $p>0$ techniques
from \cite{MM1} and a result of Campana (\cite{Ca2}). 
Our result refines that of Qi Zhang
but without involving characteristic $p>0$ 
thereby giving another proof (more direct for
differential geometry) of this result of Mihai Paun,
which we can now state
with a slight strengthening at least in
the direction of the Albanese map. See \cite{P}
for the proof. The last part of the theorem
follows from the standard exact sequence of
fundamental groups for a fibration, see for
example \cite{XG}, whose theme in relation to
the Albanese map in the full
orbifold context will be taken up later.

\begin{theorem}[M. Paun]
Let $X$ be a projective manifold with  
nef $-K_X$. Then $X$ has almost abelian
fundamental group.
Furthermore, the Albanese map of $X$
is surjective, has connected general fibers
$F$ which are special, imposes no orbifold
structure on the base and induces $\pi_1(X)$ as 
an extension of $\pi_1(Alb(X))$ by a
quotient of $\pi_1(F)$, which is necessarily
finite. \BOX
\end{theorem}

It is of interest to generalize this theorem to the
log case as there are already intensive studies of
the problem for log-surfaces, by D.-Q. Zhang for
example (see \cite{FKL} and the references therein).
For now, we can state the following result from 
the analog of the above result for a standard orbifold;
the general case is treated in \cite{Lu03}.
We refer the reader to \cite{Lu01} for the standard
definition of the fundamental group of a standard
smooth orbifold. Recall that $ X^{|\bd|}:= X\m A_{red}$ 
for a smooth orbifold $\xd=X\m A$. 

\begin{theorem}\label{pi}
Let $\xd=X\m A$ be a smooth projective
standard orbifold. Suppose that $ X^{|\bd|}$
is special and that there is
an effective divisor $A'$ with $A'_{red}=[A]$ such that
$-K_\xd+\eps A'$ is nef for all sufficiently small
$\eps>0$. This holds for example when $-K_{X^{|\bd|}}$, 
$A_{red}-A$  and $A'$ are nef and, in the case $A_{red}=A$,
when $-K_\xd$ and $A'$ are nef.
Then the orbifold fundamental group $G$ of 
$\xd$ is almost abelian.
\end{theorem}

\noindent
{\bf Proof:} We assume $A_{red}=A$ for simplicity
since the general case follows the same proof.
Since $-K_\xd+\eps A'$ is nef for all sufficiently small
$\eps>0$, $-K_\xd$ is nef and Theorem~\ref{Alb}
yields that the quasi-Albanese map $\alpha$ of $\xd$
surjects onto the complement of a 
codimension two subset of $J=Alb(X\smallsetminus A)$
and has connected general fibers. Hence, we have an
exact sequence: 
$$
1\ra K\ra G \ra \pi_1(J)\ra 1,
$$
where $K$ is the fundamental group of the fiber product
of the universal covering map of the semiabelian variety
$J$ with $\alpha$. Since every element of $K$ has finite
order by the definition of $\alpha$, whose rank is the
rank of abelianization of $G$, it is sufficient to 
establish that $K$ is almost abelian. For this purpose,
we choose a small loop around each component of
$A$ and consider them as elements of the the fundamental
group of $X\smallsetminus A$ for which a base point is
fixed. We let $l_1,\dots,l_k$ be all the elements of $K$ 
that are given by these loops.

Let $m$ be a positive multiple of the coefficient
of the primary decomposition of $A'$. By definition of $A'$,
$\bar A=A-A'/m$ is a standard orbifold structure. We consider
$\bar \xd=X\m \bar A$ and its orbifold fundamental group 
$\pi_1(\bar\xd)$. Using the nefness of 
$-K_{\bar \xd}$, the proof
of \cite{P} still works via an orbifold version of
the theorem of Cheeger and Colding (which is 
a standard corollary of their theorem for which we refer
to \cite{Lu03} for the details) to give the almost
abelianness of $\pi_1(\bar\xd)$ via the quasi-Albanese map
above. It follows by the definition of 
$\pi_1(\bar\xd)$ that the quotient group $K/N_m$ 
is almost abelian where $N_m$ is the normal subgroup of $K$ 
generated by $l_i^{m_i}$ with $m_i$ is the multiplicity of
the component of $\bar A$ corresponding to $l_i$. Since
the $l_i$'s in $K$ have finite order, we can arrange
to have $N_m=\{1\}$ by multiplying these orders into $m$.
Hence $K$ is almost abelian, as required. \BOX\\

One can deduce easily from the above considerations
a conjecture in \cite{KZ}
that the smooth locus of a Fano variety
whose singularities are log terminal has finite 
fundamental group, see \cite{Lu03}. This gives an alternative proof
of a result of S. Takayama \cite{Ta} that such Fano varieties
are simply connected.\\[1mm]

We now make some remarks on the weak abundance conjecture, which 
for an effective smooth orbifold $\xd$ states  
the following:
(+) If $K_\xd$ is pseudoeffective, then it is $\QQ$-effective.\\

It is well known that (+) is equivalent to Mori's orbifold
ansatz, Conjecture~\ref{Mo}, for effective orbifolds modulo 
the flip conjecture  (i.e., the existence and
termination of log-flips which nowadays is embedded in a
linear algebra program of Shokurov). Our main theorem 
effectively reduces (+) to the same for special $\xd$
via Theorem~\ref{Cnm+d}. We note that (+) is known 
by Theorem~(7.4.1) in
the case of vanishing numerical Kodaira dimension, or 
equivalently $\ks(K_\xd)=0$ (see the remark before 
Theorem~(7.3.5) in the introduction of \cite{Na1}),
so that $\kappa(\xd)=0$ in this case, which is the
classical case of special orbifolds.

We now assume the subadditivity conjecture 
(!) that nonnegative fiber and
base Kodaira dimensions implies the same for the total space
as in our proof of Remark~\ref{rm}.
Now assume (+) for dimension 
less than that of $\xd$. To show (+) for $\xd$, we may thus assume
that $\xd$ is special. One deduces directly from (!) the existence
of a mermorphic fibration $\fd:\xd\dra Y$ with 
$\kappa(Y, \fd)\ge 0$ dominating all other such fibrations.
We may assume $\fd$ is a morphism by blowing up as
$K_\xd$ remains pseudoeffective.
If $\dim Y>0$, then (+) holds for the general fibers all of which
have pseudoeffective orbifold canonical bundles so that (+) holds
also for $\xd$ by (!). If $\dim Y=0$ or equivalently if
$\kappa'(\xd)=-\infty$, then (+) would follow from 
Conjecture~\ref{rc}. If we call a smooth orbifold $\xd$ with
$\kappa'(\xd)=-\infty$ {\bf very special}, we can summarize 
our observation here as follows.

\begin{rem}\label{rm10} 
Assume the part of the subadditivity conjecture given by 
{\rm (!)} above.
Then the weak abundance conjecture {\rm (+)} above for arbitrary
smooth orbifold $\xd$ is equivalent to the same for
very special $\xd$ and follows from the conjecture 
that very special orbifolds have
Mori's fibrations which in turn follows from
Conjecture~\ref{rc} or the orbifold version of
the Conjecture of Kollar
stated just before Remark~\ref{rm2}.
\end{rem}

As for (!), it can also be reduced to the case of very
special varieties analogous to (6.2.1)$_r$ in \cite{Mori}.
We only demonstrate this by reducing orbifold Cnm, namely
Conjecture~\ref{Iit}, to the case of special varieties:

\begin{rem}\label{rm11}
Assume Conjecture~\ref{Iit} for special fibrations whose
total space is special. Then Conjecture~\ref{Iit} 
holds in general.
\end{rem}

\noindent
{\bf Proof:} With the said assumption, we first establish
Cnm for a special fibration $f: X\ra Y$. We suppress boundary
divisors for simplicity of notation but the proof works in 
general. Consider the basic fibration $h: X\dra Z$ of $X$. 
As $f$ is special, there is a fibration $g: Y\dra Z$ (almost
holomorphic) such that $h=g\circ f$. Theorem~\ref{Cnm+d}
says that $\kappa(X) \ge \kappa(X_z) +\dim Z$. Since the
general fiber $X_z$ of $h$ is special, our assumption says
that $\kappa(X_z)\ge \kappa(X_y)+ \kappa(Y_z)$. These coupled
with $\kappa(Y)\le \kappa(Y_z)+\dim Z$ from the easy addition
law gives us Cnm for $f$. 

Now drop the assumption that $f$
is special. We have by Proposition~\ref{GTR} a special 
almost holomorphic fibration $e: X\dra W$ and a fibration
$W\dra Y$ whose general fibers $W_y$ are of general type.
Hence, Kollar's additivity (Theorem~\ref{kollar}) gives
$\kappa(W)\ge \dim W_y +\kappa(Y)$. Also, since $e$ is
special, we have established above that 
$\kappa(X)\ge \kappa(X_w)+\kappa(W)$. So our result
follows from $\kappa(X_y)\le \kappa(X_w)+\dim W_y$
by the easy addition law. \BOX\\

\section{Appendix A. Orbifold morphism}

The most natural definition of an orbifold morphism is via
orbifold metrics, which we did in \cite{BL2} and in
\cite{Lu01} under the name of Z-orbifold maps. 
We only give this definition for standard
smooth orbifolds in the following but it should be clear how it
is generalized for an arbitrary smooth orbifold.

Given an $n$-dimensional complex manifold $Y$
with a divisor $D^Y=\sum_i (1-1/m_i)D_i$ 
where $D^Y_{red}=\sum_i D_i$
is a simple normal crossing divisor and $m_i\in \NN$. 
Then $D^Y$ defines an orbifold
structure on $Y$ and $Y$ so endowed is written as
$Y^\bd:=(Y;  D^Y)=Y\m D^Y.$ Locally, $Y^\bd$ can be
realized as a quotient of open subsets of $\CC^n$ by a finite
subgroup of $\CC^{*n}$ of cyclic indices prescribed by the
corresponding $m_i$'s. Pushing down the flat metric on $\CC^n$
by these local quotient maps and patching them up via a partition of
unity on $Y$ gives us an orbifold metric $g$ for $Y^\bd$ and we
will define any other metric on 
$Y\smallsetminus D^Y_{red}$ to be an orbifold
metric for $Y^\bd$ if its ratio with $g$ extends to a continuous
function on $Y$.
\begin{definition}
A holomorphic map $f: Y_1\ra Y_2$ between the
complex manifolds $Y_1$ and $Y_2$ endowed with orbifold structures
as above is called a Z-orbifold map if the natural
pull back of an orbifold metric on $Y_2^{orb}$ to $Y_1$ can be
dominated by an orbifold metric on $Y_1^{orb}$.
\end{definition}

It is easy to see
that these definitions are independent of the choices involved. 
It is not difficult to verify that orbifold morphism between smooth
standard orbifolds are precisely the Z-orbifold maps.
\\


For the next definition, if 
$B=\sum_{i\in I} a_iD_i$ is a $\QQ$-divisor, we
call $B'=\sum_{i\in I'} a_iD_i$ 
a partial sum in $B$ for any $I'\subset I$.

\begin{defn} \label{d1}
Given a surjective morphism $f:X\ra Y$ of manifolds $X, Y$
endowed with smooth orbifold structures $A, B$ respectively.
We say that $f$ gives rise to a {\bf good K-map} 
$\fd:X^\bd\dra \yd$ if for $B^*=B\cap \Delta(\fd)$ and 
$\bar B:=f^{-1}(B^*_{red}),$ we have 
$$f^{-1}(B_{red})\subset A_{red}\ \ \ \mbox{  and }\ \ \ 
f^*(B-B^*_{red})\leq A\cap\bar B - \bar B$$
as $\QQ$-divisors. If further
$$f^*(B'-{B'}_{red})\leq A\cap\bar B' - \bar B'$$
for $\bar B':=f^{-1}({B'}_{red})$ for all 
partial sum $B'$ in $B^*$, then we say that
$\fd$ is a {\bf  good orbifold morphism} or simply
a good morphism for short.
\end{defn}

Given log-manifolds $X\m A$ and $Y\m B$, we observe that a good
morphism $\fd: X^\bd\ra \yd$ is precisely a surjective 
log-morphism of
log-manifolds. The following is a straightforward generalization
of generically finite surjective K-maps defined in 
Definition~\ref{dk}.

\begin{defn} \label{d2}
Given a surjective morphism $f:X\ra Y$ of manifolds $X, Y$
endowed with smooth orbifold structures $A, B$ respectively.
We say that 
$\fd:X^\bd\dra \yd$ is a {\bf K-map}
if $$f^*K_Y+B\leq L_f+A\cap f.$$
\end{defn}

\begin{lem}\label{K-map}
If $\fd:X^\bd\dra \yd$ between smooth orbifolds 
$X^\bd$ and $\yd$ is a good K-map, then $\fd$ is a 
K-map. If further $\fd$ is a good morphism, then it is an
orbifold morphism. A (good) K-map $f^\bd$ 
is a (good) morphism outside its
exceptional divisor $E(f)$. \BOX
\end{lem}

\noindent
{\bf Proof:} 
The first two statements are direct consequences of the fact that
$f^*(K_Y+B_{red})\leq L_f+\bar B$ because 
$f^{|\bd|}: Y\m f^{-1}(B_{red})\ra Y\m B_{red}$ is
a log-morphism. The second statement is clear from 
the local expression for $\fd$. \BOX\\

We have many different types of 
maps. We conclude with their easily observed behaviors
under composition.

\begin{lem} Let $f:X\ra Y$ and $g:Y\ra Z$ be morphisms,
$X^\bd, \yd$ and $Z^\bd$ orbifolds and $h=g\circ f$.
If $\fd$ and $\gd$ are both (good) K-maps and $g$
is generically finite and surjective or if they are both
(good) morphisms, then so is $h^\bd$. \BOX
\end{lem}


\medskip
\noindent
Steven Shin-Yi Lu\\
Max Planck Institute of Mathematics\\
7 Vivatsgasse, 53111 Bonn\\
Germany\\

\noindent
slu@math.utoronto.ca\\
slu@mpim-bonn.de

\end{document}